\documentclass[11pt]{amsart}
\usepackage[latin1]{inputenc}
\usepackage[T1]{fontenc}
\usepackage{latexsym}
\usepackage{xspace}
\usepackage{amstext}
\usepackage{amsfonts}
\usepackage{amssymb}
\usepackage{amsmath}
\usepackage{amscd}
\usepackage[mathcal]{euscript}
\usepackage{fancyhdr}
\usepackage{geometry}
\usepackage{tikz}
\tikzset{node distance=1.5cm, auto}

\usepackage[all]{xy}
\usepackage{enumerate}

\newtheorem{prop}{Proposition}[section]

\newtheorem{thm}[prop]{Theorem}
\newtheorem{lemma}[prop]{Lemma}

\theoremstyle{definition}
\newtheorem{defi}[prop]{Definition}
\newtheorem{exmp}[prop]{Example}

\newtheorem{rmk}[prop]{Remark}

\newcommand {\dem}{\textbf{Proof:} }

\newcommand{\prove}[1]{\textbf{Proof of {#1}:} }

\newcommand {\ov}{\overline}

\newcommand {\kg}{(k G)^*}

\newcommand {\ri}{\rightarrow}

\newcommand*{\tens}{\otimes}

\newcommand {\dl}{\Delta_l}

\newcommand {\vare}{\varepsilon}
\newcommand {\wide}{\widetilde}
\newcommand {\prodt}[2]{\underset{#1}{\stackrel{#2}{\prod}}} 
\newcommand {\hte}{H^{\tens n}}
\newcommand {\coe}{h^1\tens \dots\tens h^n}
\newcommand {\coen}{h^1 , \dots , h^n}

\newcommand {\prodn}[3]{ ({#1}^{1}_{(#2)}, \dots , {#1}^{#3}_{(#2)})}

\newcommand {\inv}[1]{{#1}^{-1}}
\newcommand {\omg}[2]{\omega(p_{#1},p_{#2})}

\newcommand {\func}[5]{\begin{array}{cccl} #1 & #2 & \longrightarrow & #3\\ & #4 & \mapsto &
#5\end{array}}
\newcommand {\sub}[2]{#1_{(#2)}}
\newcommand {\subp}[2]{#1^{#2}}
\newcommand {\rmr}[1]{\rm (#1)}
\newcommand {\find}{\hfill $\blacksquare$ }

\newcommand {\prodah}{\underline{\wide{A}\#_{\omega}H}}

\begin{document}

\title[Cohomology for partial actions of Hopf algebras]{Cohomology for partial actions of Hopf algebras}

\author[E. Batista]{Eliezer Batista}
\address{Departamento de Matem\'atica, Universidade Federal de Santa Catarina, Brazil}
\email{eliezer1968@gmail.com}
\author[A.D.M. Mortari]{Alda D.M. Mortari}
\address{Departamento de Matem\'atica, Universidade Federal de Santa Catarina, Brazil}
\email{aldadayana@gmail.com}
\author[M.M. Teixeira]{Mateus M. Teixeira}
\address{Istituto Federal de Santa Catarina, Brazil}
\email{mateusteixeira.mtm@gmail.com}
\thanks{\\ {\bf 2000 Mathematics Subject Classification}: Primary 16W30; Secondary 16S40, 16S35, 58E40.\\   {\bf Key words and phrases:} 
partial Hopf actions, cohomology for partial Hopf actions, cleft extensions. }

\flushbottom

\begin{abstract} In this work, the cohomology theory for partial actions of co-co\-mmu\-ta\-ti\-ve Hopf algebras over commutative algebras
is formulated. This theory ge\-ne\-ralizes the cohomology theory for Hopf algebras introduced by Sweedler and the cohomology theory 
for partial group actions, introduced by Dokuchaev and Khrypchenko. Some nontrivial examples, not coming from groups are 
constructed. Given a partial action of a co-commutative Hopf algebra $H$ over a commutative algebra $A$, we prove that there exists a new 
Hopf algebra $\widetilde{A}$, over a commutative ring $E(A)$, upon which $H$ still acts partially and which gives rise to the 
same cohomologies as the original algebra $A$. We also study the partially cleft extensions of commutative algebras by 
partial actions of cocommutative Hopf algebras and prove that these partially cleft extensions can be viewed as a cleft 
extensions by Hopf algebroids.
\end{abstract}

\maketitle

\section{Introduction}

The history of Hopf algebras began within the context of algebraic topology with the seminal paper by Heinz Hopf, published 
in 1941, describing the algebraic properties of the cohomology ring of a group manifold \cite{Hopf}. The subject of group 
cohomology soon became increasingly more independent of its topological background assuming a more algebraic formulation 
\cite{Ade,W}. The first formulation of a cohomology theory of co-commutative Hopf algebras acting over commutative algebras 
was done by Moss E. Sweedler in 1968 \cite{Swe}, which, in certain sense, became paradigmatic for further developments in 
this area.

Partial group actions, in its turn, had its beginning with the work of Ruy Exel in the classification of certain class of
C*-algebras which had an action of the unit circle but which cannot be described as a usual crossed product \cite{Ex1}. 
A more algebraic formulation for partial actions was done by Mikhailo Dokuchaev and Ruy Exel in \cite{DE} and then, partial
actions drew the attention of algebraists and allowed further developments in several directions. One of the developments 
particularly relevant for our discussion here is the notion of twisted partial actions of groups \cite{DES1} and its 
globalization \cite{DES2}. There, one can see the definition of partial $2$-cocycles in order to define twisted actions 
and partial crossed products, this suggested the existence of a general cohomology theory in which these partial $2$-cocycles 
could be placed. This cohomological theory for partial group actions was achieved by Mikhailo Dokuchaev and Mykola Khrypchenko 
in \cite{DK1}. This theory is constructed upon partial actions of groups over commutative algebras. As the very notion of a 
partial action of a group is deeply related with actions of inverse semigroups \cite{Ex2}, the same authors, \cite{DK2}, 
could place their cohomology theory for partial group actions within the context of cohomology of inverse semigroups, developed 
by Hans Lausch \cite{Lau}.

Partial actions came into the Hopf algebra context by the work of Stefaan Caenepeel and Kris Janssen in \cite{CJ}. 
This work allowed the generalization of classical results in Hopf algebra theory as well several ideas developed for 
partial group actions, as glo\-ba\-li\-zation theorem \cite{AB}, Morita equivalence between the partial smash product and the 
invariant subalgebra \cite{AB2}, duality for partial actions \cite{AB3}, partial representations \cite{ABV}, etc. For a more 
detailed account on recent developments of partial actions of groups and Hopf algebras, see \cite{Bat} and references therein.
Of particular interest for the present article are the notions of twisted partial actions of Hopf algebras, partial crossed 
products and partially cleft extensions of algebras by Hopf algebras \cite{ABDP1}. The aim of this paper is exactly to f
ormulate a cohomology theory for partial actions of Hopf algebras, in the same spirit of \cite{DK1}, such that the partial
$2$-cocycles defined in \cite{ABDP1} can be placed properly. This cohomological theory is obtained for the case of partial 
actions of a co-commutative Hopf algebra $H$ acting partially over a commutative algebra $A$. Moreover, one can, without 
loss of generality, replace the original algebra $A$ by a commutative and co-commutative Hopf algebra $\widetilde{A}$ 
over a base algebra $E(A)\subseteq A$ and still obtain the same cohomology theory.

As we have already learned in \cite{ABV}, the theory of partial actions of Hopf algebras in fact is deeply related to the
theory of representations of Hopf algebroids. This opens a totally new landscape to be explored, for example, in this work
we prove that partially cleft extensions can be understood as cleft extensions by Hopf algebroids in the sense of Gabriella 
B\"{o}hm and Tomasz Brzezinski \cite{BB} and then one can raise new questions on how to put this cohomological theory for 
partial actions in the context of cohomology for Hopf algebroids \cite{BS1,KoP}. 

This article is organized in the following way. In Section 2, we give some mathematical preliminaries, recalling the notion
of a partial action of a Hopf algebra over a unital algebra and giving some examples of such partial actions. Special attention
is required for examples \ref{groupactionsoverfield} and \ref{groupgradingsoverfield}, which will serve as basis for our
specific examples of cohomologies given in Subsection 3.4. Section 3 is dedicated to the construction of our cohomological 
theory for partial actions of a co-commutative Hopf algebra $H$ over a commutative algebra $A$. We start with the study of a
system of idempotents in the commutative convolution algebras $\mbox{Hom}_k (H^{\otimes n},A)$, for a natural $n$. In subsection
3.2 we define the cochain complex, $C^{\bullet}_{par}(H,A)$ associated for the partial action of $H$ upon $A$. The ideas 
involved follows the classical construction due to M.E. Sweedler \cite{Swe} but in order to overcome the complexities coming
from partial actions we define some auxiliary operators which help us to prove that the coboudary operator is a morphism 
of abelian groups (Theorem \ref{deltaismorphism}) and that it is nilpotent in this context (Theorem \ref{deltaisnilpotent}).
Example \ref{mishamykola} considers the case of $H=kG$, for $G$ a given group, in this case, we recover the cohomology theory
for partial group actions developed by M. Dokuchaev and M. Khrypchenko in \cite{DK1}. More specific examples, namely, the 
cohomology theory for partial group actions and partial group gradings over the base field are given in Subsection 3.4. More 
specifically, in Example \ref{australian} we calculate the first cohomology groups of a partial grading of the base field by 
the Klein four group, this extends extends what has been done in \cite{ABDP2}. 

In Section 4, given a partial action of $H$ on $A$, we first define the reduced cochain complex $\widetilde{C}_{par}^n (H,A) \cong C^n_{par}(H,A) /A^{\times}$ this new cochain complex produces the same cohomology groups as the original cochain complex $C^n_{par}$. Next, we consider the algebra $\wide{A}$, which is a quotient of the free commutative algebra generated by the images of all cochains $f\in \widetilde{C}_{par}^n (H,A)$. This defines a commutative and co-commutative Hopf algebra  over the commutative algebra $E(A)$, which is the subalgebra of $A$ generated by elements of the form $h\cdot 1_A$. One proves that $\widetilde{A}$ generates the same cohomologies as the original algebra $A$, that is, for any $n\in \mathbb{N}$ , we have the isomorphisms $H^n_{par} (H,\wide{A}) \cong H^n_{par} (H,A)$. 
Then, without loss of generality, we can consider only the cohomological theory for partial actions of a co-commutative Hopf 
algebra $H$ over a commutative and co-commutative Hopf algebra $\wide{A}$. Section 5 is devoted to analise twisted partial 
actions and partial crossed products \cite{ABDP1}. For the case of a co-commutative Hopf algebra $H$ and a commutative 
algebra $A$, all twisted partial actions are, in fact partial actions. Nonetheless, we still can get nontrivial partial 
crossed products by means of chosing a partial $2$-cocycle. In fact, the partial crossed products are classified by the 
second partial cohomology group $H^2_{par} (H,A)$. The new feature which appears in the context of partial actions is that 
the crossed product $\underline{\wide{A}\#_{\omega} H}$ has a structure of a Hopf algebroid over the same base algebra 
$E({A})$ (Theorem \ref{crossedproductishopfalgebroid}). This suggests that the cohomology theory for partial actions 
can be viewed as a cohomological theory for Hopf algebroids.

Section 6 is devoted to study partially cleft extensions and then we have Theorem 
\ref{partiallycleftisclefthopfalgebroid}, which we consider to be the main result of this paper. Basically, we proved that
given a partially cleft extension $B$ of a commutative co-commutative Hopf algebra $\wide{A}$ by a co-commutative Hopf 
algebra $H$, there exists a Hopf algebroid over the base subalgebra $E({A})$, namely, the partial smash product 
$\underline{E({A})\# H}$, such that $B$ is a $\underline{E({A})\# H}$-cleft extension of $\wide{A}$ in the sense
of G. B\"{o}hm and T. Brzezinski \cite{BB}. Finally, in Section 7 we give some perspectives for future developments.

For sake of a concise exposition, we omit several proofs which consists in basic manipulation of properties of partial
actions. The reader interested in more details can find them in the apendice of this paper.

\section{Mathematical preliminaries}

Given a field $k$, $k$-bialgebra or a Hopf algebra $H$ and a $k$-algebra $A$, for each $n\geq 0$ we have the convolution algebras $\mbox{Hom}_k (H^{\otimes n} ,A)$, with convolution product given by
\[
f\ast g (h^1 \otimes \cdots \otimes h^n )=f(h^1_{(1)} \otimes \cdots \otimes h^n_{(1)} ) g (h^1_{(2)} \otimes \cdots \otimes h^n_{(2)} ) ,
\]
and unit
\[
\mathbf{1} (h^1 \otimes \cdots \otimes h^n ) =\epsilon_H (h^1) \ldots \epsilon_H (h^n) \mathbf{1}_A .
\]

In particular, for $n=0$, we have
\[
\mbox{Hom}_k (H^{\otimes 0} ,A) =\mbox{Hom}_k (k ,A)\cong A.
\]

The following result can be easily obtained, we leave the details of the proof to the reader.

\begin{prop}\label{prop2.1}
Let $H$ be a cocommutative bialgebra, or Hopf algebra and $A$ be a commutative algebra. Then, for every $n\geq 0$ the 
convolution algebras $Hom_{k}(H^{\tens n},A)$ are commutative.
\end{prop}

Henceforth, for sake of simplicity, we will denote $f\hspace{-0.05cm}(h^1 \otimes \cdots \otimes h^n )$ just by $f\hspace{-0.05cm}(h^1\hspace{-0.1cm} , \ldots , h^n)$.

The main concept underlying this work is the notion of a partial action of a Hopf algebra $H$ over an algebra $A$. 

\begin{defi}\label{def2.2} \cite{CJ}
A partial action of a Hopf algebra $H$ over an algebra $A$ is a linear map 
$\cdot: H\otimes A \rightarrow A$, such that, for every $a,b \in A$ and $h,l \in H$, we have
\begin{enumerate}
\item[(PA1)] $1_H \cdot a = a$;
\item[(PA2)] $h\cdot (ab) = (h_{(1)}\cdot a)(h_{(2)}\cdot b)$;
\item[(PA3)] $h\cdot (l\cdot a) = (\sub{h}{1}\cdot 1_A)(\sub{h}{2}l\cdot a)$. 
\end{enumerate}
We say that the partial action is symmetric if, in addition, we have 
\begin{enumerate}
\item[(PA3')] $h\cdot (l\cdot a) = (\sub{h}{1}l\cdot a)(\sub{h}{2}\cdot 1_A)$ .
\end{enumerate}
The algebra $A$ is said to be a partial $H$-module algebra.
\end{defi}

Note that, if $H$ is a cocommutative Hopf algebra and $A$ is a commutative algebra, then every partial action of $H$ on $A$ is automatically symmetric. 

\begin{exmp}\label{groupactions} \cite{AB} Let $G$ be a group, recall that a partial action of $G$ over an algebra $A$ is a pair $(\{ A_g \}_{g\in G}, \{ \alpha_g :A_{g^{-1}}\rightarrow A_g \}_{g\in G})$, where $A_g$ is an ideal of $A$ for each $g\in G$ and $\alpha_g$ is an isomorphism of (not necessarily unital) algebras. We say that the partial action of $G$ on $A$ is unital if, for each $g\in G$, $A_g =1_g A$, where $1_g$ is a central idempotent in $A$ and $\alpha_g$ is a unital isomorphism between $A_{g^{-1}}$ and $A_g$. For the case where $H=kG$, the group algebra of $G$, symmetric partial actions of $kG$ are in one to one correspondence with unital partial actions of $G$. This correspondence can be easily seen: Given a unital partial action 
$(\{ A_g =1_g A \}_{g\in G}, \{ \alpha_g :A_{g^{-1}}\rightarrow A_g \}_{g\in G})$, one defines $\cdot :kG \otimes A \rightarrow A$, by $\delta_g \cdot a =\alpha_g (1_{g^{-1}}a)$. Conversely, given a symmetric partial action of $kG$ over $A$, define, for each $g\in G$, the idempotents $1_g =\delta_g \cdot 1_A$, by them, construct the ideals $A_g =1_g A$ and the isomorphisms $\alpha_g =\delta_g \cdot \underline{\quad} |_{A_{g^{-1}}}$.
\end{exmp} 

\begin{exmp}\label{inducedactions} \cite{AB} Given a Hopf algebra $H$, a left $H$-modules algebra $B$ and a central idempotent $e\in B$, one can define a partial action of $H$ on $A=eB$. Denoting by $\triangleright$ the action of $H$ over $B$, the induced  partial action is given by $h\cdot ea =e(h\triangleright (ea))$, for every $a\in B$ and $h\in H$.
\end{exmp}

The next two examples will be explored with more details throughout this paper for giving examples of cohomologies.

\begin{exmp}\label{groupactionsoverfield}
Consider a group $G$, let us see the partial actions of the Hopf algebra $H = kG$ over $A=k$, the base field. A partial 
action$\cdot: kG\otimes k \rightarrow k$, associates to each $g\in G$ the linear transformation 
$\delta_g \cdot \underline{\quad} :K \rightarrow k$, this is the same as  defining a linear functional 
$\lambda: kG \rightarrow k$. Denoting $\lambda(\delta_g)$ simply by $\lambda_g$, one can write $\delta_g \cdot a = \lambda_g a$, 
for every $a\in k$. Using the functional $\lambda$, The axiom (PA1) says that $\lambda_e = 1$, where $e$ is the neutral element 
of the group $G$. Axiom (PA2), in its turn, implies that $\lambda_g = \lambda_g\lambda_g$, for every $g\in G$, and consequently
$\lambda_g = 1$ or $\lambda_g = 0$. Define
\[
H=\{ g\in G \; | \; \lambda_g =1 \} ,
\]
it is clear that $e\in G$. Axiom (PA3) says that $\lambda_g\lambda_h = \lambda_g\lambda_{gh}$, this implies that for $g,h\in H$, we have $gh\in H$. Finally, putting $h=g^{-1}$ in the previous identity, we conclude that $g\in H$ implies that $g^{-1}\in H$, therefore $H$ is a subgroup of $G$. It is easy to see that the action is global if, and only if $H=G$. Then, we can label the partial actions of $k G$ over $k$ by the subgroups of $G$. 
\end{exmp}

\begin{exmp}\label{groupgradingsoverfield}

Let $G$ be a finite abelian group and consider 
$H = (kG)^* = \langle p_g |\ g\in G\rangle$, the dual of the group algebra, with bialgebra structure given by 
\[
p_g p_h = \delta_{g,h}p_g, \quad \mathbf{1}= \sum_{g\in G} p_g, \quad \Delta(p_g) = \sum_{g\in G} p_h\otimes p_{h^{-1}g}, \quad \epsilon(p_g) = \delta_{g,e}.
\]

As in the previous example, partial actions of $(kG)^*$ over $k$ are associated to a linear functional 
$\lambda:\kg\rightarrow k$ defined by $\lambda(p_g)= \lambda_{p_g} = p_g\cdot 1$. In this case, the axioms for a partial action (PA1), (PA2) and (PA3) become, respectively 
$$
\sum_{g\in G}\lambda_{p_g}  =  1; \ \ \ \ 
\lambda_{p_g}  =  \sum_{h\in G}\lambda_{p_h}\lambda_{p_{h^{-1}g}};\ \ \ \ 
\lambda_{p_g}\lambda_{p_h}  =  \lambda_{p_{gh^{-1}}}\lambda_{p_h} = \lambda_{p_{gh^{-1}}}\lambda_{p_g}.
$$

Defining $L = \{g\in G|\ \lambda_{p_g} \neq 0\}$, one can see that $L$ is a subgroup of $G$: First, as 
$\underset{g\in G}{\sum}\lambda_{p_g} = 1$ then there exists some $g\in G$ such that $\lambda_{p_g} \neq 0$, 
and therefore $L\neq \emptyset$. Moreover, given $g,\ h \in L$, the third equation says that 
$0\neq \lambda_{p_g}\lambda_{p_h} = \lambda_{p_{gh^{-1}}}\lambda_{p_g}$, which implies that $\lambda_{p_{gh^{-1}}} \neq 0$, 
and therefore $gh^{-1}\in L$.
 
In order to analyze the possible values of $\lambda_{p_g}$, for $g\in L$, take $g=h$ in the third equation, then  
$\lambda_{p_g}\lambda_{p_g} = \lambda_{p_{gg^{-1}}}\lambda_{p_g}
 = \lambda_{p_e}\lambda_{p_g}$. This implies that, $\lambda_{p_g} = \lambda_{p_e}, \ \forall\ g\in L$. Finally, from the first 
 equation,
\[
1 = \underset{g\in G}{\sum} \lambda_{p_g} = \underset{g\in L}{\sum} \lambda_{p_e} = |L|\lambda_{p_e},
\]
and therefore $\lambda_{p_g} = \lambda_{p_e} = \frac{1}{|L|}$, for all $g \in L$. We leave to the reader the verification that 
the action is global if, and only if, $L=\{ e \}$.

We conclude that the partial actions of $\kg$ over the base field $k$ are classified by subgroups of $G$ and given by
\[
\lambda_{p_g} = \left\{\begin{array}{cl} \frac{1}{|L|} & ,\ g\in L\\
     0 & ,\ c.c.                         
      \end{array}\right.
\]
\end{exmp}

\section{Cohomology for partial actions}\label{sec3}

In his 1968's seminal article, M.E. Sweeder presented a cohomology theory for commutative algebras which are modules over a
given cocommutative Hopf algebra. Basically, the cochain complexes $C^n (H,A)$ are defined as the abelian groups consisting 
of the invertible elements of the commutative convolution algebras $\mbox{Hom}_k (H^{\otimes n}, A)$. The main difference between
the cohomology theory of global and partial actions is that in the partial case one needs to find appropriated unital ideals 
in the convolution algebras in order to define correctly the cochain complexes. These ideals are constructed upon a system 
of idempotents for the convolution algebras. 

Henceforth, $H$ is always a cocommutative Hopf algebra acting partially on a commutative algebra $A$ with partial action 
$\cdot :H\otimes A \rightarrow A$.

\subsection{A system of idempotents for the convolution algebras}

We start in\-tro\-du\-cing some idempotent elements of $\mbox{Hom}_k(\hte ,A)$ which are very important throughout this paper. 
As the convolution algebras are commutative, for each $n\geq 0$, an idempotent is automatically a central idempotent. 
Moreover, the convolution pro\-duct of a finite number of idempotents is also an idempotent. In what follows, we introduce 
a nested system of idempotents in the convolution algebra related to the partial action.

\begin{prop}\label{prop3.1} \begin{enumerate}
\item For each $n\geq 1$, the linear maps
$$\func{\widetilde{e}_n :}{H^{\tens n}}{A}{(h^1\tens \cdots \tens h^n)}{(\subp{h}{1}\ldots\subp{h}{n})\cdot 1_A}$$
are idempotent in the corresponding convolution algebras $\mbox{Hom}_k(\hte,A)$.
\item Let $n<m$ and $e\in \mbox{Hom}_{k}(\hte, A)$ an idempotent, then
$$e_{n,m} = e \tens \underset{m-n}{\underbrace{\epsilon \tens \cdots \tens \epsilon}} \in Hom_{k}(H^{\tens m}, A)$$
is an idempotent in $\mbox{Hom}_k (H^{\otimes m} , A)$.
\end{enumerate}
\end{prop}

\noindent \dem The proof of item (1) follows from basically from item (PA2) of definition of a partial action. Item (2) is straightforward.
\find

\begin{defi}\label{e_{l,n}}
For $n\geq 1$ and $1 \leq l \leq n$ arbitrary, we define:
$$\widetilde{e}_{l,n} := \widetilde{e}_{l} \tens \underset{n-l}{\underbrace{\epsilon \tens \cdots \tens \epsilon}}.$$
\end{defi}

Note that, for $l=n$ in the definition above $\widetilde{e}_{n,n}=\widetilde{e}_{n}$.

\begin{defi}\label{e_n}
For $n\geq 1$, we define
$$e_n := \widetilde{e}_{1,n}*\widetilde{e}_{2,n}*\cdots * \widetilde{e}_{n} \in \mbox{Hom}_{k} (\hte, A).$$
\end{defi}

The next proposition gives us a useful characterization of the idempotents $e_n \in Hom_{k} (\hte, A)$.

\begin{prop}\label{prop3.4}
For any $(\coe) \in  \hte$ we have that
$$e_n (\coen) = \subp{h}{1}\cdot(\subp{h}{2}\cdot(\cdots\ \cdot(\subp{h}{n}\cdot 1_A)\cdots)).$$
\end{prop}

\noindent\dem The proof of this result consists basically in elementary manipulations of itens (PA2) and 
(PA3) of the definition of a partial action.\find

\subsection{Cochain complexes}

Based on the idempotents defined in the previous subsectin, one can define the cochain complexes for the partial action of $H$ on $A$. For each $n>0$ define the following ideals of $\mbox{Hom}_k (\hte ,A)$:
\[
I(\hte , A)= e_n \ast \mbox{Hom}_k (\hte ,A) =  \{ e_n * g : g \in Hom_{k} (\hte ,A)  \}.
\]

As $e_n$ is a central idempotent, the ideal $I(\hte , A)$ can be considered as a unital algebra with unit $e_n$. 
An element $f\in I(\hte , A)$ is said to be (convolution) invertible in this ideal if there is another element
$g\in I(\hte , A)$ such that $f\ast g= g\ast f =e_n$.

\begin{defi}
Let $H$ be a cocommutative bialgebra and  $A$ be a partial $H$-modules algebra with partial action
$\cdot: H\otimes A \rightarrow A$. For $n>0$ , a ``partial'' $n$-cochain of $H$ taking values in $A$ 
is an invertible element in the ideal $I(\hte,A)$. Denoting by $C^{n}_{par}(H,A)$ the set of $n$ cochains, 
we have that $C^n_{par} (H,A) = (I(\hte,A))^{\times}$. 
For $n=0$ we say that a $0$-cochain is an invertible element in the agebra $A$, that is  $C^{0}_{par}(H,A) = A^{\times}$.
\end{defi}

Note that $C^{n}_{par}(H,A)$ is an abelian group with relation to the convolution product, while $C^{0}_{par}(H,A) = A^{\times}$, is an abelian group with the ordinary multiplication in $A$ and the unit $e_0 = 1_A$. 

\begin{defi}\label{delta}
For an arbitrary $f\in C^{n}_{par}(H,A)$, $(h^1\otimes\cdots\otimes h^{n+1})\in H^{\otimes n+1}$, we define the ``partial'' coboundary operator
\begin{eqnarray*}
& \, & (\delta_n f)(h^1, \ldots, h^{n+1})  =   (h^1_{(1)}\cdot f(h^2_{(1)},\ldots, h^{n+1}_{(1)}))*\\
   & & *\prodt{i=1}{n}f^{(-1)^i}(h^1_{(i+1)},\ldots, h^i_{(i+1)}h^{i+1}_{(i+1)},\ldots, h^{n+1}_{(i+1)})  *f^{(-1)^{n+1}}(h^1_{(n+2)},\ldots,h^n_{(n+2)})
\end{eqnarray*}

If $n=0$ and $a\in A^\times$, we have $(\delta_0 a)(h) = (h\cdot a)a^{-1}$.
\end{defi}
   
The challenge is to prove that the coboundary operators are well defined, that is, for every $f\in C^{n}_{par}(H,A)$, the map $\delta_n f$ is indeed in $C^{n+1}_{par}(H,A)$. Moreover, one needs to prove that the sequence
\[
C^0_{par} (H,A) \stackrel{\delta_0}{\rightarrow} C^1_{par} (H,A) \stackrel{\delta_1}{\rightarrow} \cdots 
\stackrel{\delta_{n-1}}{\rightarrow} C^n_{par} (H,A) \stackrel{\delta_n}{\rightarrow} C^{n+1}_{par} (H,A) 
\stackrel{\delta_{n+1}}{\rightarrow} \cdots 
\]
is a cochain complex, that is, each $\delta_n$ is a homomorphism of abelian groups between 
$C^n_{par} (H,A)$ and $C^{n+1}_{par} (H,A)$ satisfying  $\delta_{n+1}\circ \delta_n (f) = e_{n+2}$, for each $f\in C^n_{par} (H,A)$. For this purpose, we introduce some auxiliary operators which will help us to describe the coboundary operators in a more intrinsic way and whose properties will lead us to the desired results.

\begin{defi} \label{auxiliares} 
\begin{enumerate} 
\item For each $n\geq 0$, define the map  
$$\func{E^n:}{C^n_{par}(H,A)}{\mbox{Hom}_k(H^{\otimes n+1},A)}{f}{E^n(f)},$$
given by $E^n(f)(h^1,\ldots, h^{n+1}):= h^1\cdot f(h^2,\ldots, h^{n+1})$.
\item For $n<m$, define the map,
$$\func{i_{n,m} :}{C^n_{par}(H,A)}{\mbox{Hom}_k (H^{\otimes m} ,A)}{f}{i_{n,m}(f)}$$
given by $i_{n,m}(f)(h^1, \ldots, h^{m}) := f(\coen)\epsilon(h^{n+1})\ldots\epsilon(h^{m})$.
\item For $i\in\{1,\ldots, n\}$, define the map 
$$\func{\mu_i:}{H^{\tens n+1}}{\hte}{(\subp{h}{1}\tens \cdots\tens \subp{h}{n+1})}{(\subp{h}{1}\tens \cdots\tens \subp{h}{i}\subp{h}{i+1}\tens \cdots\tens \subp{h}{n+1}).}$$
\end{enumerate}
\end{defi}

With these auxiliary operators, the coboundary operator can be rewritten as
$$\func{\delta_n:}{ C^{n}_{par}(H,A)}{ C^{n+1}_{par}(H,A)}{f}{\delta_n(f) := E^n(f)* \prodt{i=1}{n} f^{(-1)^{i}}\circ \mu_i * i_{n,n+1}(f^{(-1)^{n+1}})}$$
and the properties of $\delta_n$ are based upon the properties of these operators.\\

\begin{lemma} \label{p1} 
\begin{enumerate}[(i)]
\item For $f,\ g \in C^n_{par}(H,A)$, we have $E^n(f\ast g) = E^n(f)\ast E^n(g)$.
\item $E^n(e_n) = e_{n+1}$.
\item For $n<m$ and $f,\ g \in C^n_{par}(H,A)$, we have $i_{n,m}(f\ast g)= i_{n,m}(f)\ast i_{n,m}(g)$.
\item $i_{n,m}(e_n )\ast e_m =e_m$.
\item For $f,\ g \in C^n_{par}(H,A)$, we have $(f*g)\circ\mu_i = (f\circ\mu_i)*(g\circ\mu_i)$,  $\forall \ i\in \{1,\ldots, n\}$.
\item $(e_n\circ\mu_n)*i_{n,n+1}(e_n)= e_{n+1}$.
\item $(e_n\circ\mu_i)*e_{n+1} = e_{n+1}$, $\forall\ i\in\{1,\cdots, n-1\}$.
\end{enumerate}
\end{lemma}

\noindent\dem The proof of these items is done by a long, but straightforward computation evaluating on elements of the 
corresponding tensor product. \find

\begin{thm} \label{deltaismorphism} For any $n\geq 0$, $f\in C^n_{par} (H,A)$, the linear map
$\delta_n (f): H^{\otimes n+1}\rightarrow A$ is an element of $C^{n+1}_{par} (H,A)$. Moreover, the map 
$\delta_n :C^n_{par} (H,A) \rightarrow C^{n+1}_{par} (H,A)$ is a morphism of abelian groups.
\end{thm}

\noindent\dem If $f\in C^n_{par} (H,A)$, then$f=f\ast e_n$ and it is invertible by convolution. Consider the expression of  $\delta_n (f)$, 
\[
\delta_n(f) := E^n(f)* \prodt{i=1}{n} f^{(-1)^{i}}\circ \mu_i * i_{n,n+1}(f^{(-1)^{n+1}}) .
\]
Items (i), (iii) and (v) of Lemma \ref{p1} and using the fact that the convolution algebra $\mbox{Hom}_k (H^{\otimes n+1} ,A)$ 
is commutative, we conclude that $\delta_n (f\ast g)=\delta_n (f) \ast \delta_n (g)$, in particula
r$\delta_n (f)=\delta_n (f\ast e_n) =\delta_n (f) \ast \delta_n (e_n)$. By item (ii), we know that $E^n (e_n )=e_{n+1}$ and by
items (iv), (vi) and (vii) we see that the unit $e_{n+1}$ absorbs the other factors, leading to $\delta_n (e_n )=e_{n+1}$. 
Then we have $\delta_n (f)=\delta_n (f) \ast e_{n+1}$. A straightforward calculation leads us to 
$\delta_n (f^{-1})= (\delta_n (f))^{-1}$ Therefore, we prove that $\delta_n$ is well defined and it is a morphism of 
abelian groups. \find

In order to prove that $\delta_{n+1} \circ \delta_n (f) =e_{n+2}$, for every $f\in C^{n}_{par} (H,A)$, we need the following lemma.

\begin{lemma} \label{p2} Let $f\in C^{n-1}_{par} (H,A)$, then 
\begin{enumerate}[(i)]
\item $E^n(i_{n-1,n}(f)) = i_{n,n+1}(E^{n-1}(f))$.
\item $(e_{n-1}\circ\mu_i\circ\mu_{i+1})*e_{n+1} = e_{n+1}$, $\forall\ i\in\{1,\cdots, n-1\}$.
\item $(e_{n-1}\circ\mu_i\circ\mu_{i+j})*e_{n+1} = e_{n+1}$, $\forall \ i\in\{1,\cdots, n-1\},\ j\in \{2,\cdots, n-i\}$.
\item $E^n(f\circ \mu_i) = E^{n-1}(f)\circ \mu_{i+1}$, $\forall\ i\in\{1,\ldots, n-1\}$.
\item $E^n\circ E^{n-1}(f) = i_{1,n+1}(\widetilde{e}_1)*(E^{n-1}(f)\circ\mu_1)$.
\item $i_{n,n+1}(f\circ \mu_i)*i_{n-1,n}(f^{-1})\circ\mu_i = i_{n,n+1}(e_{n-1}\circ\mu_i)$, $\forall\ i\in\{1,\ldots, n-1\}$.
\item $(i_{n-1,n}(f)\circ \mu_n)*i_{n-1,n+1}(f^{-1}) = i_{n-1,n+1}(e_{n-1})$.
\item $(f\circ\mu_i\circ\mu_i)*(f^{-1}\circ\mu_i\circ\mu_{i+1}) = e_{n-1}\circ\mu_i\circ\mu_{i}$, $\forall\ i\in\{1,\ldots, n-1\}$.
\item $(f\circ\mu_i\circ\mu_{i+j})*(f^{-1}\circ\mu_{i+j-1}\circ\mu_i) = e_{n-1}\circ\mu_i\circ\mu_{i+j}$, $\forall\ i\in\{1,\ldots, n-2\}$, $\forall\ j \in \{2,\ldots, n-i\}$.
\end{enumerate}
\end{lemma}

\noindent\dem The proofs of these items come after a long but straightforward calculation. 
\vspace{-0.2cm}
\begin{flushright}
 \find
\end{flushright}

With this lemma, one can prove the following result.

\begin{thm} \label{deltaisnilpotent} For any  $f\in C^{n}_{par}(H,A)$, we have that $\delta_{n+1}\circ\delta_n (f) = e_{n+2}$.
\end{thm}
\noindent\dem For $f\in C^{n}_{par}(H,A)$, we have
$$\begin{array}{rcl}
& & \delta_{n+1}(\delta_n(f)) = 
 E^{n+1}(\delta_n(f))* \prodt{i=1}{n+1} (\delta_n (f) )^{(-1)^{i}}\circ \mu_i *i_{n+1,n+2}((\delta_n (f))^{(-1)^{n+2}})\\
& = & E^{n+1}(E^n(f)* \prodt{j=1}{n} f^{(-1)^{j}}\circ \mu_j * i_{n,n+1}(f^{(-1)^{n+1}})) \\
& & \;  *\prodt{i=1}{n+1} (E^n(f^{(-1)^{i}})* \prodt{j=1}{n} f^{(-1)^{i+j}}\circ \mu_j *i_{n,n+1}(f^{(-1)^{n+i+1}}))\circ \mu_i  \\
& & \; * i_{n+1,n+2}(E^n(f^{(-1)^{n+2}})* \prodt{j=1}{n} f^{(-1)^{n+j+2}}\circ \mu_j * i_{n,n+1}(f^{(-1)^{2n+3}})) .
\end{array}$$

Then, using the items of Lemma \ref{p1} and Lemma \ref{p2}, one can show that this expression is equal to $e_{n+2}$.
\find

Therefore, we ended up with a cochain complex $(C^n_{par} (H,A) ,\delta_n)_{n\in \mathbb{N}}$, which allows us to define a cohomology theory.

\subsection{Cohomologies}

\begin{defi} Let $H$ be a cocommutative Hopf algebra acting partially over a commutative algebra $A$ and consider the cochain complex $(C^n_{par} (H,A) ,\delta_n)_{n\in \mathbb{N}}$ defined in the previous sections. For $n>0$, define the abelian groups
$Z^n(H,A) =\ker \delta_n$, $B^n(H,A) = Im \ \delta_{n-1}$ e $H^n(H,A) = \ker \delta_n/\ Im \ \delta_{n-1}$, respectively,
as the groups of partial $n$-cocycles, partial $n$-coboundaries and partial $n$-cohomologies of $H$ taking values in $A$.  
$n\geq 1$. For $n=0$, define $H^0(H,A) = Z^0(H,A) = \ker \delta_0$.
\end{defi}

Let us characterize the partial cocycles and the partial coboundaries for $n=0$, $1$ and $2$. 

For $n=0$, we have by definition
$$H^0(H,A) = Z^0(H,A) = \{a \in A^{\times}|h\cdot a = (h\cdot 1_A)a, \ \forall h\in H\}. $$

Thus the partial $0$-cocycles are the elements of $A$ invariant under the partial action as defined in \cite{AB2}.

For $n=1$, the partial $1$-coboundaries are
$$B^1(H,A) = \mbox{Im} \delta_0  = \{ f\in C^1_{par}(H,A)|\ \exists a\in A^{\times}, \quad f(h) = \delta_0 (a)(h) \} ,$$
 this means 
\[
B^1(H,A) = \{f\in C^1_{par}(H,A) \; |\; \exists a\in A^{\times}, \quad f(h) = (h\cdot a)a^{-1}\} .
\]   

Also, for $f\in C^1_{par}(H,A)$, we have
$$\begin{array}{ccl}\delta_1 (f)(h,l) & = & E^2(f)(h_{(1)},l_{(1)})f^{-1}(h_{(2)}l_{(2)})f(h_{(3)})\epsilon(l_{(3)})\\ 
  & = & (h_{(1)}\cdot f(l_{(1)}))f^{-1}(h_{(2)}l_{(2)})f(h_{(3)})\epsilon(l_{(3)}).
  \end{array}$$
  
Then, forall $h,\ l\in H$, the partial $1$-cocycles are 
\begin{eqnarray*} 
& & Z^1(H,A)= \{f\in C^1_{par}(H,A)|\ \delta_1(f)(h,l) = e_2(h,l)\}\\
&=& \{f\in C^1_{par}(H,A)|\  (h_{(1)}\cdot f(l_{(1)}))f^{-1}(h_{(2)}l_{(2)})f(h_{(3)})= h\cdot(l\cdot 1_A)\}\\
&=& \{f\in C^1_{par}(H,A)|\  (h_{(1)}\cdot f(l_{(1)}))f(h_{(3)})= (h_{(1)}\cdot (l_{(1)} \cdot 1_A))f(h_{(2)}l_{(2)})\} .
\end{eqnarray*}

Due to the fact that for a $1$-cocycle $f$ we have $f=e_1 \ast f$, then the condition of $1$-cocycle can also be rewritten as
\[
(h_{(1)}\cdot f(l_{(1)}))f(h_{(3)})= (h_{(1)}\cdot  1_A)f(h_{(2)}l_{(2)}) ,
\]

For $n=2$, we have the partial $2$-coboundaries
\begin{eqnarray*}& & B^2(H,A)  =  \{g\in C^2_{par}(H,A)|\ \delta_1(f)(h,l) = g(h,l)\}\\ 
 & = & \{g\in C^2_{par}(H,A)|\ g(h,l) = (h_{(1)}\cdot f(l_{(1)}))f^{-1}(h_{(2)},l_{(2)})f(h_{(3)})\}.
 \end{eqnarray*}

Also, for $f\in C^2_{par}(H,A)$, we have
$$\begin{array}{cl} 
 & \delta_2 (f)(h,l,m) = E^2(f)*\prodt{i=1}{2}f^{(-1)^{i}}\circ\mu_i*i_{2,3}(f^{(-1)^3})(h,l,m)\\
= & (h_{(1)}\cdot f(l_{(1)},m_{(1)}))f^{-1}(h_{(2)}l_{(2)},m_{(2)})f(h_{(3)},l_{(3)}m_{(3)})f^{-1}(h_{(4)},l_{(4)})\epsilon(m_{(4)}). 
\end{array}$$

Then, the partial $2$ cocycles are

\begin{eqnarray*} & Z^2(H,A) = \{f\in C^2_{par}(H,A)|\ \delta_2(f)(h,l,m) = e_3(h,l,m),\ \forall\  h,l,m\in H\}\\
= & \{f\in C^2_{par}(H,A)|\  (h_{(1)}\cdot f(l_{(1)},m_{(1)}))f^{-1}(h_{(2)}l_{(2)},m_{(2)})f(h_{(3)},l_{(3)}m_{(3)})\\
&\hspace{4.4cm} f^{-1}(h_{(4)},l_{(4)})  =  h\cdot(l\cdot(m\cdot 1_A)),\ \forall\  h, l, m\in H\}\\
=& \{f\in C^2_{par}(H,A)|  (h_{(1)}\cdot f(l_{(1)},m_{(1)}))f(h_{(2)},l_{(2)}m_{(2)}) = \\
& \hspace{1.2cm} (h_{(1)} \cdot (l_{(1)} \cdot (m_{(1)}\cdot 1_A)))f(h_{(2)}l_{(2)},m_{(2)})f(h_{(3)},l_{(3)}),\ \forall\  h, l, m\in H\} .
\end{eqnarray*}

Again by absorption of units, one can rewrite the condition of $2$-cocycle as
\[
(h_{(1)}\cdot f(l_{(1)},m_{(1)}))f(h_{(2)},l_{(2)}m_{(2)}) 
=  (h_{(1)} \cdot  1_A) f(h_{(2)}l_{(2)},m_{(2)})f(h_{(3)},l_{(3)}) ,
\]
which is the form presented in \cite{ABDP1}.

\begin{exmp} In the case of a global action of $H$ over $A$, which is equivalent to say that $h\cdot 1_A =\epsilon (h)1_A$, $\forall h\in H$, the cochain complexes are simply given by $C^n (H,A)=\mbox{Hom}_k (H^{\otimes n},A)^{\times}$. Then we recover exactly the cohomology theory obtained by Sweedler in \cite{Swe}.
\end{exmp}

\begin{exmp} \label{mishamykola} Let $G$ be a group and $H=kG$, the group algebra of $G$. Using the canonical basis $\{ \delta_g \in kG \; |\; g\in G \}$, the axioms (PA1), (PA2)  and (PA3) of partial action read 
\begin{enumerate}
\item[(PA1)] $\delta_e\cdot a = a$, for every $a\in A$;
\item[(PA2)] $\delta_g\cdot(ab) = (\delta_g\cdot a)(\delta_g\cdot b)$, for every $g\in G$ and $a,b\in A$;
\item[(PA3)] $\delta_g\cdot(\delta_h\cdot a) = (\delta_g \cdot 1_A)(\delta_{gh}\cdot a)$, for every $g,h\in G$ and $a\in A$.
\end{enumerate}

In order to calculate the partial $n$-cocycles, partial $n$-coboundaries and partial $n$-cohomologies, we denote the coboundary operator by $\partial_n$ instead of $\delta_n$ to avoid confusion with the elements $\delta_g \in kG$.

For $n=0$, we have
$$H^0(kG,A) = Z^0(kG,A) = \{a \in A^{\times}|(\delta_g\cdot a)a^{-1} = (\delta_g\cdot 1_A), \ \forall \delta_g\in kG\}, $$

For $n=1$, the $1$-coboundaries are
\begin{eqnarray*}
B^1(kG,A) & = &\{f\in C^1_{par}(kG,A)|\ \exists a\in A^{\times}, \quad f(\delta_g) = \partial_0(\delta_g)(a) \} \\ 
&=& \{f\in C^1_{par}(kG,A)|\ f(\delta_g) = (\delta_g\cdot a)a^{-1}\} .
\end{eqnarray*} 

Also, we have, for every $f\in C^1_{par}(kG,A)$,
$$\partial_1 (f)(\delta_g,\delta_h) = (\delta_{g}\cdot f(\delta_h))f^{-1}(\delta_{gh})f(\delta_g).$$

Then, for all  $\delta_g, \delta_h\in kG$, we obtain the $1$-cocycles
\begin{eqnarray*}
 Z^1(kG,A) & = & \{f\in C^1_{par}(kG,A)|\ \partial_1(f)(\delta_g,\delta_h) = \delta_g\cdot(\delta_h\cdot 1_A) \}\\
&=& \{f\in C^1_{par}(kG,A)|\  (\delta_g\cdot f(\delta_h))f(\delta_g)= (\delta_g \cdot 1_A) f(\delta_{gh})\}.
\end{eqnarray*}

For $n=2$, the $2$ coboundaries are
$$\begin{array}{rcl}B^2(kG,A)  & = &\{i\in C^2_{par}(kG,A)|\ \partial_1(f)(\delta_g,\delta_h) = i(\delta_g,\delta_h)\}\\ 

&=& \{i\in C^2_{par}(kG,A)|\ i(\delta_g,\delta_h) = (\delta_g\cdot f(\delta_h))f^{-1}(\delta_{gh})f(\delta_g)\}.
   
  \end{array}
 $$
 
Moreover, for $f\in C^2_{par}(kG,A)$
$$\partial_2 (f)(\delta_g,\delta_h,\delta_l) = (\delta_g\cdot f(\delta_h,\delta_l))f^{-1}(\delta_{gh},\delta_l)f(\delta_g,\delta_{hl})f^{-1}(\delta_g,\delta_h).
$$

Then, for all $\delta_g,\delta_h,\delta_l\in kG $, the partial $2$-cocycles are
$$\begin{array}{cl}
&Z^2(kG,A)  =  \{f\in C^2_{par}(kG,A)|\ \partial_2(f)(\delta_g,\delta_h,\delta_l) = e_3(\delta_g,\delta_h,\delta_l)\}\\
=& \{f\in C^2_{par}(kG,A)|\  (\delta_g\cdot f(\delta_h,\delta_l))f(\delta_g,\delta_{hl}) = (\delta_g \cdot 1_A) f(\delta_{gh},\delta_l)f(\delta_g,\delta_h)\}.
\end{array}$$

This cohomology for partial actions of the group algebra $kG$ corresponds to the cohomology for partial group actions
described in \cite{DK1}. Recall from Example \ref{groupactions} that there is a one-to-one correspondence between partial
actions f $kG$ and unital partial actions of the group $G$, given by $A_g =1_g A$ in which $1_g =\delta_g \cdot 1_A$,
and $\alpha_g =\delta_g\cdot \underline{\quad}|_{A_{g^{-1}}}$. For elements $x_1, \ldots x_n \in G$ we define the ideals
\[
A_{(x_1,\ldots, x_n)} := A_{x_1}A_{x_1 x_2}\ldots A_{x_1 \ldots x_n},
\]
in which $A_{x_i}=1_{x_i} A$. This expression for the ideals are natural, considering the units 
\[
e_n (\delta_{x_1} ,\ldots ,\delta_{x_n})= \delta_{x_1} \cdot (\delta_{x_2} \cdot (\cdots  (\delta_{x_n} \cdot 1_A)))= 1_{x_1} 1_{x_1 x_2}\ldots 1_{x_1 \ldots x_n }.
\]

The set of these ideals forms a semilattice, because the product of two ideals of this type is also an ideal of this type, this product is commutative and each ideal is idempotent, that is 
$A_{(x_1,\ldots, x_n)}=A_{(x_1,\ldots, x_n)}A_{(x_1,\ldots, x_n)}$. This can be viewed easily by the properties of the system of idempotents presented before. 

The correspondence between the cochain complexes presented here and those presented in \cite{DK1} can be viewed more exactly by the identification of the convolution algebra $Hom_k (kG^{\otimes n},A)$ with the algebra of functions 
$Fun (G^n , A)$, moreover, the functions $f:G^n \rightarrow A$ can also be viewed as a collection of elements of $A$ indexed by $n$-tuples in $G$, that is $f(g_1, ..., g_n) =f_{g_1, ...g_n} \in A$. As the canonical basis elements $\delta_g$, for $g\in G$ are group-like, the convolution product is in fact the pointwise product, that is, for $f^1, f^2 \in Fun(G^n,A)$ and $g_1,\ldots,g^n \in G$, we have 
\[ 
f^1 \ast f^2 (g_1 , ..., g_n)= f^1 (g_1 , ..., g_n)  f^2 (g_1 , ..., g_n) = f^1_{g_1 , ..., g_n}f^2_{g_1 , ..., g_n}.
\]

Therefore, the $n$-cochains $C^n_{par}(kG, A)$ coincide with the $n$ cochains $C^n_{par}(G,A)$.

The partial $n$-cocycles, partial $n$-coboundaries and partial $n$-cohomologies in the group setting are written as.
$$H^0(G,A) = Z^0(G,A) = \{a \in A^{\times}|(\alpha_g (1_{\inv{g}}a))a^{-1} = 1_g, \ \forall g\in G\}, $$

For $n=1$ the partial $1$-coboudaries are
$$\begin{array}{rcl}B^1(G,A) & = &\{f\in C^1_{par}(G,A)|\ \exists a\in A^{\times}, \quad  f(g) = \partial_0(g)(a) \}\\ 
&=& \{f\in C^1_{par}(G,A)|\ \exists a\in A^{\times}, \quad f(g) = (\alpha_g (1_{\inv{g}}a))a^{-1}\}.
  \end{array} $$
  
Moreover, for $f\in C^1_{par}(G,A)$ we have 
$$\partial_1 (f)(g,h) = (g\cdot (1_{\inv{g}}f(h)))f^{-1}(gh)f(g) .$$

Then, the partial $1$-cocycles are
$$\begin{array}{rcl}Z^1(G,A) &=& \{f\in C^1_{par}(G,A)|\ \partial_1(f)(g,h) = e_2(g,h),\ \forall\  g, h\in G\}\\
&=& \{f\in C^1_{par}(G,A)|\  (\alpha_g (1_{\inv{g}}f(h)))f(g)= 1_g f(gh),\ \forall\   g, h\in G\} .
\end{array}$$
Note that $\delta_g\cdot(\delta_h\cdot 1_A) = (\delta_g\cdot1_A)(\delta_{gh}\cdot 1_A) = 1_g1_{gh}$, and $1_{gh}$ is absorbed by $f(gh)$.

For $n=2$, the partial $2$-coboundaries are
$$\begin{array}{rcl}B^2(G,A)  & = &\{i\in C^2_{par}(G,A)|\ \partial_1(f)(g,h) = i(g,h)\}\\ 

&=& \{i\in C^2_{par}(G,A)|\ i(g,h) = (\alpha_g (1_{\inv{g}}f(h)))f^{-1}(gh)f(g)\}.
   
  \end{array}
 $$
 
For $f\in C^2_{par}(G,A)$, we have,
$$\partial_2 (f)(g,h,l) = (g\cdot (1_{\inv{g}}f(h,l)))f^{-1}(gh,l)f(g,hl)f^{-1}(g,h) ,
$$
then,
$$\begin{array}{l}  Z^2(G,A) = \{f\in C^2_{par}(G,A)|\ \partial_2(f)(g,h,l) = e_3(g,h,l),\ \forall\  g,h,l\in G\}\\
= \{f\in C^2_{par}(G,A)|\  (\alpha_g (1_{\inv{g}}f(h,l)))f(g,hl) = 1_g f(gh,l)f(g,h),\ \forall\  g,h,l\in G\} .
\end{array}$$

Again, the appearance only of $1_g$ in the right hand side of the $2$-cocycle condition is due to absorption of units.

Therefore, the cohomology obtained is the same as \cite{DK1}.
\end{exmp}

In the next subsection we will give more specific examples of cohomologies for partial actions in which the algebra $A$ is the base field $k$.

\subsection{Cohomology for partial actions on the base field}

\begin{exmp} {\bf (Partial group actions over the base field)} Let $G$ be a group. We have already seen in Example 
\ref{groupactionsoverfield} that partial actions $k G$ over $k$ are in correspondence with subgroups $L\leq G$ by means of the linear functional 
$$\func{\lambda:}{kG}{k}{\delta_g}{\lambda_g = \lambda(\delta_g)= \left\{\begin{array}{cl} 1 & ,\ g\in L\\
     0 & ,\ c.c.                          
      \end{array}\right.} .$$

Fix the subgroup $L$ of $G$ which defines the partial action. Let us now calculate the cohomologies $H^n_{par} (kG,k)$ (we use the symbol $\partial_n$ for the coboundary map to avoid confusion with the basis elements $\delta_g \in kG$):

\begin{itemize}
\item For $n=0$, $C^0_{par}(k G, k) = k^\times =k\backslash \{ 0\}$. Let $a\in C^0$, then,
$$(\partial_0 a)(\delta_g) = (\delta_g\cdot a)a^{-1} = \lambda_g a a^{-1} = \lambda_g = \left\{\begin{array}{ll} 1, & g\in L\\ 0, & g\not\in L\end{array}\right.$$
Therefore, $H_{par}^0 = Z_{par}^0 = C^0_{par} = k^\times$.
\item For $n=1$, let $f\in Z^1_{par}(k G,k)$. Then,
$$(\partial_1 f)(\delta_g,\delta_h)= \lambda_g f(\delta_h)f^{-1}(\delta_{gh})f(\delta_g) = \lambda_g\lambda_h$$
Denote $f(\delta_g)$ simply by $f(g)$ using the identification between the convolution algebra and the algebras of functions $f:G\rightarrow k$. If $g,h \in L$, then the $1$ cocycle condition can be rewritten as $f(gh)=f(g)f(h)$, which means that $f_{|_{L}}: k G\rightarrow k^\times$ is a character of the subgroup $L$.  If $g\notin L$ then it is easy to see that for a $1$ cocycle $f$, we have $f(g)=0$. As the partial $1$-coboundaries are given by $\lambda: G \rightarrow k$ such that 
$\lambda(g) = 1$ for any $g\in L$, we have that $H^1 = Z^1/ B^1$ are given by the nontrivial 1-dimensional representations of the subgroup $L$ which determines the partial action.
\item For $n=2$, let $\omega \in Z^2_{par}(k G,k)$. Then, denoting $\omega(\delta_g,\delta_h)$ simply by $\omega_{g,h}$, we have

$$(\partial \omega)(\delta_g,\delta_h,\delta_l) =\lambda_g \omega_{h,l}\omega_{gh,l}^{-1}\omega_{g,hl} \omega_{g,h}^{-1} =  \lambda_g\lambda_h\lambda_l$$

It is easy to see from the identity above that, if $(g,h)\notin L\times L$ then $\omega_{g,h}=0$. Then, defining $\omega: G\times G \rightarrow k$ by 
\[
\omega(g,h) = \left\{\begin{array}{ll}0, & (g,h)\not\in L\times L\\ 
\omega(\delta_g,\delta_h), & (g,h)\in L\times L\end{array}\right. ,
\] 
we have that the partial $2$-cocycles relative to $G$ are in fact usual $2$-cocycles of the subgroup $L$ \cite{Ade,W}, in other words 
$$Z^2_{par}(kG ,k)=Z^2(L,k).$$
\end{itemize}
\end{exmp}

\begin{exmp} {\bf (Partial group gradings over the base field)} Let $G$ be a finite abelian group. In Example \ref{groupgradingsoverfield} we saw that the partial actions of the Hopf algebra 
$H = \kg =\ <p_g\ |\ g\in G>$ over the base field $k$ are in one-to-one correspondence with subgroups $L\leq G$, namely
$$\lambda_{p_g} = \left\{\begin{array}{cl} \frac{1}{|L|} & ,\ g\in L\\
     0 & ,\ c.c.                          
      \end{array}\right.$$
      
Let us now calculate explicitly the partial cohomologies for $\kg$. 

For $n= 0$, recalling that $\delta_a (h) = (h\cdot a)a^{-1}$, for every $a\in k^{\times}$, we have 
$\delta_a( p_g) = \lambda_{p_g}$, and this leads to $Z^0(\kg,k)= C^0(\kg,k) = H^0(\kg,k) = k^\times$. Moreover, the 
$1$-coboundaries are basically given by the functional $\lambda$.

For $n=1$, let $\omega: \kg \ri k$ be a partial $1$-cocycle, then
\begin{eqnarray*}
 \lambda_{p_g}\lambda_{p_h} = \delta(\omega)(p_g,p_h) & = & 
\underset{l,m,i \in G}{\sum}(p_{ml^{-1}}\cdot\omega(p_{hi^{-1}}))\overline{\omega}(p_lp_i)\omega(p_{m^{-1}g}) \\
& = &\sum_{l,m \in G} \lambda_{ml^{-1}} \omega(p_{hl^{-1}}))\overline{\omega}(p_l)\omega(p_{m^{-1}g}) .
\end{eqnarray*}

Recalling that $e_1 *\omega = \omega = \omega*e_1 $ e que $e_1 (p_g) = \lambda_{p_g}$, we have
$$(e_1 *\omega)(p_g) = \underset{h\in G}{\sum}\lambda_{p_h}\omega(p_{h^{-1}g}) \Rightarrow \dfrac{1}{|L|}
\underset{h\in L}{\sum}\omega(p_{gh^{-1}}) = \omega(p_g) = \dfrac{1}{|L|}
\underset{h\in L}{\sum}\omega(p_{h^{-1}g}) .$$

This means that
$$\omega(p_g) = \dfrac{1}{|L|}\underset{h\in L}{\sum}\omega(p_{hg})= \dfrac{1}{|L|}\underset{h\in L}{\sum}\omega(p_{gh}).$$

For any $g\in L$, we have
$\omega(p_g) = \dfrac{1}{|L|}\underset{h\in L}{\sum}\omega(p_{h}) ,$
which is an invariance by translations in the subgroup. Furthermore, using the normalization condition, $\omega(1) = 1$, we have that
$$|L|\omega(p_e) = \underset{g\in L}{\sum}\omega(p_g) = 1 \Rightarrow \omega(p_g) = \dfrac{1}{|L|}, \forall \ g\in L .$$

We don't have, a priori any further constraint for the values of $\omega (p_g)$, for $g\notin L$. If we impose that, $\omega (p_g)=0$, for $g\notin L$, then the only possible choice is the linear functional $\lambda:(kG)^* \rightarrow k$, which defines the partial action. Therefore, 
$$\begin{array}{rcl}
Z^1_{par}  (\kg ,k ) & = & \{ \omega :kg\rightarrow k, \; \ \; \omega (p_g )=\frac{1}{|L|} , \quad g\in L \} ,\\
B^1_{par}  (\kg ,k ) & = & \{ \lambda \} \\
H^1_{par}  (\kg ,k ) & = & \{ \omega :kg\rightarrow k, \; \ \; \omega (p_g )=\frac{1}{|L|} , \; g\in L , \quad  \omega (p_g )\neq 0 , \; g\notin L \}
\end{array}$$

For $n=2$, first recall that $e_2 = \widetilde{e}_{2,1}*\widetilde{e}_{2,2}$, in which $e_2(h,l) = h\cdot(l\cdot 1_A)$,
$\widetilde{e}_{2,1}(h,l) = (h\cdot 1_A)\epsilon(l)$ and $\widetilde{e}_{2,2}(h,l) = hl\cdot1_A$. Then, for $g, h \in G$, we have
$$\omega(p_g,p_h) = \widetilde{e}_{2,1}*\omega(p_g,p_h) = \underset{l\in G}{\sum}\lambda_{p_l}\omega(p_{l^{-1}g},p_h) = 
\dfrac{1}{|L|}\underset{l\in L}{\sum}\omega(p_{l^{-1}g},p_h) .$$

This leads to an invariance by translation on the left slot, that is 
$\omega(p_{lg},p_h) = \omega(p_g,p_h),$ for any $g,h\in G$ and $l\in L$. On the other hand, 
$$\begin{array}{ccl} \omega(p_g,p_h) & = & \wide{e}_{2,2}*\omega(p_g,p_h) = \underset{l,m\in L}{\sum}\lambda_{p_l p_m}\omega(p_{l^{-1}g},p_{m^{-1}h})\\
  & = & \dfrac{1}{|L|}\underset{l\in L}{\sum} \omega(p_{l^{-1}g},p_{l^{-1}h}) = \dfrac{1}{|L|}\underset{l\in L}{\sum} \omega(p_g,p_{l^{-1}h})
   \end{array}$$
   
This gives an invariance by translation on the right slot, that is, $\omega(p_{g},p_{lh}) = \omega(p_g,p_h)$, 
for any $g,h\in G$ and $l\in L$. As these invariances are independent, we have finally 
$$\omega(p_{lg},p_{mh}) = \omega(p_g,p_h),$$
for any $g,h\in G$ and $l,m\in L$ \cite{ABDP2}. This translation invariance is a useful tool for searching solutions of partial $2$-cocycles in specific cases. Besides the translation invariance, we have the normalization constraint, given by
$$\omega(1,p_g) = \omega(p_g,1) = \lambda_{p_g} \Rightarrow \dfrac{1}{|L|} = \underset{h\in L}{\sum}\omega(p_g,p_h) = 
\underset{h\in L}{\sum}\omega(p_h,p_g), \ \forall \ g\in L .$$

Finally, we have the cocycle condition. For $g,h,i\in L$, we have
\begin{eqnarray*}
& \, & \dfrac{1}{|L|^3}  =  \lambda_{p_g}\lambda_{p_h}\lambda_{p_i} = \delta \omega(p_g,p_h,p_i)\\
& = &\hspace{-0.3cm} \underset{\underset{t,x,y\in L}{l, m, n, r, s,}}{\sum} \lambda_{p_l}\omega(p_r,p_x)\ov{\omega}
   (\underset{\underset{\Rightarrow \ l = rm\inv{s}}{\Rightarrow\ l^{-1}m=\inv{r}s}}{\underbrace{p_{l^{-1}m}p_{\inv{r}s}}},p_{\inv{x}y}) \omega(p_{\inv{m}n},\underset{\underset{\Rightarrow\ \inv{t} = \inv{i}y\inv{s}}
{\Rightarrow s^{-1}t=\inv{y}i}}{\underbrace{p_{\inv{s}t}p_{\inv{y}i}}}) \ov{\omega}(p_{\inv{n}g},p_{\inv{t}h})\\
& = &\hspace{-0.3cm}  \dfrac{1}{|L|}\underset{\underset{x,y\in L}{m, n, r, s,}}{\sum} \omega(p_r,p_x)\ov{\omega}(p_{\inv{r}s},p_{\inv{x}y})
   \omega(p_{\inv{m}n},p_{\inv{y}i})\ov{\omega}(p_{\inv{n}g},p_{\inv{i}y\inv{s}h})\\
& = &\hspace{-0.3cm}  \dfrac{1}{|L|}\underset{m, n, s, y\in L}{\sum} \lambda_{p_{s}}\lambda_{p_{y}}
   \omega(p_{\inv{m}n},p_{\inv{y}i})\ov{\omega}(p_{\inv{n}g},p_{\inv{i}y\inv{s}h})\\
& = &\hspace{-0.3cm}  \dfrac{1}{|L|^2}\underset{n, s, y\in L}{\sum}\left(\underset{m\in L}{\sum} \dfrac{1}{|L|} 
   \omega(p_{\inv{m}n},p_{\inv{y}i})\right)\ov{\omega}(p_{\inv{n}g},p_{\inv{i}y\inv{s}h})\\
 & = &\hspace{-0.3cm}  \dfrac{1}{|L|}\underset{n, s, y\in L}{\sum} 
   \omega(p_{n},p_{\inv{y}i})\ov{\omega}(p_{\inv{n}g},p_{\inv{s}y\inv{i}h})\dfrac{1}{|L|}\\
& = &\hspace{-0.3cm}  \dfrac{1}{|L|}\underset{n, y\in L}{\sum} 
   \omega(p_{n},p_{\inv{y}i})\ov{\omega}(p_{\inv{n}g},p_{y\inv{i}h})\\
& \stackrel{(*)}{=} &  \dfrac{1}{|L|}\underset{n,x \in L}{\sum} 
   \omega(p_{n},p_{x})\ov{\omega}(p_{\inv{n}g},p_{\inv{x}h}) ,
\end{eqnarray*}
in which $(*)$ is taken putting $\inv{y}i = x$ and $\inv{i}y = y\inv{i} = \inv{x}$. Therefore,

$$\dfrac{1}{|L|^2} = \underset{n, x\in L}{\sum} 
   \omega(p_{n},p_{x})\ov{\omega}(p_{\inv{n}g},p_{\inv{x}h})$$
 \end{exmp}

The next example is a specific case of a partial grading of the base field for a fixed group $G$ and a fixed subgroup $L\leq G$ defining the partial action.

\begin{exmp}\label{australian}
Fix $G =\langle a,b\ |\ a^2 = b^2 = e \rangle = \{e, a, b, ab\}$ and $L =\langle a\rangle $, then $|L| = 2$. 
Let us calculate the partial $1$-cocycles in this case. 
The invariance by translations gives us $\omega(p_e) = \omega(p_a) = x$ and $\omega(p_b)=\omega(p_{ab})=y$,
 $\ov{\omega}(p_e) = \ov{\omega}(p_a) = \ov{x}$ e $\ov{\omega}(p_b)=\ov{\omega}(p_{ab})=\ov{y}$.

By the normalization constraint $
\underset{g\in G}{\sum} \omega(p_g) =  1 = \underset{g\in G}{\sum} \ov{\omega}(p_g) ,
$
we have
\begin{equation}\label{eq04} x+y=\frac{1}{2}      \end{equation} 
 and
\begin{equation}\label{eq05} \ov{x}+\ov{y}=\frac{1}{2}\end{equation} 

Moreover, the condition $\omega*\ov{\omega} = e$, which can be written as
$$\underset{h\in G}{\sum} \omega(p_h)\ov{\omega}(p_{\inv{h}g}) =\left\{ \begin{array}{lr} \frac{1}{|L|} & g\in L \\ 0 & g\notin L \end{array} \right. ,$$
gives us two equations,
\begin{eqnarray}\label{eq06} x\ov{x} + y\ov{y} = \frac{1}{4} ,\end{eqnarray}
\begin{eqnarray}\label{eq07} x\ov{y} + y\ov{x} = 0\end{eqnarray}

Finally, the cocycle condition,
\[
\lambda_{p_g}\lambda_{p_h} =  \underset{m\in G}{\sum} \omega(p_m)\ov{\omega}(p_{\inv{m}h})\omega(p_{m\inv{h}g})
\]
gives us
 \begin{itemize}
  \item For $g=h=e$, 
\begin{eqnarray*} & & \omega(p_e)\ov{\omega}(p_e)\omega(p_e) + \omega(p_a)\ov{\omega}(p_a)\omega(p_a)
 +  \omega(p_b)\ov{\omega}(p_b)\omega(p_b)  + \omega(p_{ab})\ov{\omega}(p_{ab})\omega(p_{ab})\hspace{-0.05cm} = \hspace{-0.05cm}\frac{1}{4}
\end{eqnarray*}  
    $$\begin{array}{cl} \Rightarrow & 2x^2\ov{x} + 2y^2\ov{y} = \frac{1}{4} \Rightarrow x^2\ov{x} + y^2\ov{y} = \frac{1}{8} \Rightarrow x^2\ov{x} + xy\ov{y} 
  - xy\ov{y}  + y^2\ov{y} = \frac{1}{8}\\   
  \stackrel{x\ov{y} = -y\ov{x}}{\Rightarrow} &  x^2\ov{x} + xy\ov{y} 
  + yy\ov{x}  + y^2\ov{y} = \frac{1}{8}\\
  \Rightarrow & \frac{1}{8} = x^2\ov{x} + xy\ov{y} 
  + yy\ov{x}  + y^2\ov{y} = (x+y)(x\ov{x} + y\ov{y})\end{array}$$ 
  
This is the product of equations (\ref{eq04}) and (\ref{eq06}), therefore, no new information is added. The same occurs for $g= e$ and $h=a$, $g = a$ and $h=e$, and $g=h=a$.

  \item For $g= e$ e $h=b$,  
  \begin{eqnarray*} & &\omega(p_e)\ov{\omega}(p_b)\omega(p_b) + \omega(p_a)\ov{\omega}(p_{ab})\omega(p_{ab}) 
 +  \omega(p_b)\ov{\omega}(p_e)\omega(p_e)  + \omega(p_{ab})\ov{\omega}(p_{a})\omega(p_{a})\hspace{-0.05cm} =\hspace{-0.05cm} 0
\end{eqnarray*} 
    $$\Rightarrow xy\ov{y} + xy\ov{y} + xy\ov{x} + xy\ov{x}= 0 \Rightarrow xy\ov{y} + yx\ov{x} = 0.$$

As $\ov{x} + \ov{y} = \frac{1}{2}$, then we have $xy=0$. The same condition is obtained for the cases $g = e$ and $h=ab$, $g = a$ and $h=b$, $g = a$ and $h=ab$, $g = b$ and $h=e$, $g = b$ and $h=a$, $g = ab$ and $h=e$, and $g = ab$ and $h=a$.

  \item For $g= b$ e $h=b$,  
\begin{eqnarray*} & & \omega(p_e)\ov{\omega}(p_b)\omega(p_e) + \omega(p_a)\ov{\omega}(p_{ab})\omega(p_{a}) 
 +  \omega(p_b)\ov{\omega}(p_e)\omega(p_b)  + \omega(p_{ab})\ov{\omega}(p_{a})\omega(p_{ab})\hspace{-0.05cm} =\hspace{-0.05cm} 0
\end{eqnarray*}  
  $$\begin{array}{cl}\Rightarrow & x^2\ov{y} + x^2\ov{y} + y^2\ov{x} + y^2\ov{x}= 0 \Rightarrow x^2\ov{y} + y^2\ov{x} = 0\\
     \Rightarrow & x^2\ov{y} + xy\ov{y} - xy\ov{y}  + y^2\ov{x} = 0\\
  \stackrel{xy\ov{x} = -xy\ov{y}}{\Rightarrow} & x^2\ov{y} + xy\ov{y} 
  + yx\ov{x}  + y^2\ov{x} = 0 \\
  \Rightarrow & 0 = x^2\ov{y} + xy\ov{y} 
  + yx\ov{x}  + y^2\ov{x} = (x+y)(x\ov{y} + y\ov{x}) \end{array}$$
  
This equation is the product of (\ref{eq04}) and (\ref{eq07}), therefore, no new information is added. The same occurs if we take $g= b$ and $h=ab$, $g = ab$ and $h=b$, and $g=h=ab$.

Resuming, we have the following equations: 
\begin{eqnarray*}
 x+y  = \frac{1}{2}, \ \ \ \  \ov{x} + \ov{y} =  \frac{1}{2},\ \ \ \ x\ov{x} + y\ov{y}  =  \frac{1}{4},\ \ \ \
x\ov{y} + y\ov{x}  =  0,\ \ \ \
xy  =  0 ,
\end{eqnarray*}
whose unique possible solution is $\omega(p_g) = \lambda_{p_g} =\left\{\begin{array}{cl} \frac{1}{|L|} & ,\ g\in L\\
											0& ,\ g\not\in L
                                                                       
                                                                      \end{array}\right.
$.
\end{itemize}

For $n=2$, first note that $\omega(p_g,p_h) = \omega(p_{ag},p_h) = \omega(p_g,p_{ah}) = \omega(p_{ag},p_{ah})$,
for every $g,h\in G$, then,
\begin{itemize}
 \item $\omega(p_{e},p_{e}) = \omega(p_{a},p_{e}) = \omega(p_{e},p_{a}) = \omega(p_{a},p_{a})$;
 \item $\omega(p_{b},p_{e}) = \omega(p_{ab},p_{e}) = \omega(p_{b},p_{a}) = \omega(p_{ab},p_{a})$;
 \item $\omega(p_{e},p_{b}) = \omega(p_{a},p_{b}) = \omega(p_{e},p_{ab}) = \omega(p_{a},p_{ab})$;
 \item $\omega(p_{b},p_{b}) = \omega(p_{ab},p_{b}) = \omega(p_{b},p_{ab}) = \omega(p_{ab},p_{ab})$.
  \end{itemize}

The normalization constraint gives us $$\underset{h\in G}{\sum}\omega(p_g,p_h) = \lambda_{p_g} = 
\underset{h\in G}{\sum}\omega(p_h,p_g) .$$

Applying the above normalization constraint respectively for $g=e, a, b, ab$, we have
\begin{itemize}
 \item For $g=e$ ($\lambda_{p_e} = 1/2$), 
\begin{eqnarray*} \omg{e}{e} + \omg{e}{a} + \omg{e}{b} + \omg{e}{ab} & = & \frac{1}{2} \\
\omg{e}{e} + \omg{a}{e} + \omg{b}{e} + \omg{ab}{e} & = & \frac{1}{2},
\end{eqnarray*}
 $$\Rightarrow  \omg{e}{e} + \omg{e}{b} = \frac{1}{4}  \ \ \ \text{ and } \ \ \ \omg{e}{b} = \omg{b}{e}$$
 
The same is obtained for $g=a$ ($\lambda_{p_a} = \frac{1}{2}$).

 \item For $g=b$ ($\lambda_{p_b} = 0$), 
\begin{eqnarray*} \omg{b}{e} + \omg{b}{a} + \omg{b}{b} + \omg{b}{ab} & = & 0 \\
 \omg{e}{b} + \omg{a}{b} + \omg{b}{b} + \omg{ab}{b} & = & 0,
 \end{eqnarray*}
 $$\Rightarrow  \omg{b}{b} = - \omg{e}{b} = -\omg{b}{e} $$
The same is obtained for  $g=ab$ ($\lambda_{p_{ab}} = 0$). 
 \end{itemize}
 
 Therefore, the only remaining independent components are $\omega (p_e ,p_e)$ and $\omega (p_e ,p_b )$. Moreover
 \[
 \omega (p_e , p_e )+ \omega (p_e ,p_b )=\frac{1}{4}.
 \]

The cocycle condition  
\[
h_{(1)}\cdot\omega(l_{(1)},m_{(1)})\omega(h_{(2)},l_{(2)}m_{(2)}) = 
 (h_{(1)}\cdot 1)\omega(h_{(2)},l_{(1)})\omega(h_{(3)}l_{(2)},m),
\]
 can be written, in our case, as 
  $$\underset{s\in G}{\sum}\omg{g}{s}\omg{h\inv{s}}{l\inv{s}} = \underset{s\in G}{\sum}\omg{g\inv{s}}{h\inv{s}}\omg{s}{l} .$$

  This condition gives us 64 equations which are, in fact redundant, that is, every $2$-cochain in this case is a $2$-cocycle

Finally, from $\omega*\ov{\omega} = e_2$, taking $x = \omg{e}{e}$ and $y = \ov{\omega}(p_{e},p_e)$, we obtain the equation \cite{ABDP2}
\[
16xy - 3(x+ y) + \dfrac{1}{2} =0.
\]
\end{exmp}

\section{The associated Hopf algebra of a partial action}

In \cite{DK1}, for a partial action of a group $G$ on a commutative algebra $A$, the authors introduced an inverse 
semigroup $\wide{A}$, given by all the invertible elements of all ideals of the form $1_{x_1}\ldots 1_{x_n}A$, for  
$x_1,\ldots,x_n\in G$ and $n\in \mathbb{N}$. Once showed that $\theta_x(1_{x^{-1}}\wide{A}) = 1_x\wide{A}$, in other words, 
$\theta$ restricted to $\wide{A}$ defines a partial action $\wide{\theta}$ of $G$ on $\wide{A}$ such that their cohomologies 
are the same, that is, $H^n_{par}(G,A) \cong H^n_{par}(G,\wide{A})$. 

This construction brings advantages because $\wide{A}$ 
possesses a richer structure than $A$ and then one can study, for example, extension theory by partial group actions 
from a wider perspective, namely, the theory of extensions of inverse semigroups \cite{DK2}. 

In our context, we can have also very similar constructions, allowing us to trade partial actions of a co-commutative Hopf algebra $H$ on a commutative algebra $A$ by a partial action of $H$ on a commutative and co-commutative Hopf algebra $\widetilde{A}$ generating the same cohomology.

In order to proceed the construction of this new Hopf algebra, one has a technical obstruction concerning the invertible elements of the algebra $A$. Indeed, the multiplicative abelian group $A^{\times}$ embeds into the abelian group $C^n_{par}(H,A)$ for each $n\in \mathbb{N}$ by means of the group monomorphisms $\phi_n :A^{\times} \rightarrow C^n_{par}(H,A)$, given by $\phi_n (a)=ae_n$. These morphisms $\phi_n$ are coherent with the coboundary morphisms, that is, for each $n\in \mathbb{N}$, we have $\delta_n \circ \phi_n (a)= \delta_n (ae_n) =ae_{n+1} =\phi_{n+1}(a)$. Therefore, one can construct a new cochain complex which gets rid of these invertible elements and yet defining the same cohomology.

\begin{defi} \label{reducedcochain} Let $H$ be a co-commutative Hopf algebra and $A$ be a commutative partial $H$-module 
algebra. We define, for $n\in \mathbb{N}$ the $n$-th reduced partial cochain group $\widetilde{C}^n_{par} (H,A)$ as the 
quotient abelian group $\frac{C^n_{par}(H,A)}{A^\times}$.
\end{defi}

\begin{prop} The reduced partial cochain complex $\widetilde{C}^\bullet_{par} (H,A)$ generates cohomology groups isomorphic to those relative to the cochain complex $C^\bullet_{par} (H,A)$.
\end{prop}

\noindent \dem Indeed, denote, for each $n\in \mathbb{N}$, the $n$-th reduced cohomology group by $\widetilde{H}^n_{par}(H,A)$ 
and define the map $\psi^n : H^n_{par} (H,A)  \rightarrow  \widetilde{H}^n_{par} (H,A)$  by 
$\psi^n([f])  \mapsto  [fA^\times ]$

One can easily see that this map is well defined, surjective and it is a morphism of abelian groups. The injectivity comes from the fact that given a partial $n$-cochain $f\in C^n_{par} (H,A)$ and $a\in A^\times$, we have
$af\ast f^{-1} =ae_n =\delta_{n-1} (ae_{n-1})$,
then $f$ and $af$ are cohomologous. Therefore the cohomology groups $H^n_{par} (H,A)$ and $\widetilde{H}^n_{par} (H,A)$ are isomorphic.
\find

\begin{rmk} We will denote the reduced $n$-cochains again by $f$ instead of $fA^\times$ in order to make the notation cleaner. It is clear also that at level zero we have 
$\widetilde{C}^0_{par} (H,A)=\{ 1_A \}$.
\end{rmk}

Now, define the algebra $\widetilde{A}$ as the quotient $\widetilde{A}=\frac{\widehat{A}}{\mathcal{I}}$, in which 
$\widehat{A}$ is the free commutative unital algebra 
\[
\widehat{A} = k[ X_{1_A},  X_{f(h^1, \ldots, h^n)}\ |\ n\geq 1,\ h^1,\ldots, h^n \in H,\ f\in \widetilde{C}^n_{par}(H,A) ].
\]
The set of variables runs over the distinct $f\in \widetilde C^n_{par} (H,A)$, that is, if $f$ and $g$ are two partial $n$-cochains such that $f=g$, then, for every $h^1 \otimes \cdots \otimes h^n \in H^{\otimes n}$ we have $X_{f(h^1 , \ldots h^n)}=X_{g(h^1 , \ldots h^n)}$. The ideal $\mathcal{I}$ is taken exactly to recover certain properties from the original algebra $A$ and from the partial partial action of $H$. This ideal is generated by elements of the type
\begin{equation} \label{unity}
X_{1_A} -1;
\end{equation}
\begin{equation} \label{linearity}
X_{f(h^1, \ldots , \sum_i \lambda_i h^j_i , \ldots , h^n )}-\sum_i \lambda_i X_{f(h^1, \ldots , h^j_i , \ldots h^n )},
\text{\ for each\ } f\in \widetilde{C}^n_{par} (H,A),\ \forall n>0; 
\end{equation}
\begin{equation} \label{null}
\sum_i \lambda_i X_{e_{n_1} (h^{1, 1}, \ldots ,h^{1,n_1})} \ldots X_{e_{n_{k_i}} (h^{k_i, 1}, \ldots ,h^{k_i ,n_{k_i}})}\ ,
\end{equation}
for each zero combination $
\sum_i \lambda_i e_{n_1} (h^{1, 1}, \ldots ,h^{1,n_1}) \ldots e_{n_{k_i}} (h^{k_i, 1}, \ldots ,h^{k_i ,n_{k_i}}) =0 \in A;
$
\begin{equation} \label{convolution}
X_{f(h^1_{(1)} , \ldots , h^n_{(1)})}X_{g(h^1_{(2)} , \ldots , h^n_{(2)})} -X_{(f\ast g)(h^1 , \ldots h^n )};
\end{equation}
\begin{equation} \label{partialaction1}
X_{(\hspace{-0.03cm}h\cdot (f_{1}\hspace{-0.03cm} (h^{1,1}\hspace{-0.15cm},\ldots, h^{\hspace{-0.03cm}1,n_1}\hspace{-0.03cm} ) + \ldots +\hspace{-0.03cm} f_m\hspace{-0.03cm}(l^{m,1} \hspace{-0.15cm},\ldots ,l^{m, n_m}\hspace{-0.03cm})\hspace{-0.03cm})\hspace{-0.03cm})}
\hspace{-0.105cm} - \hspace{-0.105cm} \left( \hspace{-0.07cm} X_{(\hspace{-0.03cm}h\cdot f_{1}\hspace{-0.02cm} (h^{1,1}\hspace{-0.15cm},\ldots, h^{1,n_1})\hspace{-0.03cm})}\hspace{-0.1cm} +\hspace{-0.09cm}\ldots\hspace{-0.07cm} + \hspace{-0.1cm} 
X_{(h\cdot f_m\hspace{-0.02cm}(l^{m,1} \hspace{-0.15cm},\ldots ,l^{m, n_m})\hspace{-0.03cm})}\hspace{-0.1cm} \right);
\end{equation}
\begin{equation} \label{partialaction1,5}
X_{(1_H \cdot f(h^1 , \ldots , h^n ))}- X_{f(h^1 , \ldots , h^n )};
\end{equation}
\begin{equation} \label{partialaction2}
 X_{(h\cdot (f_{1}\hspace{-0.03cm} (h^{1,1}\hspace{-0.15cm},\ldots, h^{1,n_1} )\ldots f_m(l^{m,1}\hspace{-0.15cm} ,\ldots ,l^{m, n_m}\hspace{-0.04cm})\hspace{-0.04cm})\hspace{-0.04cm})}  
\hspace{-0.04cm} -\hspace{-0.04cm}  X_{(h_{(1)}\cdot f_{1}\hspace{-0.03cm} (h^{1,1}\hspace{-0.15cm},\ldots, h^{1,n_1} ))} \ldots X_{(h_{(m)}\cdot f_m(l^{m,1} \hspace{-0.15cm},\ldots ,l^{m, n_m}\hspace{-0.04cm})\hspace{-0.04cm})};
\end{equation}
and 
\begin{equation} \label{partialaction3}
X_{(h\cdot (k\cdot f(h^1 , \ldots , h^n )))} -X_{(h_{(1)}\cdot 1_A) } X_{(h_{(2)}k \cdot f(h^1 , \ldots , h^n ))} .
\end{equation}

\begin{rmk} \begin{enumerate} 
\item Note that $\mathcal{I}$ is indeed an ideal of the algebra $\widehat{A}$, for example, one element 
$h\cdot (f_{1} (h^{1,1},\ldots, h^{1,n_1} )\ldots f_m(l^{m,1} ,\ldots ,l^{m, n_m}))$ can be written as
\begin{eqnarray*}
& & E^{n_1 +\cdots n_m} \left( i_{n_1,n_1 +\cdots +n_m}(f_1) \, \ast \,  (\epsilon^{\otimes n_1}\otimes i_{n_2, n_2 +\cdots +n_m} (f_2)) \, \ast \cdots  \right. \\
& & \left. \cdots \ast (\epsilon^{\otimes (n_1 +\cdots +n_{m-1})}\otimes f_m) \right) (h, h^{1,1},\ldots, h^{1,n_1}, \ldots , l^{m,1} ,\ldots ,l^{m, n_m}) ,
\end{eqnarray*}
according to Definition \ref{auxiliares}. 
\item The condition (\ref{unity}) means that the unit of the algebra $A$ will play the role of the unit of the algebra $\widetilde{A}$.
\item The condition (\ref{null}) refers to every linear combination of monomials involving the units of the cochain groups $\widetilde{C}^\bullet_{par} (H,A)$ which vanish in the algebra $A$. Of course, some of these relations are in fact present among elements of the form (\ref{linearity}), but there are other vanishing linear combinations in $A$ involving partial actions of elements of $H$ upon the unit $1_A$ which needed to be ruled out in order to remember the structure of $A$.
\item Casting out elements of the form (\ref{linearity}) is needed to remember that the generators of $\widetilde{A}$ are linear maps between $H^{\otimes n}$ 
and $A$. In particular, in the quotient we have identities of the type
$X_{f(h^1, \ldots , h^i , \ldots , h^n)}=X_{f(h^1, \ldots , h^i_{(1)} , \ldots , h^n )} \epsilon_H(h^i_{(2)})=  \epsilon_H (h^i_{(1)}) X_{f(h^1 , \ldots , h^i_{(2)} , \ldots , h^n )},$
for each $i\in \{ 1, ... , n \}$, for each $f\in \widetilde{C}^n_{par} (H,A)$ for all $n>0$.
\item Casting out elements of the form (\ref{convolution}) is needed in order to make the relations coming from the convolution product between cochains still valid in $\wide{A}$.
\item Finally, we need to mod out elements of the form (\ref{partialaction1}), (\ref{partialaction1,5}), (\ref{partialaction2}) and (\ref{partialaction3}) in order to recover the linearity of the partial action of $H$ on $A$ and the identities coming from axioms (PA1), (PA2) and (PA3).
\end{enumerate}
\end{rmk}

After taking the quotient, as far as it doesn't lead to a misunderstanding, we are going to denote the classes $X_{f(h^1, \ldots , h^n)}+\mathcal{I} \in \widetilde{A}$ simply by $X_{f(h^1, \ldots , h^n)}$ 

Define also the subalgebra of $A$, $E({A}) = \langle h\cdot 1_A \ |\ h\in H \rangle $
and the unit map $\eta: E(A) \rightarrow \wide{A}$ given by
$\eta ((h^1 \cdot 1_A)\ldots (h^n \cdot 1_A )) =X_{e_1 (h^1 )}\ldots X_{e_1 (h^n )}$.

This map is well defined, because among the generators of the ideal $\mathcal{I}$ which defines the algebra $\widetilde{A}$ there are all linear combinations representing null combinations in $A$ involving the units of the cochain complex. Also, by construction it is an algebra morphism (note that $\eta (1_A )=X_{1_A}=1_{\wide{A}} \in \wide{A}$, consequently, $\widetilde{A}$ is a $E(A)$ algebra. Moreover, the unit map is injective. This can be easily seen considering the evaluation map $\widehat{\mbox{ev}}: \widehat{A}\rightarrow A$ which simply associate to each element of $\widehat{a} \in \widehat{A}$ its value at $\widehat{\mbox{ev}}(\widehat{a}) \in A$. It is easy to see that $\widehat{\mbox{ev}} (\mathcal{I})=0$, then one can define a linear map $\mbox{ev}:\widetilde{A}\rightarrow A$ with the same content. Therefore, if $\eta (a)=0$ in $\widetilde{A}$, then $a=\mbox{ev} (\eta (a))=0$. By the injectivity of the unit map, one can identify $E(A)$ with its image in 
$\eta (E(A)) \subseteq \widetilde{A}$.

Every element in the image of $E({A})$ different of $1_{\wide{A}}$ can be written as a combination of variables $X_e$, in which is the image of an idempotent in $\widetilde{C}^n_{par}(H,A)$ for some $n>0$. In fact, what we are going to prove is that one can rewrite an element of the form 
$(h^1 \cdot 1_A)\ldots (h^n \cdot 1_A)$ as a linear combination of images of the idempotent $e_n \in \widetilde{C}^n_{par}(H,A)$. Let us make induction on the number $n$ of factors $h\cdot 1_A$ involved. For $n=1$, we have 
$h\cdot 1_A = e_1(h)$.

Now, suppose that the result is valid for $r\in \mathbb{N}$, $1 \leq r < n$, that is, 
$$(h^1\cdot 1_A)(h^2\cdot 1_A)\ldots (h^r\cdot 1_A) = \sum_{i=1}^s e_r(l_i^1,\ldots,l_i^r) $$  
for some elements $l^j_i\in H$, for $i\in \{ 1,\ldots ,s \}$ and $j\in \{ 1,\ldots r \}$. Take $h^1, \ldots ,h^n \in H$, then 
$$\begin{array}{rcl}
& & (h^1\cdot 1_A)(h^2\cdot 1_A)\ldots (h^n\cdot 1_A)  =  (h^1\cdot 1_A)[(h^2\cdot 1_A)\ldots(h^n\cdot 1_A)]\\
& = & \sum_i (h^1\cdot 1_A)e_{n-1}(l^1_i ,\ldots,l^{n-1}_i )
= \sum_i (h^1\cdot 1_A)(l^1_i \cdot(l^2_i\cdot(\dots\ \cdot(l^{n-1}_i\cdot 1_A)\dots)))
\end{array}$$
$$\begin{array}{rcl}
& = & \sum_i (\sub{h^1}{1}\cdot\hspace{-0.03cm} 1_A)(\sub{h^1}{2}S(\sub{h^1}{3})l^1_i \cdot\hspace{-0.03cm}(l^2_i \cdot\hspace{-0.03cm}(\dots\ \hspace{-0.05cm} \cdot(l^{n-1}_i \cdot 1_A)\dots)\hspace{-0.03cm})\hspace{-0.03cm})\\
& = & \sum_i \sub{h^1}{1}\cdot (S(\sub{h^1}{2})l^1_i \cdot(l^2_i \cdot(\dots\ \cdot(l^{n-1}_i \cdot 1_A)\dots)))\\
& = & \sum_i e_n(\sub{h^1}{1},S(\sub{h}{2})l^1_i , l^2_i , \ldots, l^{n-1}_i ) .
\end{array}$$

This proves our claim. Moreover, for each $n>0$, we have $e_n(h^1,\ldots,h^n) \in E({A})$. Indeed,
$$\begin{array}{rcl}
& & e_n(h^1,\ldots, h^n)  =  h^1\cdot(h^2\cdot(\dots \ \cdot(h^n\cdot 1_A)\dots))\\
& = & (\sub{h^1}{1}\cdot 1_A)(\sub{h^1}{2}h^2\cdot(\dots \ \cdot(h^n\cdot 1_A)\dots))\\
& = & (\sub{h^1}{1}\cdot 1_A)(\sub{h^1}{2}\sub{h^2}{1}\cdot 1_A)(\sub{h^1}{3}\sub{h^2}{2}h^3\cdot(\dots\ \cdot(h^n\cdot 1_A)\dots))\\
& = & (\sub{h^1}{1}\cdot 1_A)(\sub{h^1}{2}\sub{h^2}{1}\cdot 1_A)\ldots(\sub{h^1}{n}\sub{h^2}{n-1}\ldots\sub{h^{n-1}}{2}h^n\cdot 1_A) \in E({A}) .
\end{array}$$

Our construction will enable us to see a richer structure on the algebra $\widetilde{A}$ with the advantage of getting the same cohomology theory as the original algebra $A$.

\begin{thm} \label{auxiliarhopf} Let $H$ be a co-commutative Hopf algebra $H$ and $A$ be commutative partial $H$-module algebra $A$. Then the algebra $\widetilde{A}$ is a commutative and co-commutative Hopf algebra  which is also a partial $H$-module algebra such that for any $n\in \mathbb{N}$ the $n$-cohomology group $H^n_{par} (H, \widetilde{A})$ is isomorphic to the $n$-cohomology $H_{par}^n(H,A)$.
\end{thm}

\noindent \dem We have already shown that $\widetilde{A}$ is a commutative algebra over $E({A})$. For the coalgebra structure, define the map
$\widehat{\Delta}:\widehat{A}\rightarrow \widehat{A}\otimes_{E({A})} \widehat{A}$, given by,
$$\begin{array}{c} 
 \widehat{\Delta} (X_{f_{1}(h^{1,1},\ldots, h^{1,n_1})}\ldots X_{f_m(l^{m,1} ,\ldots ,l^{m, n_m})}) = \\
= X_{f_1(\sub{h^{1,1}}{1},\ldots, \sub{h^{1,n_1}}{1})}\ldots X_{f_m (\sub{l^{m,1}}{1},\ldots,\sub{l^{m,n_m}}{1})}
  \tens X_{f_1(\sub{h^{1,1}}{2},\ldots, \sub{h^{1,n_1}}{2})} \ldots 
X_{f_m(\sub{l^{m,1}}{2},\ldots , \sub{l^{m,n_m}}{2})},
\end{array}$$
for $f_1\in \widetilde{C}^{n_1}(H,A), \ldots, f^m \in \widetilde{C}^{n_m}(H,A)$. And the map $\widehat{\epsilon}: \widehat{A} \rightarrow E({A})$, given by 
\[
 \widehat{\epsilon} (X_{f_{1}(h^{1,1},\ldots, h^{1,n_1} )}\ldots X_{f_m (l^{m,1} ,\ldots ,l^{m, n_m})}) 
= e_{n_1}\hspace{-0.05cm} (h^{1,1}\hspace{-0.05cm},\ldots, h^{1,n_1} ) \ldots\hspace{-0.03cm} e_{n_m}\hspace{-0.05cm} (l^{m,1}\hspace{-0.05cm} ,\ldots ,l^{m, n_m}). 
\]

Finally, we define the antipode $\widehat{S}:\widehat{A}\rightarrow \widehat{A}$ as
$$ \widehat{S}(X_{f_{1} (h^{1,1},\ldots, h^{1,n_1} )}\ldots X_{f_m(l^{m,1} ,\ldots ,l^{m, n_m})}) 
= X_{f^{-1}_{1}(h^{1,1},\ldots, h^{1,n_1} )}\ldots X_{f^{-1}_m (l^{m,1},\ldots ,l^{m, n_m})},$$
for $f_1\in \widetilde{C}^{n_1}(H,A), \ldots, f^m \in \widetilde{C}^{n_m}(H,A)$.

One needs first to show that these maps can be well defined in $\widetilde{A}$ that is, we must verify that $\mathcal{I}$ 
is a Hopf ideal. Most of the verifications are long, but straightforward. Basically, for $\widehat{\epsilon}$, as its image 
lies in $E(A) \subseteq A$, where the relations are valid, then $\widehat{\epsilon} (\mathcal{I})=0$. For $\widehat{S}$, it 
is also easy to see that $\widehat{S}(\mathcal{I})\subseteq \mathcal{I}$. Therefore, one can define algebra maps 
$\epsilon :\widetilde{A}\rightarrow E(A)$ and $S:\widetilde{A} \rightarrow \widetilde{A}$, ($S$ is an algebra map 
because $\widetilde{A}$ is commutative) with the same form.

The most involved ones are the verifications for $\widehat{\Delta}$. For this task, it is convenient to divide the process into two steps. First, we consider the ideal $\mathcal{J} \trianglelefteq \widehat{A}$ generated only by elements of the form (\ref{unity}), (\ref{linearity}) and (\ref{null}). For elements of the form (\ref{unity}) and (\ref{linearity}), it is quite straightforward, now take an element of the form (\ref{null}) that is, a linear combination
$$x= \sum_i \lambda_i X_{e_{n_1} (h^{1, 1}, \ldots ,h^{1,n_1})} \ldots X_{e_{n_{k_i}} (h^{k_i, 1}, \ldots ,h^{k_i ,n_{k_i}})} \in \mathcal{J}$$
such that $\widehat{\mbox{ev}} (x)=0$. Then, we have

$$\begin{array}{rl}
&\hspace{-0.4cm} \widehat{\Delta} (x)\hspace{-0.1cm}  =\hspace{-0.1cm} \underset{i}{\sum}\hspace{-0.02cm} \lambda_i 
X_{\hspace{-0.05cm}e_{n_1}\hspace{-0.065cm} (\hspace{-0.03cm}h^{1, 1}_{(1)}, \ldots ,h^{1,n_1}_{(1)}\hspace{-0.05cm})}\hspace{-0.04cm}
\hspace{-0.04cm}\ldots\hspace{-0.04cm} X_{\hspace{-0.08cm}e_{n_{k_i}}\hspace{-0.075cm} (\hspace{-0.03cm}h^{k_i, 1}_{(1)}\hspace{-0.05cm}, \ldots ,h^{k_i ,n_{k_i}}_{(1)}\hspace{-0.04cm})}
\hspace{-0.15cm} \otimes \hspace{-0.1cm} X_{\hspace{-0.05cm}e_{n_1}\hspace{-0.06cm} (\hspace{-0.02cm}h^{1, 1}_{(2)}\hspace{-0.03cm}, \ldots ,h^{1,n_1}_{(2)}\hspace{-0.03cm})}
 \hspace{-0.04cm}\ldots\hspace{-0.04cm} X_{\hspace{-0.08cm}e_{n_{k_i}}\hspace{-0.075cm} (\hspace{-0.03cm}h^{\hspace{-0.03cm}k_i, 1}_{(2)}\hspace{-0.1cm}, \ldots ,h^{\hspace{-0.09cm}k_i ,n_{k_i}}_{(2)}\hspace{-0.05cm})} \\

 = &\hspace{-0.4cm} \underset{i}{\sum} \lambda_i X_{\hspace{-0.05cm}e_{n_1}\hspace{-0.065cm} (\hspace{-0.03cm}h^{1, 1}_{(1)}, \ldots ,h^{1,n_1}_{(1)}\hspace{-0.05cm})}\hspace{-0.05cm}
 \ldots\hspace{-0.05cm} X_{e_{n_{k_i}} (h^{k_i, 1}_{(1)}, \ldots ,h^{k_i ,n_{k_i}}_{(1)})} 
  X_{e_{n_1} (h^{1, 1}_{(2)}, \ldots ,h^{1,n_1}_{(2)})}\hspace{-0.05cm} \ldots\hspace{-0.05cm} X_{e_{n_{k_i}} (h^{k_i, 1}_{(2)}, \ldots ,h^{k_i ,n_{k_i}}_{(2)})} \hspace{-0.15cm} \otimes \hspace{-0.1cm} 1 \\
= &\hspace{-0.4cm} \underset{i}{\sum} \lambda_i X_{e_{n_1} (h^{1, 1}_{(1)}, \ldots ,h^{1,n_1}_{(1)})} X_{e_{n_1} (h^{1, 1}_{(2)}, \ldots ,h^{1,n_1}_{(2)})} \ldots 
    X_{e_{n_{k_i}} (h^{k_i, 1}_{(1)}, \ldots ,h^{k_i ,n_{k_i}}_{(1)})}  X_{e_{n_{k_i}} (h^{k_i, 1}_{(2)}, \ldots ,h^{k_i ,n_{k_i}}_{(2)})} \hspace{-0.15cm} \otimes \hspace{-0.1cm} 1\\ 
= &\hspace{-0.4cm}\underset{i}{\sum} \lambda_i \left ( \left( X_{e_{n_1} (h^{1, 1}_{(1)}, \ldots ,h^{1,n_1}_{(1)})} X_{e_{n_1} (h^{1, 1}_{(2)}, \ldots ,h^{1,n_1}_{(2)})} -X_{e_{n_1} (h^{1, 1}, \ldots ,h^{1,n_1})} \right) \right. \\ 
&\hspace{-0.4cm} \; X_{e_{n_2} (h^{2, 1}_{(1)}, \ldots ,h^{2,n_2}_{(1)})} X_{e_{n_2} (h^{2, 1}_{(2)}, \ldots ,h^{2,n_2}_{(2)})} \ldots 
  X_{e_{n_{k_i}} (h^{k_i, 1}_{(1)}, \ldots ,h^{k_i ,n_{k_i}}_{(1)})} X_{e_{n_{k_i}} (h^{k_i, 1}_{(2)}, \ldots ,h^{k_i ,n_{k_i}}_{(2)})} \otimes 1 \\
&\hspace{-0.4cm}  + X_{e_{n_1}(h^{1, 1}, \ldots ,h^{1,n_1})}\hspace{-0.1cm} \left(\hspace{-0.05cm} X_{e_{n_2}(h^{2, 1}_{(1)}, \ldots ,h^{2,n_2}_{(1)})} X_{e_{n_2} (h^{2, 1}_{(2)}, \ldots ,h^{2,n_2}_{(2)})}\hspace{-0.08cm}-\hspace{-0.04cm}X_{e_{n_2} (h^{2, 1}, \ldots ,h^{2,n_2})}\hspace{-0.1cm} \right) \ldots \\

&\hspace{-0.4cm}  \ldots \hspace{-0.05cm} X_{\hspace{-0.05cm}e_{n_{k_i}}\hspace{-0.1cm} (h^{k_i, 1}_{(1)}\hspace{-0.03cm}, \ldots ,h^{k_i ,n_{k_i}}_{(1)}\hspace{-0.04cm})}\hspace{-0.08cm} X_{\hspace{-0.05cm}e_{n_{k_i}}\hspace{-0.1cm} 
(h^{k_i, 1}_{(2)}\hspace{-0.03cm}, \ldots ,h^{k_i ,n_{k_i}}_{(2)}\hspace{-0.05cm})} \hspace{-0.15cm} \otimes \hspace{-0.1cm} 1\hspace{-0.05cm} +\hspace{-0.05cm}\ldots 
 \hspace{-0.05cm}+\hspace{-0.1cm}X_{e_{n_1} (h^{1, 1}, \ldots ,h^{1,n_1})}X_{e_{n_2} (h^{2, 1}, \ldots ,h^{2,n_2})}\hspace{-0.07cm}\ldots \\
&\hspace{-0.4cm}  \left. \ldots \left( X_{e_{n_{k_i}} (h^{k_i, 1}_{(1)}, \ldots ,h^{k_i ,n_{k_i}}_{(1)})} X_{e_{n_{k_i}} (h^{k_i, 1}_{(2)}, \ldots ,h^{k_i ,n_{k_i}}_{(2)})} -X_{e_{n_{k_i}} (h^{k_i, 1}, \ldots ,h^{k_i ,n_{k_i}})} \right) \otimes 1 \right) \\
&\hspace{-0.4cm} + \underset{i}{\sum} \lambda_i X_{e_{n_1} (h^{1, 1}, \ldots ,h^{1,n_1})} \ldots X_{e_{n_{k_i}} (h^{k_i, 1}, \ldots ,h^{k_i ,n_{k_i}})} \otimes 1 .
\end{array}$$

Therefore, $\widehat{\Delta}(x) \in \widehat{A} \otimes\mathcal{J} +\mathcal{J} \otimes \widehat{A}$. Then, one can define a new linear map 
\[
\overline{\Delta} :\widehat{A}/\mathcal{J} \rightarrow \left( \widehat{A}/\mathcal{J} \right) \otimes_{E(A)} \left( \widehat{A}/\mathcal{J} \right) \cong \left( \widehat{A}\otimes_{E(A)} \widehat{A} \right) /\left( \widehat{A} \otimes_{E(A)}\mathcal{J} +\mathcal{J} \otimes_{E(A)} \widehat{A} \right) 
\]
with the same form. Recall that in $\widehat{A}/\mathcal{J}$ we have identities of the form
$X_{f(h^1, \ldots ,h^n )}=X_{f(h^1_{(1)}, \ldots , h^n_{(1)})}\epsilon_H (h^1_{(2)})\ldots \epsilon_H (h^n_{(2)})$.

Now define the ideal $\mathcal{I}' \trianglelefteq \widehat{A}/\mathcal{J}$ generated by the elements of the form (\ref{convolution}), (\ref{partialaction1}), (\ref{partialaction1,5}), (\ref{partialaction2}) and (\ref{partialaction3}). Take an element of the form (\ref{convolution}),
$x=X_{f(h^1_{(1)} , \ldots , h^n_{(1)})}X_{g(h^1_{(2)} , \ldots , h^n_{(2)})} -X_{(f\ast g)(h^1 , \ldots h^n )} \in \mathcal{I}'$,
then,

$$\begin{array}{rcl}
&\hspace{-0.35cm} &\hspace{-0.3cm}\overline{\Delta} (x)= X_{f(h^1_{(1)} , \ldots , h^n_{(1)})}X_{g(h^1_{(3)} , \ldots , h^n_{(3)})} \otimes X_{f(h^1_{(2)} , \ldots , h^n_{(2)})}X_{g(h^1_{(4)} , \ldots , h^n_{(4)})}\\
&\hspace{-0.35cm} &\hspace{-0.3cm} -X_{(f\ast g)(h^1_{(1)} , \ldots h^n_{(1)} )} \otimes X_{(f\ast g)(h^1_{(2)} , \ldots h^n_{(2)} )}\\
&\hspace{-0.35cm} = &\hspace{-0.3cm} X_{f(h^1_{(1)} , \ldots , h^n_{(1)})}X_{g(h^1_{(2)} , \ldots , h^n_{(2)})} \otimes X_{f(h^1_{(3)} , \ldots , h^n_{(3)})}X_{g(h^1_{(4)} , \ldots , h^n_{(4)})} \\
&\hspace{-0.35cm} &\hspace{-0.3cm} -X_{(f\ast g)(h^1_{(1)} , \ldots h^n_{(1)} )} \otimes X_{f(h^1_{(2)} , \ldots , h^n_{(2)})}X_{g(h^1_{(3)} , \ldots , h^n_{(3)})} \\
&\hspace{-0.35cm} &\hspace{-0.3cm} +X_{(f\ast g)(h^1_{(1)} , \ldots h^n_{(1)} )} \otimes X_{f(h^1_{(2)} , \ldots , h^n_{(2)})}X_{g(h^1_{(3)} , \ldots , h^n_{(3)})} \\
&\hspace{-0.35cm} &\hspace{-0.3cm} -X_{(f\ast g)(h^1_{(1)} , \ldots h^n_{(1)} )} \otimes X_{(f\ast g)(h^1_{(2)} , \ldots h^n_{(2)} )} \\
&\hspace{-0.35cm} = &\hspace{-0.3cm} X_{f(h^1_{(1)} , \ldots , h^n_{(1)})}X_{g(h^1_{(2)} , \ldots , h^n_{(2)})} \otimes X_{f(h^1_{(3)} , \ldots , h^n_{(3)})}X_{g(h^1_{(4)} , \ldots , h^n_{(4)})} \\
&\hspace{-0.35cm} &\hspace{-0.3cm} -X_{(f\ast g)(h^1_{(1)} , \ldots h^n_{(1)} )} \epsilon_H (h^1_{(2)}) \ldots\epsilon_H (h^n_{(2)})\otimes X_{f(h^1_{(3)} , \ldots , h^n_{(3)})}X_{g(h^1_{(4)} , \ldots , h^n_{(4)})} \\
&\hspace{-0.35cm} &\hspace{-0.3cm} +X_{(f\ast g)(h^1_{(1)} , \ldots h^n_{(1)} )} \otimes \left( X_{f(h^1_{(2)} , \ldots , h^n_{(2)})}X_{g(h^1_{(3)} , \ldots , h^n_{(3)})}  - X_{(f\ast g)(h^1_{(2)} , \ldots h^n_{(2)} )}\right) \\
&\hspace{-0.35cm} = &\hspace{-0.3cm} \left( X_{f(h^1_{(1)} , \ldots , h^n_{(1)})}X_{g(h^1_{(2)} , \ldots , h^n_{(2)})} -X_{(f\ast g)(h^1_{(1)} , \ldots h^n_{(1)} )} \epsilon_H (h^1_{(2)}) \ldots \epsilon_H (h^n_{(2)})\right) \\ 
&\hspace{-0.35cm} &\hspace{-0.3cm} \; \otimes X_{f(h^1_{(3)} , \ldots , h^n_{(3)})}X_{g(h^1_{(4)} , \ldots , h^n_{(4)})} \\
&\hspace{-0.35cm} &\hspace{-0.3cm} +X_{(f\ast g)(h^1_{(1)} , \ldots h^n_{(1)} )} \otimes \left( X_{f(h^1_{(2)} , \ldots , h^n_{(2)})}X_{g(h^1_{(3)} , \ldots , h^n_{(3)})}  - X_{(f\ast g)(h^1_{(2)} , \ldots h^n_{(2)} )}\hspace{-0.1cm}\right)\hspace{-0.05cm} .
\end{array}$$

Therefore $\overline{\Delta} (x)\in \mathcal{I}' \otimes \left( \widehat{A}/\mathcal{J} \right)   + \left( \widehat{A}/\mathcal{J} \right) \otimes \mathcal{I}'$. With similar strategies, one can prove the same for elements of the form (\ref{partialaction1}), (\ref{partialaction1,5}), (\ref{partialaction2}) and (\ref{partialaction3}). Therefore, there exists a well defined algebra map $\Delta :\widetilde{A} \rightarrow \widetilde{A}\otimes_{E(A)} \widetilde{A}$ with the same form on generators.

It is easy to see that $(\widetilde{A}, \mu, \eta, \Delta, \epsilon)$ gives a commutative and co-commutative bialgebra over the base algebra $E(A)$.

Let us verify the antipode axioms, $(I*S)= (S*I)=\eta\circ \epsilon$. Indeed, for $f_1\in \widetilde{C}^{n_1}(H,A), \ldots f^m \in \widetilde{C}^{n_m}(H,A)$ we have 
$$\begin{array}{rcl}
& & \hspace{-0.4cm}(S*I)(X_{f_{1} (h^{1,1},\ldots, h^{1,n_1} )}\ldots X_{f_m(l^{m,1} ,\ldots ,l^{m, n_m})}) \\
&\hspace{-0.45cm} =&\hspace{-0.4cm} \mu(S\tens I)\circ\Delta (X_{f_{1} (h^{1,1},\ldots, h^{1,n_1} )}\ldots X_{f_m(l^{m,1} ,\ldots ,l^{m, n_m})})\\
&\hspace{-0.45cm} =&\hspace{-0.4cm}  S(X_{f_1(\sub{h^{1,1}}{1},\ldots, \sub{h^{1,n_1}}{1})}\hspace{-0.07cm}\ldots X_{f_m(\sub{l^{m,1}}{1},\ldots,\sub{l^{m,n_m}}{1})} X_{f_1(\sub{h^{1,1}}{2},\ldots, \sub{h^{1,n_1}}{2})}\hspace{-0.07cm}\ldots X_{f_m(\sub{l^{m,1}}{2},\ldots,\sub{l^{m,n_m}}{2})}\\
&\hspace{-0.45cm} =&\hspace{-0.4cm}  X_{f^{-1}_1(\sub{h^{1,1}}{1},\ldots, \sub{h^{1,n_1}}{1})}\ldots X_{f^{-1}_m(\sub{l^{m,1}}{1},\ldots,\sub{l^{m,n_m}}{1})} X_{f_1(\sub{h^{1,1}}{2},\ldots, \sub{h^{1,n_1}}{2})}\ldots X_{f_m(\sub{l^{m,1}}{2},\ldots,\sub{l^{m,n_m}}{2})}\\
&\hspace{-0.45cm} =&\hspace{-0.4cm}  X_{(f^{-1}_1*f_1)(h^{1,1},\ldots, h^{1,n_1})}\ldots X_{(f^{-1}_m *f_m)(l^{m,1} ,\ldots ,l^{m, n_m})}\\
&\hspace{-0.45cm} =&\hspace{-0.4cm}  X_{e_{n_1}(h^{1,1},\ldots, h^{1,n_1})} \ldots X_{e_{n_m}(l^{m,1} ,\ldots ,l^{m, n_m})}=
  \eta \circ\epsilon (X_{f_{1} (h^{1,1},\ldots, h^{1,n_1} )}\ldots X_{f_m(l^{m,1} ,\ldots ,l^{m, n_m})})
\end{array}$$

Analogously, we have the equality $(I\ast S)=\eta \circ \epsilon$. Therefore, $\wide{A}$ is a commutative and co-commutative Hopf algebra over $E({A})$. 

One can define a partial action of $H$ on $\wide{A}$, $\bullet: H\tens\wide{A}\rightarrow\wide{A}$. First, define a linear map $\blacktriangleright :H\otimes \widehat{A} \rightarrow \widehat{A}$ given by  
\[
h\hspace{-0.1cm}\blacktriangleright\hspace{-0.1cm}(\hspace{-0.04cm}X_{f_{1}\hspace{-0.03cm} (h^{1,1}\hspace{-0.06cm},\ldots, h^{1,n_1}\hspace{-0.03cm} )}\ldots
X_{f_m(l^{m,1}\hspace{-0.05cm} ,\ldots ,l^{m, n_m}\hspace{-0.03cm})}\hspace{-0.03cm})\hspace{-0.1cm} =\hspace{-0.1cm}  
X_{(\sub{h}{1}\cdot f_{1}\hspace{-0.03cm} (h^{1,1}\hspace{-0.05cm},\ldots, h^{1,n_1}\hspace{-0.03cm})\hspace{-0.03cm})}
\hspace{-0.03cm}\ldots X_{(\hspace{-0.03cm}\sub{h}{m}\cdot f_m(l^{m,1}\hspace{-0.05cm},\ldots ,l^{m, n_m}\hspace{-0.03cm})\hspace{-0.03cm})}.
\]

For each $h\in H$, one can prove that $h\blacktriangleright \mathcal{I} \subseteq \mathcal{I}$. For example, taking an element 
\[
x=X_{f(h^1_{(1)}, \ldots , h^n_{(1)})}X_{g(h^1_{(2)}, \ldots , h^n_{(2)})} -X_{(f\ast g)(h^1 , \ldots h^n)}, 
\]
we have
$$\begin{array}{rcl}
h\blacktriangleright x & = & X_{h_{(1)}\cdot f(h^1_{(1)}, \ldots , h^n_{(1)})}X_{h_{(2)}\cdot g(h^1_{(2)}, \ldots , h^n_{(2)})} -X_{h\cdot (f\ast g)(h^1 , \ldots h^n)} \\
& = & X_{h_{(1)}\cdot f(h^1_{(1)}, \ldots , h^n_{(1)})}X_{h_{(2)}\cdot g(h^1_{(2)}, \ldots , h^n_{(2)})} -X_{h\cdot (f (h^1_{(1)}, \ldots , h^n_{(1)})g(h^1_{(2)}, \ldots , h^n_{(2)}))} \in \mathcal{I}
\end{array}$$

Then, there is a well defined map $\bullet :H\otimes \widetilde{A} \rightarrow \widetilde{A}$, again, given by. 
\[
h\bullet(\hspace{-0.04cm}X_{f_{1} (h^{1,1}\hspace{-0,1cm},\ldots, h^{1,n_1} )}\ldots X_{f_m(l^{m,1}\hspace{-0,1cm} ,\ldots ,l^{m, n_m}\hspace{-0.03cm})}\hspace{-0,05cm})\hspace{-0,1cm} = 
\hspace{-0,1cm}X_{(\sub{h}{1}\cdot f_{1} (h^{1,1}\hspace{-0,1cm},\ldots, h^{1,n_1} ))}\ldots X_{(\sub{h}{m}\cdot f_m(l^{m,1} \hspace{-0,1cm},\ldots ,l^{m, n_m})\hspace{-0,03cm})}.
\]
It is straightforward  to show that $\bullet$ is a partial action of $H$ on $\widetilde{A}$. This follows directly from the fact that $\cdot$ is a partial action of $H$ on $A$. 

Finally, it remains to verify that $A$ and $\widetilde{A}$ generate the same cohomology groups, that is, for any $n\in \mathbb{N}$ we have $H^n_{par} (H,\wide{A}) \cong H^n_{par} (H,A)$. In fact, what we are going to prove is that $H^n_{par} (H,\wide{A}) \cong \widetilde{H}^n_{par} (H,A)$, which implies our result.

First note that, for any $n\in \mathbb{N}$ and  $h^1 \otimes \cdots \otimes  h^n \in H^{\otimes n}$ we have
\[
X_{e_n (h^1, \ldots , h^n )}=X_{h^1 \cdot (\cdots (h^n \cdot 1_A) \cdots )}=h^1 \bullet (\cdots (h^n \bullet X_{1_A})\cdots )=h^1 \bullet (\cdots (h^n \bullet 1_{\widetilde{A}})\cdots ).
\]
For each $n\in \mathbb{N}$, a reduced partial $n$-cochain $f\in \widetilde{C}^n_{par} (H,A)$ generates a partial $n$-cochain $\widetilde{f} \in C^n_{par} (H,\widetilde{A})$ given by $\widetilde{f}(h^1 , \ldots , h^n )=X_{f(h^1 , \ldots , h^n )}$. On the other hand, each $n$-cochain $g\in C^n_{par} (H,\widetilde{A})$, in order to be convolution invertible, must be of the form 
$g (h^1, \ldots h^n ) = X_{\overline{g} (h^1 ,\ldots , h^n )}$, for some $\overline{g} \in \widetilde{C}^n_{par}(H,A)$, for each $h^1 \otimes \cdots \otimes h^n \in H^{\otimes n}$. Therefore, one can define, for each $n\in \mathbb{N}$, two mutually inverses well defined morphisms of abelian groups
$\Phi : \widetilde{H}^n_{par} (H,A)  \rightarrow  H^n_{par} (H,\widetilde{A})$ given by $\Phi([f]) \mapsto  [\widetilde{f}]$
and $\Psi :   H^n_{par} (H,\widetilde{A}) \rightarrow  \widetilde{H}^n_{par} (H,A)$ given by $\Psi([g]) \mapsto [\overline{g}]$.

These maps produced the isomorphism between the cohomology groups $H^n_{par} (H,\widetilde{A})$ and $\widetilde{H}^n_{par} (H,A)$, and consequently between $H^n_{par} (H,\widetilde{A})$ and 
$H^n_{par} (H,{A})$.
\find

\begin{rmk} For the classical case of a global action of a co-commutative Hopf algebra $H$ on a commutative algebra $A$, \cite{Swe}, one can still construct this Hopf algebra $\wide{A}$, and in this case, as $h\cdot 1_A =\epsilon_H (h) 1_A$, the base subalgebra $E({A})$ coincides with the base field. Therefore, the Hopf algebra $\wide{A}$ is a commutative and co-commutative Hopf algebra over $k$ which gives the same classical cohomological theory as $A$. The properties of this Hopf algebra and its role in the classical cohomology theory is still an interesting topic to be explored
\end{rmk}

\section{Twisted partial actions and crossed products}

In reference \cite{ABDP1}, the authors introduced the notion of a twisted partial action of a Hopf algebra $H$ over an algebra $A$ and described the construction of the crossed product by a 2 cocycle. In particular two crossed products are isomorphic if their associated cocycles are related by a linear map which has properties similar to a convolution invertible $2$-coboundary. Nevertheless, we still did not have a cohomology theory underlying those crossed products. In what follows, we shall see that in the case of co-commutative Hopf algebras acting  partially over commutative algebras the crossed products are indeed classified by the second cohomology group defined before.

 \begin{defi} \label{twistedpartial}\cite{ABDP1} Let $H$ be a Hopf algebra and $A$ be a  unital algebra (with unit $1_A$). Let $\cdot: H\otimes A\ri A$ and $\omega: H\otimes H\ri A$ two linear maps. The pair $(\cdot ,\omega)$ is called a twisted partial action of $H$ over $A$ if,
\begin{enumerate}
\item[(TPA1)] $1_H \cdot a =a$, for every $a\in A$.
\item[(TPA2)] $h\cdot (ab)=(h_{(1)}\cdot a)(h_{(2)} \cdot b)$, for every $h\in H$ and $a,b \in A$.
\item[(TPA3)] $(h_{(1)}\cdot (l_{(1)}\cdot a))\omega(h_{(2)},l_{(2)}) = \omega(h_{(1)},l_{(1)})(\sub{h}{2}l_{(2)}\cdot a)$, for every $h,l\in H$ and $a\in A$.
\item[(TPA4)] $\omega(h,l) = \omega(h_{(1)},l_{(1)})(h_{(2)}l_{(2)}\cdot 1_A)$, for every $h,l \in H$.
\end{enumerate}

In this case, we say that $(A, \cdot,\omega)$ is a twisted partial $H$-module algebra..
\end{defi}

\begin{defi} \cite{ABDP1} Let $H$ be a Hopf algebra and $(A, \cdot ,\omega )$ be a twisted partial $H$-module algebra as above. como acima. Define over $A\otimes H$ a multiplication given by 
$$(a\otimes h)(b\otimes l) = \sum a(h_{(1)}\cdot b)\omega(h_{(2)},l_{(1)})\otimes h_{(3)}l_{(2)}$$
for every $a,b \in A$ e $h,l\in H$. We define the partial crossed product as $\underline{A\#_{\omega}H} = (A\otimes H)(1_A\otimes 1_H)$.
\end{defi}

\begin{prop} \cite{ABDP1} Given a Hopf algebra $H$ and $(A, \cdot ,\omega )$ a twisted partial $H$-modules algebra, the partial crossed product $A\#_{\omega}H$ is unital if, and only if, 

\begin{equation}\label{normalizedcocycle}
\omega(h,1_H) = \omega(1_H,h) = h\cdot 1_A, \qquad \forall h\in H.
\end{equation}

Moreover, the crossed product is associative if, and only if
\begin{equation}\label{cocyclecondition}
(h_{(1)}\cdot \omega(l_{(1)},m_{(1)}))\omega (h_{(2)},l_{(2)}m_{(2)}) = \omega (h_{(1)},l_{(1)})\omega (h_{(2)}l_{(2)},m), \quad \forall  h,l,m\in H. 
\end{equation}
\end{prop}

A linear map $\omega :H^{\otimes 2}\rightarrow A$ satisfying (\ref{normalizedcocycle}) and (\ref{cocyclecondition}) of the above Proposition is called a normalized cocycle.

We denote by $a\# h$ the element 
\[
(a\otimes h)(1_A \otimes 1_H) =a(h_{(1)}\cdot 1_A)\otimes h_{(2)} \in \underline{A\#_{\omega}H} .
\]

One can easily deduce that 
\[
a\# h =a(h_{(1)}\cdot 1_A )\# h_{(2)} .
\]

There is an injective algebra morphism $i:A\rightarrow \underline{A\#_{\omega} H}$, given by $i(a)=a\#1$, this enhances the crossed product $\underline{A\#_{\omega} H}$ with a left $A$-module structure. Also one can show that the linear map 
\[
\begin{array}{rccl} \rho: & \underline{A\#_{\omega} H} & \rightarrow & \underline{A\#_{\omega} H} \otimes H \\
\, & a\# h & \mapsto & a\# h_{(1)} \otimes h_{(2)} \end{array}
\]
defines a right $H$-comodule algebra structure on $\underline{A\#_{\omega} H}$. This left $A$-module and right $H$-comodule structures on the partial crossed product will be important in order to relate crossed products with extensions of $A$ by $H$.

\begin{defi} \cite{ABDP1} Let $A = (A,\cdot,\omega)$ be a twisted partial $H$-module algebra. We say that the twisted partial action is symmetric if 
\begin{enumerate}[(i)]
\item The linear maps $\tilde{e}_{1,2} , \tilde{e}_2: H\otimes H \ri A$, given by $\tilde{e}_{1,2} (h,l) = (h\cdot 1_A)\epsilon(l)$ and $\tilde{e}_2(h,l) = hl\cdot1_A$ are central idempotents in the convolution algebra $Hom_k(H\otimes H, A)$;
\item The map $\omega$ satisfies the cocycle condition (\ref{cocyclecondition}) and it is an invertible element in the ideal $\langle \tilde{e}_{1,2} \ast \tilde{e}_2 \rangle \subset Hom(H\otimes H,A)$.
\item For any $h,l\in H$, we have $e_2 (h,l)= (h\cdot (l\cdot 1_A)) = \sum (h_{(1)}\cdot 1_A)(h_{(2)}l\cdot1_A)= (\tilde{e}_{1,2} \ast \tilde{e}_2 ) (h.l)$.
\end{enumerate} 

The algebra $A$ is called, in this case, a symmetric twisted partial $H$-module algebra.
\end{defi}

For the case of a co-commutative Hopf algebra $H$ and a commutative algebra $A$, every symmetric twisted partial 
action of $H$ over $A$ is in fact a partial action.

\begin{prop} Let $H$ be a co-commutative Hopf algebra and $A$ be a commutative symmetric twisted partial $H$-module 
algebra, then $A$ is a partial $H$-module algebra. 
\end{prop}
 
\noindent\dem Indeed, from axiom (TPA3) of Definition \ref{twistedpartial}, 
\[
\sum(h_{(1)}\cdot (l_{(1)}\cdot a))\omega (h_{(2)},l_{(2)}) = \sum\omega (h_{(1)},l_{(1)}) (h_{(2)}l_{(2)}\cdot a),
\]
we conclude that
$$\begin{array}{ccl} \sum(h\cdot (l\cdot a)) & = & \sum\omega(h_{(1)},l_{(1)})(h_{(2)}l_{(2)}\cdot a)\inv{\omega}(h_{(3)},l_{(3)})\\
    & = & \sum\omega(h_{(1)},l_{(1)})\inv{\omega}(h_{(2)},l_{(2)})(h_{(3)}l_{(3)}\cdot a)\\
    & = & \sum (h_{(1)}\cdot(l_{(1)}\cdot 1_A))(h_{(2)}l_{(2)}\cdot a)\\
    & = & \sum (h_{(1)}\cdot 1_A)(h_{(2)}l_{(1)}\cdot 1_A)(h_{(2)}l_{(2)}\cdot a)\\
    & = & \sum (h_{(1)}\cdot 1_A)(h_{(2)}l\cdot a)
\end{array}$$

By the commutativity of the convolution algebra $Hom(H\otimes H,A)$, we also conclude that $h\cdot(l\cdot a) = (h_{(1)}l\cdot a)(h_{(2)}\cdot1_A)$. Therefore, $A$ is a partial $H$-module algebra.
\find

In the case of $H$ being a co-commutative Hopf algebra and $A$ being a partial $H$ module algebra, we still can define a partial crossed product for each $2$-cocycle $\omega \in Z^2_{par} (H, A)$. In fact, all possible partial crossed products which can be constructed in this case are classified by the second cohomology $H^2_{par} (H, A)$

\subsection{Partial crossed products and $H^2_{par} (H,A)$}

Theorem 4.1 of reference \cite{ABDP1} gives a necessary and sufficient condition on two different symmetric partial actions of a Hopf algebra $H$ over an algebra $A$ in order to get the associated crossed products isomorphic.

\begin{thm} \label{coboundaryisomorphism} \cite{ABDP1} Let $A$ be a unital algebra and $H$ a Hopf algebra with two symmetric twisted partial actions, $h\otimes a \mapsto h\cdot a$ and $h\otimes a \mapsto h\bullet a$, with cocycles $\omega$ and $\sigma$, respectively. Suppose that there is an isomorphism 
\[
\Phi : \underline{A\#_{\omega} H} \rightarrow \underline{A\#_{\sigma} H}
\]
which is also a left $A$-module and a right $H$-comodule map. Then there exists linear maps $u,v\in \mbox{Hom}_k (H,A)$ such that for all $h,k\in H$, $a\in A$
\begin{enumerate}[(i)]
\item $u\ast v (h)=h\cdot 1_A$;
\item $u(h)=u(h_{(1)})(h_{(2)}\cdot 1_A)=(h_{(1)}\cdot 1_A)u(h_{(2)})$;
\item $h\bullet a =v(h_{(1)})(h_{(2)}\cdot a)u(h_{(3)})$
\item $\sigma (h,k)=v(h_{(1)})(h_{(2)}\cdot v(h_{(2)})) \omega (h_{(3)} ,k_{(2)})u(h_{(4)}k_{(3)})$;
\item $\Phi (a\#_{\omega} h) =au(h_{(1)})\#_{\sigma} h_{(2)}$.
\end{enumerate}

Conversely, given maps $u,v\in \mbox{Hom}_k (H,A)$, satisfying (i),(ii), (iii) and (iv) and, in addition $u(1_H )=v(1_H )=1_A$, then the map $\Phi$, as presented in (v), is an isomorphism of algebras.
\end{thm}

For the case of a cocommutative Hopf algebra $H$ and a commutative algebra $A$, itens (i) and (ii) imply that $u$ is the convolution inverse of $v$ in the ideal $e_1 \ast \mbox{Hom}_k (H,A)$. Item (iii), in its turn, implies that the two  partial actions $\bullet$ and $\cdot$ are equal. At last, but not least, item (iv) can be rewritten as
\[
\sigma \ast \omega^{-1} (h,k) =(h_{(1)} \cdot v(k_{(1)}))u(h_{(2)}k_{(2)})v(h_{(3)})=\delta_1 (v)(h,k) .
\]

Therefore, one can rewrite Theorem \ref{coboundaryisomorphism} as:

\begin{thm} Let $H$ be a co-commutative Hopf algebra and $A$ be a partial $H$ module algebra. Then, given two 
partial $2$-cocycles $\omega , \sigma \in Z^2_{par} (H, A)$, the associated partial crossed products 
$\underline{A\#_{\omega} H}$ and $\underline{A\#_{\sigma} H}$ are isomorphic if, and only if, $\omega$ and $\sigma$ 
are cohomologous, that is, they belong to the same class in the cohomology group $H^2_{par} (H,A)$.
\end{thm}

In order to conclude that the second partial cohomology fully classify all the isomorphism classes of partial crossed products, it remains to check that every class in $H^2_{par} (H,A)$ contains a normalized two cocycle.

\begin{prop} Given a partial $2$-cocycle $\omega\in Z^2_{par}(H,A)$, there exists a normalized  $2$-cocycle $\wide{\omega}\in Z^2(H,A)$, which is cohomologous to $\omega$.
\end{prop}

\noindent \dem Indeed, take $\omega$ a $2$-cocycle, then $\omega$ satisfies
$$(h_{(1)}\cdot \omega(k_{(1)},l_{(1)}))\omega(h_{(2)},k_{(2)}l_{(2)}) = \omega(h_{(1)},k_{(1)})\omega(h_{(2)}k_{(2)},l) .$$
Putting $h=1_H$ in the expression above , we have 
$$\begin{array}{rcl}& & (1_H\cdot \omega(k_{(1)},l_{(1)}))\omega(1_H,k_{(2)}l_{(2)}) = \omega(1_H,k_{(1)})\omega(1_Hk_{(2)},l)\\
  & \Rightarrow & \omega(k_{(1)},l_{(1)})\omega(1_H,k_{(2)}l_{(2)}) = \omega(1_H,k_{(1)})\omega(k_{(2)},l)\\
  & \Rightarrow & \ov{\omega}(k_{(1)},l_{(1)})\omega(k_{(2)},l_{(2)})\omega(1_H,k_{(3)}l_{(3)})
                = \omega(1_H,k_{(1)})\omega(k_{(2)},l_{(1)})\ov{\omega}(k_{(3)},l_{(2)})\\
  & \Rightarrow & k_{(1)}\cdot(l_{(1)}\cdot 1_A)\omega(1_H,k_{(2)}l_{(2)}) = \omega(1_H,k_{(1)})(k_{(2)}\cdot (l_{(1)}\cdot 1_A))\\
  & \stackrel{k  = 1_H}{\Rightarrow} & (l_{(1)}\cdot 1_A)\omega(1_H,l_{(2)}) = \omega(1_H,1_H)\cdot (l\cdot 1_A)\\
   \end{array}$$

As $\inv{\omega}(1_H,1_H)\omega(1_H,1_H) = 1_H\cdot 1_A = 1_A$, we conclude that $\omega(1_H,1_H)\in A^{\times}$. 
Then, one can define $\wide{\omega}(h,k) = \omega(h,k)\inv{(\omega(1_H ,1_H ))}$. It is easy to see that 
$\wide{\omega}(1_H,l) = (l\cdot 1_A)$.

On the other hand, putting $l = 1_H$ in the $2$-cocycle condition, we have 
$$\begin{array}{cl}\hspace{-0.4cm} &\hspace{-0.6cm} (h_{(1)}\cdot \omega(k_{(1)},1_H))\omega(h_{(2)},k_{(2)}) = \omega(h_{(1)},k_{(1)})\omega(h_{(2)}k_{(2)},1_H)\\
  \hspace{-0.4cm} \Rightarrow &\hspace{-0.6cm} h_{(1)}\cdot \omega(k_{(1)},1_H))\omega(h_{(2)},k_{(2)})\ov{\omega}(h_{(3)},k_{(3)}) 
               \hspace{-0.05cm} = \hspace{-0.05cm} \ov{\omega}(h_{(1)},k_{(1)})\omega(h_{(2)},k_{(2)})\omega(h_{(3)}k_{(3)},1_H)\\
  \hspace{-0.4cm} \Rightarrow &\hspace{-0.6cm} h_{(1)}\cdot \omega(k_{(1)},1_H))(h_{(2)}\cdot(k_{(2)}\cdot 1_A)) 
                = (h_{(1)}\cdot(k_{(1)}\cdot 1_A))\omega(h_{(2)}k_{(2)},1_H)\\
  \hspace{-0.4cm} \Rightarrow &\hspace{-0.6cm} h_{(1)}\cdot \omega(k_{(1)},1_H)(k_{(2)}\cdot 1_A)) 
                = (h_{(1)}\cdot1_A)(h_{(2)}k_{(1)}\cdot 1_A)\omega(h_{(3)}k_{(2)},1_H)\\
  \hspace{-0.4cm} \Rightarrow &\hspace{-0.6cm} h\cdot \omega(k,1_H) = (h_{(1)}\cdot1_A)\omega(h_{(2)}k,1_H)\\
  \hspace{-0.1cm} \stackrel{k=1_H}{\Rightarrow} &\hspace{-0.4cm} h\cdot \omega(1_H,1_H) = (h_{(1)}\cdot1_A)\omega(h_{(2)},1_H)\\
  \hspace{-0.4cm} \Rightarrow &\hspace{-0.6cm} h\cdot \omega(1_H,1_H) = \omega(h,1_H) .
\end{array}$$

Therefore, $ h\cdot \wide{\omega}(1_H,1_H) = h\cdot (1_H \cdot 1_A)=h\cdot 1_A = \wide{\omega}(h,1_H)$.

Finally, let us verify that $\wide{\omega}$ is cohomologous to $\omega$, that is, there exists $\phi \in C^1_{par}(H,A)$ such that $\wide{\omega}\ast \inv{\omega} = \delta_1 \phi$. 
Indeed, on one hand, note that  
\[
\wide{\omega}(h_{(1)},k_{(1)})\inv{\omega}(h_{(2)},k_{(2)}) = (h\cdot(k\cdot1_A))\inv{(\omega(1_H,1_H))}.
\]

On the other hand, 
\[
\delta\phi(h,k) = (h_{(1)}\cdot \phi(k_{(1)}))\inv{\phi}(h_{(2)}k_{(2)})\phi(h_{(3)}).
\]

Then, if we define  $\phi(k) = (k\cdot 1_A)\inv{(\omega(1_H,1_H))}$, we have
$$\begin{array}{rcl}\delta\phi(h,k) 
  &= & h_{(1)}\cdot((k_{(1)}\cdot 1_A)\inv{(\omega(1_H,1_H))})(h_{(2)}k_{(2)}\cdot 1_A)(\omega(1_H,1_H))\\
  &  & \;    (h_{(3)}\cdot 1_A)\inv{(\omega(1_H,1_H))}\\
  &=& (h_{(1)}\cdot(k_{(1)}\cdot 1_A))\inv{(\omega(1_H,1_H))}(h_{(2)}k_{(2)}\cdot 1_A)(h_{(3)}\cdot 1_A)\\
  &=& (h_{(1)}\cdot(k_{(1)}\cdot 1_A))\inv{(\omega(1_H,1_H))} 
         \end{array}$$
This concludes our proof.
\find

\subsection{The Hopf algebroid structure of the partial crossed product}

We saw that the cohomology theory for a co-commutative Hopf algebra $H$ acting partially over a commutative algebra $A$ is equivalent to the cohomology theory of the same Hopf algebra $H$ acting on a commutative and co-commutative Hopf algebra $\wide{A}$ whose base ring is the commutative algebra $E({A})$. This replacement gives us a deeper understanding about the structure of crossed products. In fact, we shall see that the crossed product has a structure of a Hopf algebroid over the base algebra $E({A})$. Let us recall briefly the definition of a Hopf algebroid, for a detailed presentation, se the reference \cite{B}.

\begin{defi} \cite{B} Given a $k$ algebra $A$, a left (resp. right) bialgebroid over $A$ is given by the data $(\mathcal{H}, A, s_l ,t_l ,\Delta_l , \epsilon_l )$ (resp. $(\mathcal{H}, A, s_r , t_r , \Delta_r , \epsilon_r )$) such that:

\begin{enumerate}
\item $\mathcal{H}$ is a $k$ algebra.
\item The map $s_l$ (resp. $s_r$) is a morphism of algebras between $A$ and $\mathcal{H}$, and the map $t_l$ (resp. $t_r$) is an anti-morphism of algebras between $A$ and $\mathcal{H}$. Their images commute, that is, for every $a,b\in A$ we have $s_l (a)t_l(b)=t_l(b)s_l(a)$ (resp. $s_r(a)t_r(b)=t_r(b)s_r(a)$).
By the maps $s_l,t_l$ (resp. $s_r , t_r$) the algebra $\mathcal{H}$ inherits a structure of $A$ bimodule given by $a\triangleright h \triangleleft b =s_l(a)t_l(b)h$ (resp. $a\blacktriangleright h \blacktriangleleft b =h s_r(b)t_r(a)$).
\item The triple $(\mathcal{H},\Delta_l , \epsilon_l )$ (resp. $(\mathcal{H}, \Delta_r , \epsilon_r )$) is an $A$ coring relative to the structure of $A$ bimodule defined by $s_l$ and $t_l$ (resp. $s_r$, and $t_r$).
\item The image of $\Delta_l$ (resp. $\Delta_r $) lies in the Takeuchi subalgebra
\[
\mathcal{H}  {}_A\times \mathcal{H} =\{ \sum_i h_i \otimes k_i \in \mathcal{H}\otimes_{A, \triangleright \, \triangleleft} \mathcal{H} \, |\, \sum_i h_i t_l(a) \otimes k_i =\sum_i h_i \otimes k_i s_l(a) \, \; \forall a\in A \} ,
\]
respectively,
\[
\mathcal{H} \times_A \mathcal{H} =\{ \sum_i h_i \otimes k_i \in \mathcal{H}\otimes_{A, \blacktriangleright \, \blacktriangleleft} \mathcal{H} \, |\, \sum_i s_r (a) h_i  \otimes k_i =\sum_i h_i \otimes t_r (a) k_i \, \; \forall a\in A \} ,
\]
and it is an algebra morphism.
\item For every $h,k\in \mathcal{H}$, we have $\epsilon_l (hk)=\epsilon_l (hs_l(\epsilon_l (k)))=\epsilon_l (ht_l(\epsilon_l (k)))$,
respectively, $\epsilon_r (hk) =\epsilon_r (s_r(\epsilon_r (h))k)=\epsilon_r (t_r (\epsilon_r (h))k)$.

\end{enumerate}

Given two anti-isomorphic algebras $A_l$ and $A_r$ (ie, $A_l\cong A_r^{op}$), a left $A_l$ bialgebroid $(\mathcal{H}, A_l, s_l,t_l,\Delta_l , \epsilon_l )$ and a right $A_r$ bialgebroid $(\mathcal{H}, A_r, s_r, t_r, \Delta_r , \epsilon_r )$, a Hopf algebroid structure on $\mathcal{H}$ is given by an antipode, that is, an algebra anti-homomorphism $\mathcal{S}:\mathcal{H}\rightarrow \mathcal{H}$ such that
\begin{enumerate}[(i)]
\item $s_l \circ \epsilon_l \circ t_r =t_r$, $t_l\circ \epsilon_l \circ s_r =s_r$, $s_r\circ \epsilon_r \circ t_l =t_l$ and $t_r\circ \epsilon_r \circ s_l =s_l$;
\item $(\Delta_l \otimes_{A_r} I)\circ \Delta_r =(I\otimes_{A_l} \Delta_r )\circ \Delta_l$ and
 $(I\otimes_{A_r} \Delta_l )\circ \Delta_r =(\Delta_r \otimes_{A_l} I )\circ \Delta_l$; 
\item $\mathcal{S}(t_l(a)ht_r(b'))=s_r(b')\mathcal{S}(h) s_l(a)$,
for all $a\in A_l$, $b'\in A_r$ and $h\in \mathcal{H}$;
\item $\mu_{\mathcal{H}} \circ (\mathcal{S} \otimes_{A_l} I)\circ \Delta_l =s_r \circ \epsilon_r$ and $\mu_{\mathcal{H}} \circ (I\otimes_{A_r} S)\circ \Delta_r =s\circ \epsilon_l$.
\end{enumerate} 
\end{defi}

In our case, both algebras, $A_l$ and $A_r$, coincide with the commutative algebra $E({A})$ and the crossed product $\underline{\wide{A}\#_{\omega} H}$ will play the role of the Hopf algebroid $\mathcal{H}$ of the previous definition.

\begin{thm} \label{crossedproductishopfalgebroid} Let $H$ be a co-commutative Hopf algebra and $A$ be a commutative partial $H$-module algebra. Consider the commutative and co-commutative Hopf algebra $\wide{A}$, constructed in Theorem \ref{auxiliarhopf}, over the commutative algebra $E({A})$, which is also a partial $H$-module algebra. Then, the crossed product $\underline{\wide{A}\#_{\omega} H}$, in which $\omega$ is a partial $2$-cocycle in 
$H^2_{par} (H,A)=H^2_{par} (H, \wide{A})$, is a Hopf algebroid over the base algebra $E({A})$.
\end{thm}

\noindent \dem The source and target maps, both left and right, are defined by the restriction to $E({A})$ of the canonical inclusion of $\wide{A}$ into the crossed product
$$\func{s_l, t_l, s_r, t_r:}{E({A})}{\underline{\wide{A}\#_{\omega} H}}{a}{a\#1_H}.$$

We have already seen that this inclusion is an algebra map and by the commutativity of $\wide{A}$, the images of source and 
target maps commute among themselves. Note that, even though the left and right sources and targets are equal, their associated 
bimodule structures are different nonetheless. Indeed, for $a,a' \in E({A})$ and $b\# h \in \underline{\wide{A}\#_{\omega} H}$
we have

\[
a\triangleright (b\# h) \triangleleft a' =(a\# 1_H)(a'\# 1_H)(b\# h)=aa'b\# h ,
\]
and
\[
a\blacktriangleright (b\# h) \blacktriangleleft a' =(b\# h)(a\# 1_H)(a'\# 1_H)=b(h_{(1)}\cdot aa')\# h_{(2)},
\]

The left and right comultiplication maps are defined, respectively, as
$$\func{\Delta_l:}{\wide{A}\#_{\omega}H}{\wide{A}\#_{\omega}H \otimes \wide{A}\#_{\omega}H}{a\#h}{a_{(1)}\#h_{(1)}\otimes_{E({A}),\triangleright \, \triangleleft} a_{(2)}\#h_{(2)}}$$
and 
$$\func{\Delta_r:}{\wide{A}\#H}{\wide{A}\#H \otimes \wide{A}\#H}{a\#h}{a_{(1)}\#h_{(1)}\otimes_{E({A}),\blacktriangleright \, \blacktriangleleft} a_{(2)}\#h_{(2)}} ,$$
in which the tensor product $\otimes_{E({A}),\triangleright \, \triangleleft}$ (resp. $\otimes_{E({A}),\blacktriangleright \, \blacktriangleleft}$) is balanced with respect to the  $E({A})$-bimodule structure implemented by $s_l ,t_l$ (resp. $s_r, t_r$). It is easy to see that 
\begin{eqnarray*} 
(\Delta_l \otimes_{E({A}),\blacktriangleright \, \blacktriangleleft} I)\circ \Delta_r &=&(I\otimes_{E({A}),\triangleright \, \triangleleft} \Delta_r )\circ \Delta_l \\ 
(I\otimes_{E({A}),\blacktriangleright \, \blacktriangleleft} \Delta_l )\circ \Delta_r &=&(\Delta_r \otimes_{E({A}),\triangleright \, \triangleleft} I )\circ \Delta_l . 
\end{eqnarray*}

One can see also that, for any $a\in E({A})$, 
\[
\Delta_l (s_l (a))=s_l (a) \otimes (1_A \# 1_H ), \quad \mbox{ and } \quad \Delta_l (t_l (a))= (1_A \# 1_H ) \otimes t_l (a) . 
\]

This is because that any element of $E({A})$ is a linear combination of monomials of the form $(h^1 \cdot 1_A )\ldots (h^n \cdot 1_A )$, for $h^1 ,\ldots , h^n \in H$, and then 
\begin{eqnarray*}
& & \Delta_l (s_l ((h^1 \cdot 1_A )\ldots (h^n \cdot 1_A ))) \\
& = & (h_{(1)}^1 \cdot 1_A )\ldots (h_{(1)}^n \cdot 1_A ) \# 1_H \otimes (h_{(2)}^1 \cdot 1_A )\ldots (h_{(2)}^n \cdot 1_A ) \# 1_H\\
& = & (h_{(1)}^1 \cdot 1_A )\ldots (h_{(1)}^n \cdot 1_A ) \# 1_H \otimes ((h_{(2)}^1 \cdot 1_A )\ldots (h_{(2)}^n \cdot 1_A ) \# 1_H )(1_A \# 1_H )\\
& = & (h_{(1)}^1 \cdot 1_A )\ldots (h_{(1)}^n \cdot 1_A ) \# 1_H \otimes s_l ((h_{(2)}^1 \cdot 1_A )\ldots (h_{(2)}^n \cdot 1_A ))(1_A \# 1_H )\\
& = & (h_{(1)}^1 \cdot 1_A )\ldots (h_{(1)}^n \cdot 1_A ) \# 1_H \otimes ((h_{(2)}^1 \cdot 1_A )\ldots (h_{(2)}^n \cdot 1_A ))\triangleright (1_A \# 1_H )\\
& = & ((h_{(1)}^1 \cdot 1_A )\ldots (h_{(1)}^n \cdot 1_A ) \# 1_H) \triangleleft ((h_{(2)}^1 \cdot 1_A )\ldots (h_{(2)}^n \cdot 1_A ))\otimes 1_A \# 1_H \\
& = & t_l ((h_{(2)}^1 \cdot 1_A )\ldots (h_{(2)}^n \cdot 1_A ))((h_{(1)}^1 \cdot 1_A )\ldots (h_{(1)}^n \cdot 1_A ) \# 1_H) \otimes 1_A \# 1_H \\
& = & ((h_{(2)}^1 \cdot 1_A )\ldots (h_{(2)}^n \cdot 1_A ) \# 1_H )((h_{(1)}^1 \cdot 1_A )\ldots (h_{(1)}^n \cdot 1_A ) \# 1_H) \otimes 1_A \# 1_H \\
& = & ((h_{(1)}^1 \cdot 1_A )\ldots (h_{(1)}^n \cdot 1_A )(h_{(2)}^1 \cdot 1_A )\ldots (h_{(2)}^n \cdot 1_A ) \# 1_H) \otimes 1_A \# 1_H \\
& = & (h^1 \cdot 1_A )\ldots (h^n \cdot 1_A ) \# 1_H \otimes 1_A \# 1_H \\
& = & s_l ((h^1 \cdot 1_A )\ldots (h^n \cdot 1_A )) \otimes 1_A \# 1_H.
\end{eqnarray*}

From this we deduce that $\Delta_l$ is a morphism of $E({A})$-bimodules with the bimodule structure given by $\triangleright$ and $\triangleleft$.

The same occurs for $\Delta_r$, that is, for any $a\in E({A})$, we have
\[
\Delta_r (t_r (a))=t_r (a) \otimes (1_A \# 1_H ), \quad \mbox{ and } \quad \Delta_r (s_r (a))= (1_A \# 1_H ) \otimes s_r (a) , 
\]
and then $\Delta_r$ is a morphism of $E({A})$-bimodules, with the bimodule structure given by $\blacktriangleright$ and $\blacktriangleleft$.

The\ \ image\ \ of\ \ the\ \ left\ \ comultiplication\ \ $\Delta_l$ \ \ lies\ \ in\ \ the\ \ left\ \ Takeuchi\ \ product\\ $\underline{\wide{A}\#_{\omega} H} {}_{E({A})}\times \underline{\wide{A}\#_{\omega} H}$. 
Indeed, for $a\# h \in \underline{\wide{A}\#_{\omega} H}$ and $b\in E({A})$, we have
\begin{eqnarray*}
& \, &\hspace{-0.45cm} (a_{(1)}\# h_{(1)})t_l(b)\otimes a_{(2)}\#h_{(2)}= (a_{(1)}\#h_{(1)})(b\# 1_H)\otimes a_{(2)}\#h_{(2)}\\
   &\hspace{-0.65cm} = &\hspace{-0.45cm} a_{(1)}(h_{(1)}\cdot b)\omega(h_{(2)},1_H)\# h_{(3)} \otimes a_{(2)}\#h_{(4)}\hspace{-0.05cm} = \hspace{-0.05cm} a_{(1)}(h_{(1)}\cdot b)(h_{(2)}\cdot 1_A)\# h_{(3)}\otimes a_{(2)}\#h_{(4)}\\
   &\hspace{-0.65cm} = &\hspace{-0.45cm} (h_{(1)}\cdot b)a_{(1)}\# h_{(2)}\otimes a_{(2)}\#h_{(3)} =  ((h_{(1)}\cdot b)\# 1_H )(a_{(1)} \# h_{(2)}) \otimes a_{(2)} \# h_{(3)}\\
   &\hspace{-0.65cm} = &\hspace{-0.45cm} t_l (h_{(1)}\cdot b)(a_{(1)} \# h_{(2)}) \otimes a_{(2)} \# h_{(3)}  =  a_{(1)} \# h_{(2)} \otimes s_l (h_{(1)}\cdot b)(a_{(2)} \# h_{(3)})\\
   &\hspace{-0.65cm} = &\hspace{-0.45cm} a_{(1)} \# h_{(2)} \otimes ((h_{(1)}\cdot b)\# 1_H)(a_{(2)} \# h_{(3)}) =  a_{(1)}\# h_{(1)}\otimes a_{(2)}(h_{(2)}\cdot b)\#h_{(3)}\\ 
   &\hspace{-0.65cm} = &\hspace{-0.45cm} a_{(1)}\# h_{(1)}\otimes (a_{(2)}\#h_{(2)})(b\# 1_H) =  a_{(1)}\# h_{(1)}\otimes (a_{(2)}\#h_{(2)})s_l(b).
\end{eqnarray*}

Moreover, $\Delta_l$ is a morphism of algebras. Take any $a\#h,\ b\#l \in \underline{\wide{A}\#_{\omega} H}$, then
\begin{eqnarray*}
& & \Delta_l((a\#h)(b\#l)) = \Delta_l (a(h_{(1)}\cdot b)\omega(h_{(2)},l_{(1)})\#h_{(3)}l_{(2)})\\
   & = & (a(h_{(1)}\cdot b)\omega(h_{(2)},l_{(1)}))_{(1)}\#(h_{(3)}l_{(2)})_{(1)}
    \otimes (a(h_{(1)}\cdot b)\omega(h_{(2)},l_{(1)}))_{(2)}\#(h_{(3)}l_{(2)})_{(2)}\\
   & = & a_{(1)}(h_{(1)}\cdot b_{(1)})\omega(h_{(3)},l_{(1)})\# h_{(5)}l_{(3)}\otimes  a_{(2)}(h_{(2)}\cdot b_{(2)})\omega(h_{(4)},l_{(2)})\#h_{(6)}l_{(4)}\\
\end{eqnarray*} 

On the other hand, 
\begin{eqnarray*}
& & \dl(a\#h)\dl(b\#l) =  (a_{(1)}\# h_{(1)}\otimes a_{(2)}\# h_{(2)})(b_{(1)}\# l_{(1)}\otimes b_{(2)}\# l_{(2)})\\
   & = & (a_{(1)}\# h_{(1)})(b_{(1)}\#l_{(1)})\otimes (a_{(2)}\# h_{(2)})(b_{(2)}\#l_{(2)})\\
   & = & a_{(1)}(h_{(1)}\cdot b_{(1)})\omega(h_{(2)},l_{(1)})\# h_{(3)}l_{(2)}
    \otimes a_{(2)}(h_{(4)}\cdot b_{(2)})\omega(h_{(5)},l_{(3)})\# h_{(6)}l_{(4)} .
\end{eqnarray*}

The equality follows from the co-commutativity of $H$.

Analogously, one can prove that the image of the right comultiplication lies in the right Takeuchi product $\underline{\wide{A}\#_{\omega} H} \times_{E({A})} \underline{\wide{A}\#_{\omega} H}$, and it is an algebra morphism.

The left and right counits are defined, respectively, as
$$\func{\epsilon_l:}{\wide{A}\#_{\omega}H}{E({A})}{a\#h}{\epsilon_l(a\#h):= \epsilon_{\wide{A}}(a)(h\cdot 1_A)}$$
and
$$\func{\epsilon_r:}{\wide{A}\#_{\omega}H}{E({A})}{a\#h}{\epsilon_r(a\#h):= S(h)\cdot \epsilon_{\wide{A}}(a)} .$$

First, both are morphisms of $E({A})$-bimodules, with their respectives structures. Take $a,a'\in E({A})$ and $b\# h \in \underline{\wide{A}\#_{\omega} H}$, then
\begin{eqnarray*}
\epsilon_l (a\triangleright (b\# h)\triangleleft a') & = & \epsilon_l ((aa'b \# h))
 =  \epsilon_{\wide{A}}(aa')\epsilon_{\wide{A}}(b)(h\cdot 1_A)\\
& = & a\epsilon_{\wide{A}}(b)(h\cdot 1_A) a'
 =  a\epsilon_l (b\# h) a' ,
\end{eqnarray*}
and
\begin{eqnarray*}
& \, & \epsilon_l (a\blacktriangleright (b\# h)\blacktriangleleft a') =  \epsilon_l (b(h_{(1)}\cdot aa') \# h_{(2)})\\ 
& = & S (h_{(2)})\cdot (\epsilon_{\wide{A}}(b)  \epsilon_{\wide{A}}(h_{(1)}\cdot aa')) =  (S (h_{(3)})\cdot \epsilon_{\wide{A}}(b)) (S (h_{(2)})\cdot (h_{(1)}\cdot aa')) \\
& = & (S (h_{(3)})\cdot \epsilon_{\wide{A}}(b)) (S (h_{(1)})h_{(2)}\cdot aa')  =  (S(h) \cdot b) aa' =  a \epsilon_r (b\# h) a'.
\end{eqnarray*}

One can easily verify the compatibility relations with the left and right counits with the left and right source and targets, that is, $s_l \circ \epsilon_l \circ t_r =t_r$, $t_l\circ \epsilon_l \circ s_r =s_r$, $s_r\circ \epsilon_r \circ t_l =t_l$ and $t_r\circ \epsilon_r \circ s_l =s_l$.

Let us verify that $(\underline{\wide{A}\#_{\omega}H}, \Delta_l , \epsilon_l )$ and 
$(\underline{\wide{A}\#_{\omega}H}, \Delta_r , \epsilon_r )$ are corings over $E({A})$ with their respectives bimodule
structures. we have already seen that $\Delta_l$, $\epsilon_l$, $\Delta_r$ and $\epsilon_r$ are morphisms of
$E({A})$ bimodules. It is easy to see that the left and right comultiplications are coassociative. 
It remains to verify the counit axiom for both structures. Take $a\# h \in \underline{\wide{A}\#_{\omega}H}$ ,then
\begin{eqnarray*}
& \, & (\epsilon_l \otimes_{E({A}), \triangleright \triangleleft} I)\circ \Delta_l (a\# h) =  \epsilon_l (a_{(1)} \# h_{(1)}) \triangleright (a_{(2)} \# h_{(2)})\\
& = &(\epsilon_{\wide{A}} (a_{(1)})(h_{(1)}\cdot 1_A) \# 1_H )(a_{(2)} \# h_{(2)}) =  \epsilon_{\wide{A}} (a_{(1)})(a_{(2)}(h_{(1)}\cdot 1_A)\# h_{(2)})\\
& = & a(h_{(1)}\cdot 1_A)\# h_{(2)}) = a\# h .
\end{eqnarray*}

Using the cocommutativity of $H$ and the commutativity of $A$ it is also easy to see that
$$
(I \otimes_{E({A}), \triangleright \triangleleft} \epsilon_l )\circ \Delta_l (a\# h) =a\# h .
$$

Therefore, $(\underline{\wide{A}\#_{\omega}H}, \Delta_l , \epsilon_l )$ is a coring over $E({A})$. For the right structure, for $a\# h \in \underline{\wide{A}\#_{\omega}H}$, we have
\begin{eqnarray*}
& \, &\hspace{-0.4cm} (\epsilon_r \otimes_{E({A}), \blacktriangleright \blacktriangleleft} I)\circ \Delta_r (a\# h) =  
\epsilon_r (a_{(1)} \# h_{(1)}) \blacktriangleright (a_{(2)} \# h_{(2)})\\
&\hspace{-0.5cm} = &\hspace{-0.4cm} (a_{(2)} \# h_{(2)}) ((S(h_{(1)})\cdot \epsilon_{\wide{A}} (a_{(1)}))\# 1_H ) =  (a_{(2)} \# h_{(1)}) ((S(h_{(2)})\cdot \epsilon_{\wide{A}} (a_{(1)}))\# 1_H )\\
&\hspace{-0.5cm} = &\hspace{-0.4cm} a_{(2)} (h_{(1)}\cdot (S(h_{(4)})\cdot \epsilon_{\wide{A}} (a_{(1)}))) \omega (h_{(2)}, 1_H )\# h_{(3)}  \\
&\hspace{-0.5cm} = &\hspace{-0.4cm}  a_{(2)} (h_{(1)}\cdot (S(h_{(4)})\hspace{-0.05cm}\cdot \epsilon_{\wide{A}} (a_{(1)}))) (h_{(2)}\cdot 1_A )\# h_{(3)} \\
&\hspace{-0.5cm} = &\hspace{-0.4cm} a_{(2)}\hspace{-0.05cm} (h_{(1)}S(h_{(4)})\cdot \epsilon_{\wide{A}}\hspace{-0.03cm} (a_{(1)})) (h_{(2)}\cdot 1_A )\# h_{(3)} 
\hspace{-0.1cm}= \hspace{-0.1cm} a_{(2)}\hspace{-0.05cm} (h_{(1)}S(h_{(2)})\hspace{-0.05cm}\cdot \epsilon_{\wide{A}}\hspace{-0.03cm} (a_{(1)})) (h_{(3)}\hspace{-0.05cm}\cdot 1_A )\# h_{(4)} \\
&\hspace{-0.5cm} = &\hspace{-0.4cm} a_{(2)} \epsilon_{\wide{A}} (a_{(1)}) (h_{(1)}\cdot 1_A )\# h_{(2)}  =  a (h_{(1)}\cdot 1_A)\# h_{(2)})=a\# h .
\end{eqnarray*}

Analogously, one can prove that  
\[
(I \otimes_{E({A}), \blacktriangleright \blacktriangleleft} \epsilon_r )\circ \Delta_r (a\# h) =a\# h .
\]  

Therefore, $(\underline{\wide{A}\#_{\omega}H}, \Delta_r , \epsilon_r )$ is a coring over $E({A})$  
  
Let us verify now that  
\[
\epsilon_l ((a\# h)(b\# k))=\epsilon_l ((a\# h) s_l(\epsilon_l (b\# k)))=\epsilon_l ((a\# h) t_l(\epsilon_l (b\# k))) , 
\]  
and 
\[ 
\epsilon_r ((a\# h)(b\# k)) =\epsilon_r (s_r(\epsilon_r (a\# h)) (b\# k))=\epsilon_r (t_r (\epsilon_r (a\# h))(b\# k)) ,
\]  
for any $a\# h , b\# k \in \underline{\wide{A}\#_{\omega}H}$. For the left counit, on one hand,
\begin{eqnarray*}
& \, & \epsilon_l((a\# h)(b\# k))  =  \epsilon_l (a(h_{(1)}\cdot b)\omega (h_{(2)},k_{(1)})\# h_{(3)}k_{(2)})\\
   & = & \epsilon_{\wide{A}}(a)\epsilon_{\wide{A}}(h_{(1)}\cdot b) 
\epsilon_{\wide{A}}(\omega (h_{(2)},k_{(1)}))(h_{(3)}k_{(2)}\cdot 1_A)\\
   & = & \epsilon_{\wide{A}}(a) (h_{(1)}\cdot \epsilon_{\wide{A}}(b))(h_{(2)}\cdot (k_{(1)}\cdot 1_A))(h_{(3)}k_{(2)}\cdot 1_A)\\
   & = & \epsilon_{\wide{A}}(a)(h_{(1)}\cdot \epsilon_{\wide{A}}(b))(h_{(2)}\cdot (k\cdot 1_A)) =  \epsilon_{\wide{A}}(a)(h\cdot (\epsilon_{\wide{A}}(b)(k\cdot 1_A))) .  
\end{eqnarray*}

On the other hand
\begin{eqnarray*}
& \, &\hspace{-0.55cm} \epsilon_l ((a\# h)s_l (\epsilon_l (b\# k)))  =  \epsilon_l ((a\# h)s_l(\epsilon_{\wide{A}}(b)(k\cdot 1_A)))\\
   &\hspace{-0.65cm} = &\hspace{-0.55cm} \epsilon_l ((a\# h)(\epsilon_{\wide{A}}(b)(k\cdot 1_A)\# 1_H)) =  \epsilon_l (a(h_{(1)}\cdot (\epsilon_{\wide{A}}(b)(k\cdot 1_A)))\omega (h_{(2)},1_H)\# h_{(3)})\\   
   &\hspace{-0.65cm} = &\hspace{-0.55cm} \epsilon_l (a(h_{(1)}\cdot (\epsilon_{\wide{A}}(b)(k\cdot 1_A)))(h_{(2)}\cdot 1_H)\# h_{(3)}) = \epsilon_{\wide{A}}(a(h_{(1)}\cdot (\epsilon_{\wide{A}}(b)(k\cdot 1_A))))(h_{(2)}\cdot 1_A)\\
   &\hspace{-0.65cm} = &\hspace{-0.55cm} \epsilon_{\wide{A}}(a)(h_{(1)}\cdot (\epsilon_{\wide{A}}(b)(k\cdot 1_A)))(h_{(2)}\cdot 1_A) =  \epsilon_{\wide{A}}(a)(h\cdot (\epsilon_{\wide{A}}(b)(k\cdot 1_A))).
\end{eqnarray*}

Therefore, $\epsilon_l ((a\# h)(b\# k))=\epsilon_l ((a\# h) s_l(\epsilon_l (b\# k)))=\epsilon_l ((a\# h) t_l(\epsilon_l (b\# k)))$. For the right counit, on one hand
\begin{eqnarray*} 
& &\hspace{-0.45cm} \epsilon_r ((a\# h)(b\# k)) =  \epsilon_r (a(h_{(1)}\cdot b) \omega (h_{(2)},k_{(1)})\# h_{(3)}k_{(2)})\\
   &\hspace{-0.65cm} = &\hspace{-0.45cm} S(h_{(3)}k_{(2)})\cdot \epsilon_{\wide{A}}(a(h_{(1)}\cdot b)\omega (h_{(2)},k_{(1)}))\\
   &\hspace{-0.65cm} = &\hspace{-0.45cm} S(h_{(3)}k_{(2)})\cdot(\epsilon_{\wide{A}}(a)(h_{(1)}\cdot \epsilon_{\wide{A}}(b))(h_{(2)}\cdot (k_{(1)}\cdot 1_A))\\
   &\hspace{-0.65cm} = &\hspace{-0.45cm} (S(h_{(5)}k_{(4)})\cdot \epsilon_{\wide{A}}(a))(S(h_{(4)}k_{(3)})\cdot (h_{(1)}\cdot \epsilon_{\wide{A}}(b)))
    (S(h_{(3)}k_{(2)})\cdot (h_{(2)}\cdot (k_{(1)}\cdot 1_A)))\\
    &\hspace{-0.65cm} = &\hspace{-0.45cm} (S(h_{(5)}k_{(4)})\cdot \epsilon_{\wide{A}}(a))(S(h_{(4)}k_{(3)})h_{(1)}\cdot \epsilon_{\wide{A}}(b)) (S(h_{(3)}k_{(2)})h_{(2)}\cdot (k_{(1)}\cdot 1_A))\\
    &\hspace{-0.65cm} = &\hspace{-0.45cm} (S(k_{(4)})S(h_{(5)})\cdot \epsilon_{\wide{A}}(a))(S(k_{(3)})S(h_{(4)})h_{(1)}\cdot \epsilon_{\wide{A}}(b))
    (S(k_{(2)})S(h_{(3)})h_{(2)}\cdot (k_{(1)}\cdot 1_A))\\
    &\hspace{-0.65cm} = &\hspace{-0.45cm} (S(k_{(4)})S(h_{(5)})\cdot \epsilon_{\wide{A}}(a))(S(k_{(3)})S(h_{(3)})h_{(4)}\cdot \epsilon_{\wide{A}}(b))
    (S(k_{(2)})S(h_{(1)})h_{(2)}\cdot (k_{(1)}\cdot 1_A))\\
    &\hspace{-0.65cm} = &\hspace{-0.45cm} (S(k_{(4)})S(h)\cdot \epsilon_{\wide{A}}(a))(S(k_{(3)})\cdot \epsilon_{\wide{A}}(b)) (S(k_{(2)})\cdot (k_{(1)}\cdot 1_A))\\
   &\hspace{-0.65cm} = &\hspace{-0.45cm} (S(k_{(4)})S(h)\cdot \epsilon_{\wide{A}}(a))(S(k_{(3)})\cdot \epsilon_{\wide{A}}(b)) (S(k_{(2)})k_{(1)}\cdot 1_A))\\
   &\hspace{-0.65cm} = &\hspace{-0.45cm} (S(k_{(4)})S(h)\cdot \epsilon_{\wide{A}}(a))(S(k_{(3)})\cdot \epsilon_{\wide{A}}(b)) (S(k_{(1)})k_{(2)}\cdot 1_A))\\
   &\hspace{-0.65cm} = &\hspace{-0.45cm} (S(k_{(2)})S(h)\cdot \epsilon_{\wide{A}}(a))(S(k_{(1)})\cdot \epsilon_{\wide{A}}(b)) \\
   &\hspace{-0.65cm} = &\hspace{-0.45cm} (S(k_{(2)})\cdot (S(h)\cdot  \epsilon_{\wide{A}}(a)))(S(k_{(1)})\cdot \epsilon_{\wide{A}}(b)) \\
   &\hspace{-0.65cm} = &\hspace{-0.45cm} S(k) \cdot ((S(h)\cdot \epsilon_{\wide{A}} (a))\epsilon_{\wide{A}}(b)) .
\end{eqnarray*}

On the other hand,
\begin{eqnarray*}
& & \epsilon_r (s_r (\epsilon_r (a\# h))(b\# k))  =  \epsilon_r ((S(h)\cdot \epsilon_{\wide{A}}(a)) \# 1_H)(b\# k))\\
   & = & \epsilon_r ((S(h)\cdot \epsilon_{\wide{A}}(a))b\# k)
    =  S(k)\cdot ((S(h)\cdot \epsilon_{\wide{A}}(a))\epsilon_{\wide{A}}(b)).
\end{eqnarray*}

Therefore, $\epsilon_r ((a\# h)(b\# k)) =\epsilon_r (s_r(\epsilon_r (a\# h)) (b\# k))=\epsilon_r (t_r (\epsilon_r (a\# h))(b\# k))$.

Finally, we define the antipode as
$$\func{\mathcal{S}:}{\underline{\wide{A}\#_{\omega} H}}{\underline{\wide{A}\#_{\omega} H}}{a\#h}{ (S_H(h_{(3)})\cdot S_{\wide{A}}(a))\omega^{-1}(S_H(h_{(2)}),h_{(4)})\#S_H(h_{(1)})}$$

Take $b, c \in E({A})$ and $a\#h \in \underline{\wide{A}\#_{\omega} H}$, then, one can prove that
\[
\mathcal{S}(t(b)(a\#h)t(c)) = s(c)\mathcal{S}(a\#h)s(b).
\]

Indeed,  
\begin{eqnarray*}
& \, & \mathcal{S}(t_l(b)(a\# h)t_r(c))  =  \mathcal{S}((b\# 1_H)(a\# h)(c\# 1_H)) =  \mathcal{S}((ba\# h)(c\# 1_H))\\
   & = & \mathcal{S}(ba(h_{(1)}\cdot c)\omega (h_{(2)},1_H)\# h_{(3)}) =  \mathcal{S}(ba(h_{(1)}\cdot c)\# h_{(2)})\\
   & = & (S(h_{(4)})\cdot S_{\wide{A}}(ba(h_{(1)}\cdot c)))\omega^{-1}(S(h_{(3)}),h_{(5)})\# S(h_{(2)})\\
   & \stackrel{(*)}{=} & (S(h_{(4)})\cdot ((h_{(1)}\cdot c)S_{\wide{A}}(a)b))\omega^{-1}(S(h_{(3)}),h_{(5)})\# S(h_{(2)})\\
\end{eqnarray*}
in which the equality $(*)$ is valid because $S_{\wide{A}}$ restricted to the base algebra $E({A})$ is equal to the identity. On the other hand.

\begin{eqnarray*}
   &   & s_r (c) \mathcal{S} (a\# h) s_l (b) =\\
   & = & (c\# 1_H)((S(h_{(3)})\cdot S_{\wide{A}}(a))\omega^{-1}(S(h_{(2)}),h_{(4)})\# S(h_{(1)}))(b\# 1_H)\\
   & = & (c(S(h_{(4)})\cdot S_{\wide{A}}(a))\omega^{-1}(S(h_{(3)}),h_{(5)})\omega (1_H,S(h_{(2)}))\# S(h_{(1)}))(b\# 1_H)\\
   & = & (c(S(h_{(3)})\cdot S_{\wide{A}}(a))\omega^{-1}(S(h_{(2)}),h_{(4)})\# S(h_{(1)}))(b\# 1_H)\\
   & = & c(S(h_{(4)})\cdot S_{\wide{A}}(a))\omega^{-1}(S(h_{(3)}),h_{(5)})( S(h_{(2)})\cdot b)\# S(h_{(1)})\\
   & = & ( S(h_{(6)})\cdot b)(S(h_{(5)})\cdot S_{\wide{A}}(a))(S(h_{(3)})h_{(4)}\cdot c) \omega^{-1}(S(h_{(2)}),h_{(7)})\# S(h_{(1)})\\
   & = & ( S(h_{(6)})\cdot b)(S(h_{(5)})\cdot S_{\wide{A}}(a))(S(h_{(3)})\cdot (h_{(4)}\cdot c))
 \omega^{-1}(S(h_{(2)}),h_{(7)})\# S(h_{(1)})\\
   & = & ( S(h_{(6)})\cdot b)(S(h_{(5)})\cdot S_{\wide{A}}(a))(S(h_{(4)})\cdot (h_{(1)}\cdot c))
 \omega^{-1}(S(h_{(3)}),h_{(7)})\# S(h_{(2)})\\
   & = & ( S(h_{(4)})\cdot ((h_{(1)}\cdot c)S_{\wide{A}}(a) b)) \omega^{-1}(S(h_{(3)}),h_{(5)})\# S(h_{(2)}).
\end{eqnarray*}

It remains to check that 
\[
\mu (\mathcal{S}\otimes_{E({A}), \triangleright \triangleleft} Id)\circ \Delta_l= s_r \circ \epsilon_r
\]
and 
\[
\mu(Id \otimes_{E({A}), \blacktriangleright \blacktriangleleft} \mathcal{S})\circ \Delta_r= s_l \circ \epsilon_l .
\]

Take $a\# h \in \underline{\wide{A}\#_{\omega} H}$, then

\begin{eqnarray*}
&  & \mu (\mathcal{S}\otimes_{E({A}), \triangleright \triangleleft} Id)\circ \Delta_l (a\# h) = \mathcal{S}(a_{(1)}\# h_{(1)})(a_{(2)}\# h_{(2)})\\
   & = & (S(h_{(3)})\cdot S_{\wide{A}}(a_{(1)})\omega^{-1}(S(h_{(2)}),h_{(4)})\# S(h_{(1)}))(a_{(2)}\# h_{(5)})\\
   & = & (S(h_{(5)})\cdot S_{\wide{A}}(a_{(1)}))\omega^{-1}(S(h_{(4)}),h_{(6)})(S(h_{(3)})\cdot a_{(2)})
     \omega (S(h_{(2)}),h_{(7)})\# S(h_{(1)})h_{(8)}\\
   & = & (S(h_{(4)})\cdot S_{\wide{A}}(a_{(1)})) (S(h_{(3)})\cdot a_{(2)})\omega^{-1}(S(h_{(2)}),h_{(5)})
    \omega (S(h_{(1)}),h_{(6)})\# S(h_{(7)})h_{(8)}\\
   & = & (S(h_{(3)})\cdot S_{\wide{A}}(a_{(1)}))(S(h_{(2)})\cdot a_{(2)}) (S(h_{(1)})\cdot (h_{(4)}\cdot 1_A))\# 1_H\\
   & = & (S(h_{(3)})\cdot S_{\wide{A}}(a_{(1)}))(S(h_{(2)})\cdot a_{(2)}) (S(h_{(1)})h_{(4)}\cdot 1_A)\# 1_H\\
   & = & (S(h_{(2)})\cdot S_{\wide{A}}(a_{(1)}))(S(h_{(1)})\cdot a_{(2)})\# 1_H\\
   & = & (S(h)\cdot (S_{\wide{A}}(a_{(1)})a_{(2)}))\# 1_H \\
   & = & s_r (S(h)\cdot \epsilon_{\wide{A}}(a))\\
   & = & s\circ \epsilon_r (a\# h) .
  \end{eqnarray*}
and
\begin{eqnarray*} 
&  &\hspace{-0.4cm}  \mu (Id\otimes_{E({A}), \blacktriangleright \blacktriangleleft} \mathcal{S})\circ \Delta_r (a\# h) =  
(a_{(1)}\# h_{(1)})\mathcal{S}(a_{(2)}\# h_{(2)})\\
   &\hspace{-0.38cm}  = &\hspace{-0.4cm} (a_{(1)}\# h_{(1)})(S(h_{(4)})\cdot S_{\wide{A}}(a_{(2)})\omega^{-1}(S(h_{(3)}),h_{(5)})\# S(h_{(2)}))\\
   &\hspace{-0.38cm} = &\hspace{-0.4cm} (a_{(1)}(h_{(1)}\cdot (S(h_{(7)})\cdot S_{\wide{A}}(a_{(2)})\omega^{-1}(S(h_{(6)}),h_{(8)})))
   \omega (h_{(2)},S(h_{(5)}))\# h_{(3)}S(h_{(4)})\\
   &\hspace{-0.38cm} = &\hspace{-0.4cm} (a_{(1)}(h_{(1)}\cdot (S(h_{(6)})\cdot S_{\wide{A}}(a_{(2)})))(h_{(2)}\cdot  \omega^{-1}(S(h_{(5)}),h_{(7)})) 
    \omega (h_{(3)},S(h_{(4)}))\# 1_H\\
   &\hspace{-0.38cm} = &\hspace{-0.4cm} (a_{(1)}(h_{(1)}\cdot (S(h_{(2)})\cdot S_{\wide{A}}(a_{(2)})))(h_{(3)}\cdot  \omega^{-1}(S(h_{(6)}),h_{(7)}))
   \omega (h_{(4)},S(h_{(5)}))\# 1_H\\
   &\hspace{-0.38cm} = &\hspace{-0.4cm} (a_{(1)}(h_{(1)}S(h_{(2)})\cdot S_{\wide{A}}(a_{(2)})))(h_{(3)}\cdot  \omega^{-1}(S(h_{(6)}),h_{(7)}))
   \omega (h_{(4)},S(h_{(5)}))\# 1_H\\   
   &\hspace{-0.38cm} = &\hspace{-0.4cm} a_{(1)}S_{\wide{A}}(a_{(2)})(h_{(1)}\cdot  \omega^{-1}(S(h_{(4)}),h_{(5)}))
    \omega (h_{(2)},S(h_{(3)}))\# 1_H\\
       \end{eqnarray*}
   \begin{eqnarray*}
   &\hspace{-0.38cm}\stackrel{(*)}{=} &\hspace{-0.4cm} \epsilon_{\wide{A}} (a) \omega^{-1}\hspace{-0.02cm}(\hspace{-0.01cm}h_{(1)}S(h_{(8)}),\hspace{-0.01cm} h_{(9)})\hspace{-0.01cm} \omega\hspace{-0.02cm} (\hspace{-0.01cm}h_{(2)} ,\hspace{-0.01cm} S(h_{(7)}) h_{(10)})
     \omega^{-1}\hspace{-0.02cm} (\hspace{-0.01cm}h_{(3)},\hspace{-0.01cm} S(h_{(6)})\hspace{-0.01cm})\omega\hspace{-0.02cm} (\hspace{-0.01cm}h_{(4)},\hspace{-0.01cm}S(h_{(5)})\hspace{-0.01cm})\hspace{-0.02cm}\#\hspace{-0.02cm} 1_H\\
    &\hspace{-0.38cm} = &\hspace{-0.4cm} \epsilon_{\wide{A}} (a) \omega^{-1}(h_{(1)}S(h_{(6)}), h_{(7)}) \omega (h_{(2)} , S(h_{(5)}) h_{(8)}) 
    (h_{(3)}\cdot (S(h_{(4)})\cdot 1_A)) \# 1_H \\
   &\hspace{-0.38cm} = &\hspace{-0.4cm} \epsilon_{\wide{A}} (a) \omega^{-1}(h_{(1)}S(h_{(5)}), h_{(6)}) \omega (h_{(2)} , S(h_{(4)}) h_{(7)}) (h_{(3)}\cdot 1_A) \# 1_H \\  
   &\hspace{-0.38cm} = &\hspace{-0.4cm} \epsilon_{\wide{A}} (a) \omega^{-1}(h_{(1)}S(h_{(2)}), h_{(7)}) \omega (h_{(3)} , S(h_{(5)}) h_{(6)}) (h_{(4)}\cdot 1_A) \# 1_H \\
   &\hspace{-0.38cm} = &\hspace{-0.4cm} \epsilon_{\wide{A}} (a) \omega^{-1}(1_H, h_{(3)}) \omega (h_{(1)} , 1_H )(h_{(2)}\cdot 1_A) \# 1_H \\
   &\hspace{-0.38cm} = &\hspace{-0.4cm} \epsilon_{\wide{A}}(a)(h\cdot 1_A)\# 1_H =  s_l \circ \epsilon_l (a\# h).
\end{eqnarray*}

The equality $(*)$ we used that, for any $h,k,l\in H$, we have
\[
h\cdot \omega^{-1} (k,l) =\omega^{-1} (h_{(1)}k_{(1)}, l_{(1)})\omega (h_{(2)}, k_{(2)}l_{(2)})\omega^{-1} (h_{(3)}, k_{(3)}) .
\]
This folllows easily from
\[
h\cdot (\omega (k_{(1)} , l_{(1)})\omega^{-1} (k_{(2)} , l_{(2)})=h\cdot (k\cdot (l\cdot 1_A )) 
\]
and 
\[
h\cdot \omega (k,l) = \omega (h_{(1)}, k_{(1)}) \omega^{-1} (h_{(2)}, k_{(2)}l_{(1)})\omega(h_{(3)}k_{(3)}, l_{(2)}) ,
\]
which, in its turn is an immediate consequence of the $2$-cocycle identity.

Therefore, $(\underline{\wide{A}\#_{\omega} H}, s_l, t_l, s_r, t_r, \Delta_l, \Delta_r, \epsilon_l, \epsilon_r , \mathcal{S})$ is a strucutre of a Hopf algebroid over $E({A})$. \find

\section{Partially cleft extensions and cleft extensions by Hopf algebroids}

In \cite{ABDP1}, the authors introduced the notion of a partially cleft extension of an algebra $A$ by e Hopf algebra $H$. Using the Hopf algebroid structure of the crossed product, one can rethink partially cleft extensions of commutative algebras by co-commutative Hopf algebras in a broader scenario, namely, the theory of cleft extensions of algebras by Hopf algebroids developed by G. B\"{o}hm and T. Brzezinski \cite{BB}.

\begin{defi} \label{partialcleft} \cite{ABDP1} Let $H$ be a Hopf algebra, and $A\subset B$ be an $H$-extension, that is $B$ is a right $H$-comodule algebra and $A=B^{coH}$. The extension $A\subset B$ is partially cleft if there exists a pair of linear maps $\gamma,\ov{\gamma}: H\ri B$ such that:
\begin{enumerate}[(i)]
\item $\gamma(1_H) = 1_B$;
\item The following diagrams are commutative  
\[ \xymatrix{
 H \ar[rr]^{\gamma}\ar[d]^{\Delta} &  & B  \ar[d]^-{\rho} & &  H \ar[rr]^{\ov{\gamma}}\ar[d]^{\Delta^{cop}} &  & 
B  \ar[d]^-{\rho}\\
 H\otimes H \ar[rr]_{\gamma \otimes Id_H}& & B\otimes H & & H\otimes H \ar[rr]_{\ov{\gamma} \otimes S}& & B\otimes H} 
\]

\item $(\gamma *\ov{\gamma})\circ \mu$ is a central element in the convolution algebra $\mbox{Hom}_k (H\otimes H,B)$, in which $\mu:H\otimes H\ri H$ is the multiplication in $H$ and $(\ov{\gamma} *{\gamma})(h)$ commutes with every element of  $A$, for each $h\in H$,

\noindent and, for all $b\in B$, $h,l\in H$, if we write $e_h = (\gamma*\ov{\gamma})(h)$ and $\wide{e}_h = (\ov{\gamma}*\gamma)(h)$, then:
\item $\sum b_{(0)}\ov{\gamma}(b_{(1)})\gamma(b_{(2)}) = b$;
\item $\gamma(h)e_l = \sum e_{h_{(1)}l}\gamma(h_{(2)})$;
\item $\ov{\gamma}(l)\wide{e}_h = \sum \wide{e}_{hl_{(1)}}\ov{\gamma}(l_{(2)})$;
\item $\sum \gamma(hl_{(1)})\wide{e}_{l_{(2)}} = \sum e_{h_{(1)}}\gamma(h_{(2)}l)$. 
\end{enumerate}
\end{defi}

Partially cleft extensions are related to partial crossed products, as one can see in the next two results.

\begin{prop} \cite{ABDP1} If $(A,\cdot,(\omega,\omega^{-1}))$ is a symmetric twisted partial $H$- module algebra with a $2$-cocycle $\omega$, then,$A\subset \underline{A\#_{\omega} H}$ is a partially cleft $H$-extension. 
\end{prop}

For the crossed product $\underline{A\#_{\omega} H}$, the cleaving maps $\gamma,\ov{\gamma} :H\rightarrow \underline{A\#_{\omega} H}$ are given by
\begin{equation}\label{cleaving}
\gamma (h) =1_A \# h , \qquad \mbox{ and } \qquad \ov{\gamma} (h) =\omega^{-1} (S(h_{(2)}), h_{(3)})\# S(h_{(1)}) .
\end{equation}

\begin{thm} \cite{ABDP1}
Let $B$ be $H$-comodule algebra and $A=B^{coH}$. Then the $H$-extension $A\subset B$ is partially cleft if, and only if there is a symmetric twisted partial action $\cdot :H\otimes A \rightarrow A$ with a $2$-cocycle $\omega :H\otimes H \rightarrow A$ such that $B$ is isomorphic to the partial crossed product $\underline{A\#_{\omega}H}$.
\end{thm}

In the case of a co-commutative Hopf algebra $H$ acting partially over a commutative algebra $A$, we have already seen that we can replace $A$ by a commutative and co-commutative Hopf algebra $\wide{A}$ with the same cohomology theory, moreover, the crossed product $\underline{\wide{A} \#_{\omega} H}$ has a structure of Hopf algebroid over the base algebra $E({A})$. Then it is interesting to see whether one can replace the crossed product $\underline{A\#_{\omega} H}$ by the crossed product $\underline{\wide{A}\#_{\omega} H}$ in the analysis of cleft extensions by the Hopf algebra $H$. First note that the 
$H$-comodule structure on both crossed products is the same, namely $\rho (a\# h)=a\# h_{(1)}\otimes h_{(2)}$. Furthermore $(\underline{A\#_{\omega} H})^{coH} \cong A$ and 
$(\underline{\wide{A}\#_{\omega} H})^{coH} \cong \wide{A}$, then both crossed products are $H$-extensions of their respectives algebras of coinvariants. Finally, the cleaving maps $\gamma,\ov{\gamma} :H\rightarrow \underline{A\#_{\omega} H}$, given by (\ref{cleaving}), take their values actually in $\underline{\wide{A} \#_{\omega} H}$. 

In what follows, we shall see that there exists a Hopf algebroid $\mathcal{H}$ such that the crossed product $\underline{\wide{A}\#_{\omega} H}$ can be viewed as a cleft extension of $\wide{A}$ by $\mathcal{H}$. 

\begin{defi} \label{cleftalgebroid} \cite{BB} Let $(\mathcal{H}, L, R, s_L , t_L , s_R , t_R , \Delta_L , \Delta_R , \epsilon_L , \epsilon_R , \mathcal{S})$ be a Hopf algebroid and $A$ be a right $\mathcal{H}$-comodule algebra. Denote by $\eta_R (r) =r\cdot 1_A =1_A \cdot r$ the unit map of the corresponding $R$-ring structure of $A$. Let $B$ be the subalgebra of $\mathcal{H}_R$-coinvariants in $A$. The extension $B\subset A$ is called $\mathcal{H}$-cleft if
\begin{enumerate}[(a)]
\item $A$ is an $L$-ring (with unit $\eta_L:L\rightarrow A$) and $B$ is an $L$-subring of $A$; 
\item there\ \ exists\ \ a\ \ convolution\ \ invertible\ \ left\ \ $L$-linear\ \ right\ \ $\mathcal{H}$-colinear\ \
morphism $\gamma :\mathcal{H}\rightarrow A$.
\end{enumerate}
A map $\gamma$ satisfying condition (b) is called a cleaving map.
\end{defi}

\begin{rmk} Some small remarks have to be made about this definition. 
\begin{enumerate}[(1)]
\item First is that the structure of right $\mathcal{H}$-comodule algebra on $A$ is relative to the base ring $R$, that is $\rho: A\rightarrow A\otimes_R \mathcal{H}$ is a right $R-R$-bilinear map in the sense that $\rho(\eta_R (r) a\eta_R (s))=a^{(0})\otimes_R s_R (r)a^{(1)} s_R (s)$, for every $a\in A$ and $r\in R$. 
\item  The map $\gamma :\mathcal{H} \rightarrow A$ being $\mathcal{H}$-colinear implies that it is  right $R$-linear in the sense that $\gamma (hs_R (r))=\gamma (h)\eta_R (r)$ and left $R$-linear in the sense that $\gamma (s_R (r) h)=\eta_R (r) \gamma (h)$, for any $a\in A$ and $r\in R$.
\item The notion of a convolution invertible map $\gamma :\mathcal{H} \rightarrow A$, in which $A$ is both $L$-ring and $R$-ring,  means that there is a unique $\ov{\gamma} :\mathcal{H} \rightarrow A$ such that \cite{BB}
\begin{eqnarray*}
\mu_A \circ (\gamma \otimes_R \ov{\gamma} )\circ \Delta_R & = & \eta_L \circ \epsilon_L \\
\mu_A \circ (\ov{\gamma} \otimes_L {\gamma} )\circ \Delta_L & = & \eta_R \circ \epsilon_R .
\end{eqnarray*}
\item The left $L$-linearity of the map $\gamma$ in condition (b) of Definition \ref{cleftalgebroid} means, that the cleaving map satisfies $\gamma (s_L (l)h)=\eta_L (l) \gamma (h)$, for any $a\in A$ and $l\in L$.
\end{enumerate}
\end{rmk}

In our case, for a co-commutative Hopf algebra $H$ acting partially upon a commutative and co-commutative Hopf algebra $\wide{A}$, the base algebras $L$ and $R$ will coincide with the commutative subalgebra $E({A})$ of $\wide{A}$ and then many distinctions about the left and right structures will coalesce. The Hopf algebroid is given by the partial smash product $\underline{E({A})\# H}$. For a proof that this partial smash product is in fact a Hopf algebroid over $E({A})$, see reference \cite{BV}, Theorem 3.5. The extension to be considered is the previous defined partially $H$-cleft extension $\wide{A} \subset \underline{\wide{A}\#_{\omega} H}$. Then we have the following theorem.

\begin{thm} \label{partiallycleftisclefthopfalgebroid} Let $H$ be a co-commutative Hopf algebra acting partially on a commutative and co-commutative Hopf algebra $\wide{A}$ and let $\omega$ be a partial $2$ cocycle in $H^2_{par} (H,\wide{A})$. Then the partial crossed product $\underline{\wide{A}\#_{\omega} H}$ is a right $\mathcal{H}=\underline{E({A})\# H}$-module algebra with $\wide{A} \cong (\underline{\wide{A}\#_{\omega} H})^{co\mathcal{H}}$. Moreover, the extension $\wide{A} \subset \underline{\wide{A}\#_{\omega} H}$ is $\mathcal{H}$-cleft in the sense of Definition \ref{cleftalgebroid}.
\end{thm}

\noindent \dem First, define the linear map
$$\func{\widetilde{\rho}:}{\prodah}{\prodah\otimes_{E({A}), \blacktriangleright \blacktriangleleft}\underline{E({A})\#H}}{a\#h}
{a\# h_{(1)}\otimes (h_{(2)}\cdot1_A)\# h_{(3)}}$$

Note that the expression of $\widetilde{\rho}(a\# h)$ can also be written as
\[
\widetilde{\rho}(a\# h)= a\# h_{(1)}\otimes 1_{E({A})} \# h_{(2)} =a\# h_{(1)}\otimes 1_{A} \# h_{(2)} ,
\]
that is because
\begin{eqnarray*}
& & a\#h_{(1)}\otimes (h_{(2)}\cdot1_A)\# h_{(3)} =a\#h_{(1)}\otimes (1_A \# h_{(2)})(1_A \# 1_H ) \\
& = & a\# h_{(1)}\otimes (1_A \# h_{(2)})\blacktriangleleft 1_A =a\# h_{(1)}\otimes (1_A \# h_{(2)}) .
\end{eqnarray*}

It is easy to see that $(\wide{\rho}\otimes_{E({A}), \blacktriangleright \blacktriangleleft} Id )\circ \wide{\rho} = (Id \otimes_{E({A}), \blacktriangleright \blacktriangleleft} \wide{\Delta}_r )\circ \wide{\rho}$, in which $\wide{\Delta}_r$ is the right comultiplication in $\mathcal{H}$, given by
\[
\wide{\Delta}_r (r\# h )=r \# h_{(1)} \otimes 1_A \# h_{(2)} , \qquad \forall r\in E({A}), \forall h\in H .
\]

Also, one can check that the map $\wide{\rho}$ is $E({A})$-bilinear. Indeed, for $a\# h \in \underline {\wide{A}\#_{\omega} H}$ and $r,s\in E({A})$ then
\begin{eqnarray*}
& \, & \wide{\rho} ((r\# 1_H)(a\# h)(s\# 1_H))  =  \wide{\rho} (ra(h_{(1)}\cdot s)\# h_{(2)}) \\
& = & ra(h_{(1)}\cdot s)\# h_{(2)} \otimes 1_A \# h_{(3)} =  (ra\# h_{(1)})\blacktriangleleft s \otimes 1_A \# h_{(2)}\\
& = & ra\# h_{(1)} \otimes s\blacktriangleright (1_A \# h_{(2)}) =  ra\# h_{(1)} \otimes (1_A \# h_{(2)})(s\# 1_H)\\
& = & a(h_{(1)}S(h_{(2)})\cdot r)\# h_{(3)} \otimes (1_A \# h_{(4)})(s\# 1_H) \\
& = & a(h_{(1)}\cdot (S(h_{(3)})\cdot r))\# h_{(2)} \otimes (1_A \# h_{(4)})(s\# 1_H) \\
& = & (a\# h_{(1)}) \blacktriangleleft (S(h_{(2)})\cdot r)\otimes (1_A \# h_{(3)})(s\# 1_H) \\
& = & a\# h_{(1)}  \otimes (S(h_{(2)})\cdot r)\blacktriangleright (1_A \# h_{(3)})(s\# 1_H) \\
& = & a\# h_{(1)}  \otimes (1_A \# h_{(2)})((S(h_{(3)})\cdot r)\# 1_H)(s\# 1_H) \\
& = & a\# h_{(1)}  \otimes ((h_{(2)}\cdot (S(h_{(4)})\cdot r)) \# h_{(3)})(s\# 1_H) \\
& = & a\# h_{(1)}  \otimes ((h_{(2)}S(h_{(4)})\cdot r) \# h_{(3)})(s\# 1_H) \\
& = & a\# h_{(1)}  \otimes (r\# h_{(2)})(s\# 1_H)  =  a\# h_{(1)}  \otimes (r\# 1_H)(1_A \# h_{(2)})(s\# 1_H) 
\end{eqnarray*}

Denote by $\widetilde{\epsilon}_r :\underline{E({A})\#H} \rightarrow E({A})$ the right counit of the partial smash product $\underline{E({A})\#H}$, given by 
$\widetilde{\epsilon}_r (a\# h)=S(h)\cdot 1_A$. Then, we have
\begin{eqnarray*}
& \, & (Id\otimes_{E({A}), \blacktriangleright \blacktriangleleft} \wide{\epsilon}_r )\circ \wide{\rho} (a\# h)  =  (a\# h_{(1)} )\blacktriangleleft \wide{\epsilon}_r (1_A \# h_{(2)}) \\
& = & (a\# h_{(1)} )((S(h_{(2)})\cdot 1_A)\# 1_H)  =  a (h_{(1)}\cdot (S(h_{(3)})\cdot 1_A)) \# h_{(2)} \\
& = & a (h_{(1)}S(h_{(2)})\cdot 1_A) \# h_{(3)}  =  a\# h .
\end{eqnarray*}

Finally, for $a\# h, \ b\# k \in \underline{\wide{A}\#_{\omega} H}$, we have
\begin{eqnarray*}
\wide{\rho} ((a\# h)(b\# k)) & = & \wide{\rho} (a(h_{(1)}\cdot b)\omega (h_{(2)}, k_{(1)})\# h_{(3)}k_{(2)})\\
& = & a(h_{(1)}\cdot b)\omega (h_{(2)}, k_{(1)})\# h_{(3)}k_{(2)} \otimes 1_A \# h_{(4)}k_{(3)} ,
\end{eqnarray*}
on the other hand,

\begin{eqnarray*}
& & \wide{\rho}(a\# h) \wide{\rho}(b\# k)  =  (a\# h_{(1)})(b\# k_{(1)}) \otimes (1_A \# h_{(2)})(1_A \# k_{(2)})\\
& = & a(h_{(1)}\cdot b) \omega (h_{(2)},k_{(1)})\# h_{(3)}k_{(2)} \otimes (h_{(4)}\cdot 1_A)\# h_{(5)}k_{(3)} \\
& = & a(h_{(1)}\cdot b) \omega (h_{(2)},k_{(1)})\# h_{(3)}k_{(2)} \otimes (h_{(4)}k_{(3)}S(k_{(4)})\cdot 1_A)\# h_{(5)}k_{(5)}\\
& = & a(h_{(1)}\cdot b) \omega (h_{(2)},k_{(1)})\# h_{(3)}k_{(2)} \otimes (h_{(4)}k_{(3)}\cdot (S(k_{(4)})\cdot 1_A))\# h_{(5)}k_{(4)}\\
& = & a(h_{(1)}\cdot b) \omega (h_{(2)},k_{(1)})\# h_{(3)}k_{(2)} \otimes (S(k_{(4)})\cdot 1_A)\blacktriangleright (1_A \# h_{(4)}k_{(3)})\\
& = & (a(h_{(1)}\cdot b) \omega (h_{(2)},k_{(1)})\# h_{(3)}k_{(2)}) \blacktriangleleft (S(k_{(4)})\cdot 1_A)\otimes  1_A \# h_{(4)}k_{(3)}\\
& = & a(h_{(1)}\cdot b) \omega (h_{(2)},k_{(1)})  (h_{(3)}k_{(2)}\cdot (S(k_{(5)})\cdot 1_A)) \# h_{(4)}k_{(3)} \otimes  1_A \# h_{(5)}k_{(4)}\\
& = & (a(h_{(1)}\cdot b) \omega (h_{(2)},k_{(1)})  (h_{(3)}k_{(2)}S(k_{(3)})\cdot 1_A) \# h_{(4)}k_{(4)}) \otimes  (1_A \# h_{(5)}k_{(5)}\\
& = & a(h_{(1)}\cdot b) \omega (h_{(2)},k_{(1)})  (h_{(3)}\cdot 1_A) \# h_{(4)}k_{(2)} \otimes  1_A \# h_{(5)}k_{(3)}\\
& = & a(h_{(1)}\cdot b) \omega (h_{(2)},k_{(1)})  \# h_{(3)}k_{(2)} \otimes  1_A \# h_{(4)}k_{(3)}
\end{eqnarray*}

Therefore, the crossed product $\underline{\wide{A}\#_{\omega} H}$ is a right $\mathcal{H}$-comodule algebra. It is obvious that $i(\wide{A}) \subseteq (\underline{\wide{A}\#_{\omega} H})^{co\mathcal{H}}$, now take $\sum_i a_i \# h_i \in (\underline{\wide{A}\#_{\omega} H})^{co\mathcal{H}}$, then
\[
\sum_i a_i \# h_{i(1)} \otimes 1_A \# h_{i(2)} =\sum_i a_i \# h_i \otimes 1_A \# 1_H \in \wide{A} \otimes H 
\otimes_{E({A})} E({A}) \otimes H.
\]

Applying\ \ $Id\otimes \epsilon_H \otimes Id \otimes Id$\ \ to\ \ this\ \ identity\ \ and\ \ doing\ \ the\ \ identification $\wide{A}\otimes_{E({A})}E({A})\cong \wide{A}$, we obtain
\[
\sum_i a_i \# h_i = \sum_i a_i \epsilon_H (h_i) \# 1_H .
\]

Therefore $i(\wide{A})= (\underline{\wide{A}\#_{\omega} H})^{co\mathcal{H}}$.

In order to see that the $\mathcal{H}$-extension is cleft, one needs only to define the cleaving map, once item 
(a) of Definition \ref{cleftalgebroid} is automatically satisfied, since $L=R=E({A})$. Define the maps
\[
\begin{array}{rccl} \wide{\gamma} : & \underline{E({A})\# H} & \rightarrow & \underline{\wide{A} \#_{\omega} H} \\
\, & r\# h & \mapsto & r\# h \end{array}
\]
and
\[
\begin{array}{rccl} \ov{\wide{\gamma}} : & \underline{E({A})\# H} & \rightarrow & \underline{\wide{A} \#_{\omega} H} \\
\, & r\# h & \mapsto & r\omega^{-1} (S(h_{(2)}),h_{(3)})\# S(h_{(1)}) \end{array}
\]

Note that $a\# h$ in the domain of $\wide{\gamma}$ means something quite different from $a\# h$ in the image because the first is in the partial smash product, while the second lies in the partial crossed product. It is easy to see that $\wide{\gamma}$ is left $E({A})$-linear. Indeed, for $r\in E({A})$ and $a\# h\ in \underline{E({A}) \# H}$, we have
\begin{eqnarray*}
\wide{\gamma}(s_l (r)(a\# h)) & = & \wide{\gamma}((r\# 1_H )(a\# h))= \wide{\gamma} (ra\# h) \\
& = & ra \# h =(r\# 1_H)(a\# h).
\end{eqnarray*}

Also, the cleaving map is a morphism of right $\mathcal{H}$-comodules. Consider $a\# h \in \underline{E({A}) \# H}$, then
\begin{eqnarray*}
\wide{\rho} \circ \wide{\gamma} (a\# h) & = &   \wide{\rho} (a\# h) = a\# h_{(1)} \otimes 1_A \# h_{(2)} \\
& = & \wide{\gamma}( a\# h_{(1)} \otimes 1_A \# h_{(2)} =(\wide{\gamma}\otimes_{E({A}), \blacktriangleright \blacktriangleleft} Id) \circ \wide{\rho} (a\# h) .
\end{eqnarray*}

Finally, let us check that the maps $\wide{\gamma}$ and $\ov{\wide{\gamma}}$ are mutually inverse by convolution in the sense that
\begin{eqnarray*}
\mu \circ (\wide{\gamma} \otimes_{E({A}), \blacktriangleright \blacktriangleleft} \ov{\wide{\gamma}} )\circ \wide{\Delta}_r & = & i \circ \wide{\epsilon}_l \\
\mu \circ (\ov{\wide{\gamma}} \otimes_{E({A}), \triangleright \triangleleft} \wide{\gamma} )\circ \wide{\Delta}_l & = & i \circ 
\wide{\epsilon}_r .
\end{eqnarray*}

Consider $a\# h \in \underline{E({A})\# H}$, then

\begin{eqnarray*}
& &\hspace{-0.25cm} \mu \circ (\gamma \otimes_{E({A}), \blacktriangleright \blacktriangleleft} \ov{\gamma} )\circ \wide{\Delta}_r (a\# h) 
 =  \wide{\gamma} (a\# h_{(1)})\ov{\wide{\gamma}}(1_A \# h_{(2)}) \\
 &\hspace{-0.3cm} = &\hspace{-0.25cm} (a\# h_{(1)})(\omega^{-1} (S(h_{(3)}), h_{(4)})\# S(h_{(2)})) \\
&\hspace{-0.3cm} = &\hspace{-0.25cm} a (h_{(1)} \cdot  \omega^{-1} (S(h_{(6)}), h_{(7)}))\omega (h_{(2)} , S(h_{(5)})) \# h_{(3)}S(h_{(4)})\\
&\hspace{-0.3cm} = &\hspace{-0.25cm} a (h_{(1)} \cdot  \omega^{-1} (S(h_{(4)}), h_{(5)}))\omega (h_{(2)} , S(h_{(3)})) \# 1_H \\
&\hspace{-0.3cm} = &\hspace{-0.25cm} a \omega^{-1} (h_{(1)}S(h_{(9)}), h_{(10)}) \omega (h_{(2)}, S(h_{(8)})h_{(11)}) \omega^{-1} (h_{(4)}, S(h_{(7)}))
 \omega (h_{(5)} , S(h_{(6)})) \# 1_H \\
&\hspace{-0.3cm} = &\hspace{-0.25cm} a \omega^{-1} (h_{(1)}S(h_{(6)}), h_{(7)}) \omega (h_{(2)}, S(h_{(5)})h_{(8)}) (h_{(3)} \cdot ( S(h_{(4)}) \cdot 1_A)) \# 1_H \\
&\hspace{-0.3cm} = &\hspace{-0.25cm} a \omega^{-1} (h_{(1)}S(h_{(2)}), h_{(7)}) \omega (h_{(3)}, S(h_{(5)})h_{(6)}) (h_{(4)} \cdot 1_A) \# 1_H \\
&\hspace{-0.3cm} = &\hspace{-0.25cm} a(h\cdot 1_A) \# 1_H  =  i (\wide{\epsilon}_l (a\# h)).
\end{eqnarray*}

Also, we have
$$\begin{array}{rcl}
& & \mu \circ (\ov{\wide{\gamma}} \otimes_{E({A}), \triangleright \triangleleft} \wide{\gamma} )\circ \wide{\Delta}_l (a\# h) =  \mu \circ (\ov{\wide{\gamma}} \otimes_{E({A}), \triangleright \triangleleft} \wide{\gamma} )
(a\# h_{(1)}\otimes 1_A \# h_{(2)})\\
& = & \mu \circ (\ov{\wide{\gamma}} \otimes_{E({A}), \triangleright \triangleleft} \wide{\gamma} )
((1_A \# h_{(1)})\triangleleft a \otimes 1_A \# h_{(2)})\\
& = & \mu \circ (\ov{\wide{\gamma}} \otimes_{E({A}), \triangleright \triangleleft} \wide{\gamma} )
(1_A \# h_{(1)} \otimes a\triangleright (1_A \# h_{(2)}))\\
& = & \mu \circ (\ov{\wide{\gamma}} \otimes_{E({A}), \triangleright \triangleleft} \wide{\gamma} )
(1_A \# h_{(1)} \otimes a \# h_{(2)})\\
& = & \ov{\wide{\gamma}}(1_A \# h_{(1)}) \wide{\gamma}(a \# h_{(2)}) \\
& = & (\omega^{-1} (S(h_{(2)}) ,h_{(3)})\# S(h_{(1)}))(a \# h_{(4)})\\
& = & a\omega^{-1} (S(h_{(4)}) ,h_{(5)}) (S(h_{(3)})\cdot a) \omega (S(h_{(2)}) ,h_{(6)}) \# S(h_{(1)}) h_{(7)}\\
& = & a\omega^{-1} (S(h_{(5)}) ,h_{(6)}) \omega (S(h_{(4)}) ,h_{(7)}) (S(h_{(3)})\cdot a)\# S(h_{(1)}) h_{(2)}\\
& = & (S(h_{(2)})\cdot (h_{(3)}\cdot 1_A))(S(h_{(1)})\cdot a)\# 1_H \\
& = & (S(h)\cdot a)\# 1_H 
 =  i (\wide{\epsilon}_r (a\# h)) .
\end{array}$$

Therefore, $\wide{A}\subset \underline{\wide{A}\#_{\omega}H}$ is a $\underline{E({A})\#H}$-cleft extension. \find

\section{Conclusions and outlook}
\begin{itemize}
\item[-] All the cohomology theory done in this paper was done over co-commutative Hopf algebras acting partially 
over commutative algebras. This can be ge\-ne\-ralized for co-commutative Hopf algebra objects and commutative algebra 
objects in braided monoidal categories.
\item[-] One topic of interest is to relate this cohomology for partial actions and the cohomology for its globalization,
then constructing a bridge between this theory and the classical Sweedler's theory.
\item[-] The last theorem placed the notion of a partially cleft extension within the context of cleft extensions for
Hopf algebroids. This suggests, perhaps, that this entire cohomological theory can be understood properly as a 
cohomological theory of Hopf algebroids.
\item[-] Another topic to be explored in further research can be the obstruction theory for the existence of 
partially cleft extensions and its relation with the third cohomology group in the same spirit of \cite{Sch}.
\end{itemize}

\section*{Acknowledgements}

The authors would like to thank Mikhailo Dokuchaev and Mykola Khrypchenko for their precious suggestions and 
comments on this work. We also thank Joost Vercruysse for pointing us some problems in the construction of the Hopf 
algebra $\widetilde{A}$.

\section*{Appendix A. Some additional proofs}

In this appendix, we provide with details the proofs of Section \ref{sec3} that were omitted in the final journal 
version of the paper. We leave them here for convenience of the interested reader.
\vspace{0.3cm}

\noindent\prove{Proposition \ref{prop3.1}}

\rmr{1} $\wide{e}_n(h^1\tens \dots \tens h^n) = 
(\subp{h}{1}\ldots\subp{h}{n})\cdot 1_A$ is an idempotent in 
$Hom_{k}(H^{\tens n},A)$. 

In fact, consider $n\geq 1$ and $(h^1\otimes \dots\otimes h^n) \in \hte$, then
\begin{eqnarray*} & & \widetilde{e}_n*\widetilde{e}_n(\coen)  =  \widetilde{e}_n \prodn{h}{1}{n} \widetilde{e}_n \prodn{h}{2}{n} \\
& = & ((\subp{h}{1}_{(1)}\ldots\subp{h}{n}_{(1)})\cdot 1_A)((\subp{h}{1}_{(2)}\ldots\subp{h}{n}_{(2)})\cdot  1_A)\\
& = & ((\subp{h}{1}\ldots\subp{h}{n})_{(1)}\cdot 1_A)((\subp{h}{1}\ldots\subp{h}{n})_{(2)}\cdot 1_A)\\
& \stackrel{(PA2)}{=} & (\subp{h}{1}\ldots\subp{h}{n})\cdot 1_A
 =  \widetilde{e}_n(\coen).
\end{eqnarray*}

This proves our statement.
\vspace{0.3cm}

\rmr{2} $e_{n,m} = e \tens \underset{m-n}{\underbrace{\epsilon \tens \dots \tens \epsilon}}$ is an 
idempotent in $\mbox{Hom}_k (H^{\otimes m} , A)$.

Take any $n < m$ and $(h^1\tens \dots\tens h^m) \in H^{\tens m}$, then
\begin{eqnarray*} & & e_{n,m}*e_{n,m}(h^1,\ldots, h^m) = e_{n,m}\prodn{h}{1}{m}e_{n,m}\prodn{h}{2}{m}\\
& = & e(\subp{h}{1}_{(1)},\ldots,\subp{h}{n}_{(1)})\epsilon(\subp{h}{n+1}_{(1)})\ldots\epsilon(\subp{h}{m}_{(1)})e(\subp{h}{1}_{(2)},\ldots,\subp{h}{n}_{(2)})\epsilon(\subp{h}{n+1}_{(2)})\ldots
\epsilon(\subp{h}{m}_{(2)})\\
& = & e(\subp{h}{1}_{(1)},\ldots,\subp{h}{n}_{(1)})e(\subp{h}{1}_{(2)},\ldots,\subp{h}{n}_{(2)})\epsilon(\subp{h}{n+1}_{(1)})\ldots\epsilon(\subp{h}{m}_{(1)})\epsilon(\subp{h}{n+1}_{(2)})\ldots
\epsilon(\subp{h}{m}_{(2)})\\
& \stackrel{(*)}{=} & e*e(\subp{h}{1},\ldots,\subp{h}{n})\epsilon(\subp{h}{n+1})\ldots\epsilon(\subp{h}{m})
 =  e(\subp{h}{1},\ldots,\subp{h}{n})\epsilon(\subp{h}{n+1})\ldots\epsilon(\subp{h}{m})\\
& = & e_{n,m}(\subp{h}{1},\ \ldots,\ \subp{h}{m}).
\end{eqnarray*}

The identity $(*)$ follows from $\epsilon(\sub{h}{1})\epsilon(\sub{h}{2}) = \epsilon(\sub{h}{1}\epsilon(\sub{h}{2})) = \epsilon(h)$, and therefore
$e_{n,m}$ is idempotent. \find

\vspace{0.3cm}

\noindent\prove{Proposition \ref{prop3.4}} $e_n (\coen) = \subp{h}{1}\cdot(\subp{h}{2}\cdot(\dots\ \cdot(\subp{h}{n}\cdot 1_A)\dots))$. 

Indeed, take $(\coe) \in  \hte$, then
\begin{eqnarray*} & &\hspace{-0.4cm} e_n(\coen) =  \widetilde{e}_{1,n}*\widetilde{e}_{2,n}*\cdots * \widetilde{e}_{n}(\coen)\\
  & = & \hspace{-0.4cm}  \widetilde{e}_{1,n}\prodn{h}{1}{n}\widetilde{e}_{2,n}\prodn{h}{2}{n}\ldots  \widetilde{e}_{n}\prodn{h}{n}{n}\\
  & = & \hspace{-0.4cm} \widetilde{e}_{1}(h_{(1)}^{1})\epsilon(h_{(1)}^{2})\ldots \epsilon(h_{(1)}^{n})\widetilde{e}_{2}(h_{(2)}^{1},h_{(2)}^{2})\epsilon(h_{(2)}^{3})\ldots \epsilon(h_{(2)}^{n})\\
  &   & \hspace{-0.4cm} \ldots\widetilde{e}_{n-1}(h_{(n-1)}^{1}, \ldots , h_{(n-1)}^{n-1})\epsilon(h_{(n-1)}^{n})\widetilde{e}_{n}(h_{(n)}^{1}, \ldots , h_{(n)}^{n})\\
  & = &\hspace{-0.4cm}   (\sub{h^1}{1}\cdot 1_A)\epsilon(\sub{h^2}{1})\ldots\epsilon(\sub{h^n}{1})(\sub{h^1}{2}\sub{h^2}{2}\cdot 1_A)
\epsilon(\sub{h^3}{2})\ldots\epsilon(\sub{h^n}{2})
   \ldots(\sub{h^1}{n}\ldots\sub{h^n}{n}\cdot 1_A)\\
  & = &\hspace{-0.4cm}  (\sub{h^1}{1}\hspace{-0.1cm}\cdot\hspace{-0.05cm} 1_A)(\sub{h^1}{2}\sub{h^2}{1}\hspace{-0.1cm}\cdot\hspace{-0.05cm} 1_A)\ldots(\sub{h^1}{n-1}\sub{h^2}{n-2}\ldots\sub{h^{n-1}}{1}\hspace{-0.1cm}\cdot\hspace{-0.05cm} 1_A)
    (\sub{h^1}{n}\sub{h^2}{n-1}\ldots\sub{h^{n-1}}{2}h^n\hspace{-0.1cm}\cdot 1_A)\\
    & = & \hspace{-0.4cm}(\sub{h^1}{1}\cdot 1_A)(\sub{h^1}{2}\sub{h^2}{1}\cdot 1_A)\ldots((\sub{h^1}{n-1}\sub{h^2}{n-2}\ldots h^{n-1})_{(1)}\cdot 1_A)\\
  &   & \hspace{-0.4cm}((\sub{h^1}{n-1}\sub{h^2}{n-2}\ldots h^{n-1})_{(2)}h^n\cdot 1_A)\\
  & \stackrel{(PA3)}{=} &\hspace{-0.4cm}  (\sub{h^1}{1}\cdot 1_A)(\sub{h^1}{2}\sub{h^2}{1}\cdot 1_A)\ldots((\sub{h^1}{n-1}\ldots \sub{h^{n-2}}{2} h^{n-1}\cdot (h^n\cdot 1_A))\\
\end{eqnarray*}
  \begin{eqnarray*}  
  & = &\hspace{-0.4cm}  (\sub{h^1}{1}\hspace{-0.1cm}\cdot 1_A)(\sub{h^1}{2}\sub{h^2}{1}\hspace{-0.1cm}\cdot 1_A)\ldots(\sub{h^1}{n-2}\ldots \sub{h^{n-2}}{1}\hspace{-0.1cm} \cdot 1_A)
   (\sub{h^1}{n-1}\ldots \sub{h^{n-2}}{2} h^{n-1}\cdot (h^n\cdot 1_A))\\
   & = & \hspace{-0.4cm} (\sub{h^1}{1}\cdot 1_A)(\sub{h^1}{2}\sub{h^2}{1}\cdot 1_A)\ldots((\sub{h^1}{n-2}\ldots h^{n-3}_{(1)}h^{n-2})_{(1)} \cdot 1_A)\\
  &   & \hspace{-0.4cm} ((\sub{h^1}{n-2}\ldots h^{n-3}_{(2)} h^{n-2})_{(2)} h^{n-1}\cdot (h^n\cdot 1_A))\\
  & \stackrel{(PA3)}{=} &\hspace{-0.4cm} (\sub{h^1}{1}\cdot 1_A)(\sub{h^1}{2}\sub{h^2}{1}\cdot 1_A)\ldots(\sub{h^1}{n-1}\ldots h^{n-3}_{(1)}h^{n-2}\cdot( h^{n-1}\cdot (h^n\cdot 1_A)))\\   
  & = & \hspace{-0.4cm}\dots \;\; = h^1\cdot(h^2 \cdot( \dots\ \cdot(h^{n-1}\cdot (h^n\cdot 1_A))\dots)).
\end{eqnarray*}
in which $(\cdots)$ between the last two equalities means applying repeatedly the process using (PA3) until we obtain the result. 
\find

\vspace{0.3cm}

\noindent\prove{Lemma \ref{p1}} Take $(h^1\tens \cdots\tens h^{n+1}) \in H^{\tens n+1}$, then,

\rmr{i} $E^n(f\ast g) = E^n(f)\ast E^n(g)$, for $f,\ g \in C^n_{par}(H,A)$.
\begin{eqnarray*}
& & E^n(f*g)(h^1, \ldots, \ h^{n+1}) =  h^1\cdot (f*g(h^{2}, \ldots,\ h^{n+1} ))\\
& = & h^1\cdot(f(\sub{h^2}{1}, \ldots, \sub{h^{n+1}}{1})g(\sub{h^2}{2},\ldots, \sub{h^{n+1}}{2}))\\
& = & (\sub{h^1}{1}\cdot(f(\sub{h^2}{1}, \ldots, \sub{h^{n+1}}{1}))(\sub{h^1}{2}\cdot (g(\sub{h^2}{2},\ldots, \sub{h^{n+1}}{2})))\\
& = & E^n(f)(\sub{h^1}{1},\ldots, \ \sub{h^{n+1}}{1})E^n(g)(\sub{h^1}{2},\ldots, \ \sub{h^{n+1}}{2}) .
\end{eqnarray*}

\rmr{ii} $E^n(e_n) = e_{n+1}$.
\begin{eqnarray*}
& & E^n(e_n)(h^1,\dots, h^{n+1})  =  h^1 \cdot e_n(h^2,\dots,h^{n+1})\\
& = & h^1\cdot(h^2\cdot(\dots \cdot(h^{n+1}\cdot 1_A)\dots))\\
& = & e_{n+1}(h^1,\dots, h^{n+1})
\end{eqnarray*}

\rmr{iii} $i_{n,m}(f\ast g)= i_{n,m}(f)\ast i_{n,m}(g)$, for $n<m$ and $f,\ g \in C^n_{par}(H,A)$.
\begin{eqnarray*}
& & i_{n,n+1}(f*g)(h^1,\dots, h^{n+1}) = f*g(h^1,\dots, h^n)\epsilon(h^{n+1})\\
& = & f(\sub{h^1}{1},\dots, \sub{h^n}{1})g(\sub{h^1}{2}, \dots, \sub{h^n}{2})\epsilon(\epsilon(\sub{h^{n+1}}{1})\sub{h^{n+1}}{2}) \\
& = & f(\sub{h^1}{1},\dots, \sub{h^n}{1})\epsilon(\sub{h^{n+1}}{1})g(\sub{h^1}{2}, \dots, \sub{h^n}{2})\epsilon(\sub{h^{n+1}}{2})\\
& = & i_{n,n+1}(f)(\sub{h^1}{1},\dots,\sub{h^{n+1}}{1})i_{n,n+1}(g)(\sub{h^1}{2},\dots,\sub{h^{n+1}}{2})\\
& = & i_{n,n+1}(f)*i_{n,n+1}(g)(h^1,\dots, h^{n+1})
\end{eqnarray*}

\rmr{iv} $i_{n,m}(e_n )\ast e_m =e_m$.
\vspace{0.3cm}

Take $h^1\tens\dots\tens h^m \in H^{\tens m}$, then,
\begin{eqnarray*}
 &\hspace{-0.2cm} &\hspace{-0.35cm} i_{n,m}(e_n )\ast e_m (h^1,\dots, h^m) = i_{n,m}(e_n )(h^{1}_{(1)},\dots,h^{m}_{(1)}) e_m(h^{1}_{(2)},\dots,h^{m}_{(2)})\\
 &\hspace{-0.2cm} = &\hspace{-0.35cm} e_n\hspace{-0.03cm}(\hspace{-0.03cm}h^{1}_{(1)},\dots,h^{n}_{(1)}\hspace{-0.03cm})\varepsilon(\hspace{-0.02cm}h^{n+1}_{(1)}\hspace{-0.02cm})\hspace{-0.02cm}\dots\vare(\hspace{-0.02cm}h^{m}_{(1)}\hspace{-0.02cm})(\hspace{-0.02cm}h^{1}_{(2)}\hspace{-0.09cm}\cdot\hspace{-0.09cm}(\hspace{-0.02cm}\dots\ \hspace{-0.1cm}\cdot\hspace{-0.09cm}(h^n_{(2)}\hspace{-0.1cm}\cdot\hspace{-0.09cm}(h^{n+1}_{(2)}\hspace{-0.1cm}\cdot\hspace{-0.09cm}(\dots\ \hspace{-0.1cm}\cdot\hspace{-0.09cm}(h^{m}_{(2)}\hspace{-0.1cm}\cdot\hspace{-0.09cm} 1_A)\hspace{-0.03cm}\dots\hspace{-0.03cm})))\hspace{-0.03cm}\dots\hspace{-0.03cm}))\\
 &\hspace{-0.2cm} = &\hspace{-0.35cm} (h^{1}_{(1)}\cdot(\dots\ \hspace{-0.1cm} \cdot(h^{n}_{(1)}\cdot 1_A)\dots))(h^{1}_{(2)}\cdot(\dots\ \hspace{-0.1cm} \cdot(h^n_{(2)}\cdot(h^{n+1}\cdot(\dots\ \hspace{-0.1cm}\cdot(h^{m}\cdot 1_A)\dots)))\dots))\\
 &\hspace{-0.2cm} \stackrel{(PA2)}{=} &\hspace{-0.35cm} h^{1}\hspace{-0.1cm}\cdot\hspace{-0.08cm}[(h^2_{(1)}\hspace{-0.1cm}\cdot\hspace{-0.08cm}(\dots\ \hspace{-0.1cm} \cdot\hspace{-0.08cm}(h^{n}_{(1)}\hspace{-0.1cm}\cdot\hspace{-0.08cm} 1_A)\dots))(h^{2}_{(2)}\hspace{-0.1cm}\cdot\hspace{-0.08cm}(\dots\ \hspace{-0.1cm} \cdot\hspace{-0.08cm}(h^n_{(2)}\hspace{-0.1cm}\cdot\hspace{-0.08cm}(h^{n+1}\hspace{-0.1cm}\cdot\hspace{-0.08cm}(\dots\ \hspace{-0.1cm}\cdot\hspace{-0.08cm}(h^{m}\hspace{-0.1cm}\cdot\hspace{-0.08cm} 1_A)\dots)))\dots))]\\
 &\hspace{-0.2cm} \stackrel{(PA2)}{=} &\hspace{-0.35cm}  h^{1}\cdot(h^2\cdot(h^3\cdot(\dots\ \cdot(h^{n}\cdot( 1_A(h^{n+1}\cdot(\dots\ \cdot(h^{m}\cdot 1_A)\dots))))\dots)))\\
 &\hspace{-0.2cm} = &\hspace{-0.35cm} h^{1}\cdot(h^2\cdot(\dots\ \cdot(h^{n}\cdot( h^{n+1}\cdot(\dots\ \cdot(h^{m}\cdot 1_A)\dots)))\dots))\\
 &\hspace{-0.2cm} = &\hspace{-0.35cm} e_m(h^1,\ldots, h^m)
 \end{eqnarray*}

\rmr{v} $(f*g)\circ\mu_i = (f\circ\mu_i)*(g\circ\mu_i)$, for $f,\ g \in C^n_{par}(H,A)$ and $\forall \ i\in \{1,\ldots, n\}$.
\vspace{0.3cm}

In fact, for all $i\in \{1,\dots, n\}$ and for $f,\ g \in C^n_{par}(H,A)$
\begin{eqnarray*}
   & & (f*g)\circ \mu_i (h^1, \dots, h^i, h^{i+1}, \dots,  h^{n+1}) = (f*g)(h^1, \dots, h^i h^{i+1}, \dots,  h^{n+1})\\
   & = & f(h^1_{(1)}, \dots, (h^i h^{i+1})_{(1)}, \dots,  h^{n+1}_{(1)})g(h^1_{(2)}, \dots, (h^i h^{i+1})_{(2)}, \dots,  h^{n+1}_{(2)})\\
   & = & f(h^1_{(1)}, \dots, h^i_{(1)} h^{i+1}_{(1)}, \dots,  h^{n+1}_{(1)})g(h^1_{(2)}, \dots, h^i_{(2)} h^{i+1}_{(2)}, \dots,  h^{n+1}_{(2)})\\
   & = & f\circ \mu_i(h^1_{(1)}, \dots, h^i_{(1)}, h^{i+1}_{(1)}, \dots,  h^{n+1}_{(1)})g\circ\mu_i(h^1_{(2)}, \dots, h^i_{(2)}, h^{i+1}_{(2)}, \dots, h^{n+1}_{(2)})\\
   & = &  (f\circ\mu_i)*(g\circ\mu_i)(h^1, \dots, h^i, h^{i+1}, \dots,  h^{n+1})
  \end{eqnarray*}

\vspace{0.3cm}
\rmr{vi} $(e_n\circ\mu_n)*i_{n,n+1}(e_n)= e_{n+1}$.
\begin{eqnarray*}
   & &(e_n\circ \mu_n)*i_{n,n+1}(e_n) (h^1, \dots,  h^{n+1}) =\\
   & = & e_n\circ\mu_n (h^1_{(1)}, \dots, h^{n}_{(1)},  h^{n+1}_{(1)})i_{n,n+1}(e_n)(h^1_{(2)}, \dots, h^{n}_{(2)},  h^{n+1}_{(2)})\\
   & = & e_n(h^1_{(1)}, \dots, h^{n}_{(1)}  h^{n+1}_{(1)})e_n(h^1_{(2)}, \dots, h^{n}_{(2)})\varepsilon( h^{n+1}_{(2)})\\
   & = & (h^1_{(1)}\cdot(h^2_{(1)}\cdot( ... \cdot( h^{n}_{(1)}  h^{n+1}\cdot 1_A)...)))(h^1_{(2)}\cdot(h^2_{(2)}\cdot( ... \cdot( h^{n}_{(2)}\cdot 1_A)...)))\\ 
   & \stackrel{(PA2)}{=} & h^1\cdot[(h^2_{(1)}\cdot( ... \cdot( h^{n}_{(1)}  h^{n+1}\cdot 1_A)...))(h^2_{(2)}\cdot( ... \cdot( h^{n}_{(2)}\cdot 1_A)...))]\\
   & = & h^1\cdot(h^2\cdot( ... \cdot(h^{n-1}\cdot[ (h^{n}_{(1)}  h^{n+1}\cdot 1_A)( h^{n}_{(2)}\cdot 1_A)])...))\\ 
   & \stackrel{(PA3)}{=} & h^1\cdot(h^2\cdot( ... \cdot(h^{n-1}\cdot (h^{n}\cdot(  h^{n+1}\cdot 1_A)))...))\\
   & = & e_{n+1}(h^{1},\dots,h^{n+1})  
    \end{eqnarray*}

\vspace{0.3cm}
\rmr{vii} $(e_n\circ\mu_i)*e_{n+1} = e_{n+1}$, $\forall\ i\in\{1,\cdots, n-1\}$.
\begin{eqnarray*}
  &\hspace{-0.3cm} &  \hspace{-0.4cm} (e_n\circ \mu_i)*e_{n+1} (h^1, \dots , h^i, h^{i+1}, \dots,  h^{n+1}) =\\ 
  &\hspace{-0.3cm} = &\hspace{-0.4cm} e_n\circ\mu_i (h^1_{(1)}, \dots, h^i_{(1)}, h^{i+1}_{(1)}, \dots,  h^{n+1}_{(1)})e_{n+1}(h^1_{(2)}, \dots, h^{i}_{(2)},  h^{i+1}_{(2)}, \dots, h^{n+1}_{(2)})\\
  &\hspace{-0.3cm} = &\hspace{-0.4cm} e_n(h^1_{(1)}, \dots, h^{i}_{(1)}  h^{i+1}_{(1)},\dots, h^{n+1}_{(1)})(h^1_{(2)}\cdot( ... \cdot ( h^{i}_{(2)} \cdot ( h^{i+1}_{(2)} \cdot (... \cdot( h^{n+1}_{(2)}\cdot 1_A)...)))...))\\
  &\hspace{-0.3cm} = &\hspace{-0.4cm} (h^1_{(1)}\hspace{-0.1cm}\cdot\hspace{-0.08cm}(... \hspace{-0.1cm} \cdot \hspace{-0.08cm}(\hspace{-0.02cm} h^{i}_{(1)}  h^{i+1}_{(1)}\hspace{-0.1cm}\cdot\hspace{-0.08cm}(...\hspace{-0.1cm}\cdot\hspace{-0.08cm}(\hspace{-0.02cm} h^{n+1}_{(1)}\hspace{-0.1cm}\cdot\hspace{-0.08cm} 1_A)...)\hspace{-0.03cm})...)\hspace{-0.03cm})(\hspace{-0.02cm}h^1_{(2)}\hspace{-0.1cm}\cdot\hspace{-0.08cm}( ...\hspace{-0.1cm} \cdot\hspace{-0.08cm} (\hspace{-0.02cm} h^{i}_{(2)}\hspace{-0.1cm} \cdot \hspace{-0.08cm}(\hspace{-0.02cm} h^{i+1}_{(2)}\hspace{-0.1cm} \cdot\hspace{-0.08cm} (... \hspace{-0.1cm}\cdot\hspace{-0.08cm}(\hspace{-0.02cm} h^{n+1}_{(2)}\hspace{-0.1cm}\cdot\hspace{-0.08cm} 1_A)...)\hspace{-0.03cm})\hspace{-0.03cm})...)\hspace{-0.03cm})\\ 
  &\hspace{-0.3cm} \stackrel{(PA2)}{=} &\hspace{-0.4cm} h^1\hspace{-0.1cm}\cdot\hspace{-0.1cm}[\hspace{-0.02cm}(\hspace{-0.03cm}h^2_{(1)}\hspace{-0.1cm}\cdot\hspace{-0.1cm}(\hspace{-0.02cm} ... \hspace{-0.1cm} \cdot\hspace{-0.1cm}(\hspace{-0.03cm} h^{i}_{(1)}  h^{i+1}_{(1)}\hspace{-0.1cm}\cdot\hspace{-0.1cm}(\hspace{-0.02cm}...\hspace{-0.1cm}\cdot\hspace{-0.1cm}(\hspace{-0.03cm} h^{n+1}_{(1)}\hspace{-0.1cm}\cdot\hspace{-0.1cm} 1_A\hspace{-0.02cm})...\hspace{-0.02cm})\hspace{-0.05cm})\hspace{-0.02cm}...\hspace{-0.02cm})\hspace{-0.05cm})
 (\hspace{-0.02cm}h^2_{(2)}\hspace{-0.1cm}\cdot\hspace{-0.1cm}(\hspace{-0.02cm} ... \hspace{-0.1cm}\cdot\hspace{-0.1cm} (\hspace{-0.02cm} h^{i}_{(2)}\hspace{-0.1cm} \cdot\hspace{-0.1cm} (\hspace{-0.02cm} h^{i+1}_{(2)}\hspace{-0.1cm} \cdot \hspace{-0.1cm} (\hspace{-0.02cm}... \hspace{-0.1cm} \cdot\hspace{-0.1cm}(\hspace{-0.02cm} h^{n+1}_{(2)}\hspace{-0.1cm}\cdot\hspace{-0.1cm} 1_A\hspace{-0.02cm})...\hspace{-0.02cm})\hspace{-0.05cm})\hspace{-0.05cm})...\hspace{-0.02cm})\hspace{-0.05cm})\hspace{-0.05cm}]\\
  &\hspace{-0.3cm} = &\hspace{-0.4cm} h^1\hspace{-0.09cm}\cdot\hspace{-0.09cm}(\hspace{-0.01cm} ...\hspace{-0.09cm} \cdot\hspace{-0.09cm}(\hspace{-0.02cm}h^{i-1}\hspace{-0.09cm}\cdot\hspace{-0.09cm}[(\hspace{-0.02cm}h^{i}_{(1)}  h^{i+1}_{(1)}\hspace{-0.09cm}\cdot\hspace{-0.09cm}(h^{i+2}_{(1)}\hspace{-0.09cm}\cdot\hspace{-0.09cm}(...\hspace{-0.09cm}\cdot\hspace{-0.09cm}(\hspace{-0.02cm} h^{n+1}_{(1)}\hspace{-0.09cm}\cdot\hspace{-0.09cm} 1_A)...)\hspace{-0.03cm})\hspace{-0.03cm})(\hspace{-0.02cm} h^{i}_{(2)}\hspace{-0.09cm} \cdot \hspace{-0.09cm}( h^{i+1}_{(2)}\hspace{-0.09cm} \cdot\hspace{-0.09cm} (... \hspace{-0.09cm}\cdot\hspace{-0.09cm}(\hspace{-0.02cm} h^{n+1}_{(2)}\hspace{-0.09cm}\cdot\hspace{-0.09cm} 1_A)...\hspace{-0.03cm})\hspace{-0.03cm})\hspace{-0.03cm})\hspace{-0.03cm}]\hspace{-0.03cm})...) \\ 
  &\hspace{-0.3cm} \stackrel{(PA3)}{=} &\hspace{-0.4cm} h^1\cdot( ... \cdot(h^{i-1}\cdot[(h^{i}_{(1)}  h^{i+1}_{(1)}\cdot(h^{i+2}_{(1)}\cdot(...\cdot( h^{n+1}_{(1)}\cdot 1_A)...)))\\
  & & \hspace{-0.3cm} (h^{i}_{(2)}h^{i+1}_{(2)} \cdot(h^{i+2}_{(2)}\cdot( ... \cdot( h^{n+1}_{(2)}\cdot 1_A)...)))(h^i_{(3)}\cdot 1_A)])...) \\ 
  &\hspace{-0.3cm} \stackrel{(PA2)}{=} &\hspace{-0.4cm} h^1\hspace{-0.1cm}\cdot\hspace{-0.1cm}(\hspace{-0.03cm} ...\hspace{-0.1cm} \cdot\hspace{-0.1cm}(\hspace{-0.03cm}h^{i-1}\hspace{-0.1cm}\cdot\hspace{-0.1cm}[\hspace{-0.03cm}(\hspace{-0.03cm}h^{i}_{(1)}  h^{i+1}\hspace{-0.1cm}\cdot\hspace{-0.1cm}[\hspace{-0.03cm}(\hspace{-0.03cm}h^{i+2}_{(1)}\hspace{-0.1cm}\cdot\hspace{-0.1cm}(\hspace{-0.03cm}...\hspace{-0.1cm}\cdot\hspace{-0.1cm}(\hspace{-0.03cm} h^{n+1}_{(1)}\hspace{-0.1cm}\cdot\hspace{-0.1cm} 1_A\hspace{-0.03cm})...\hspace{-0.03cm})\hspace{-0.03cm})(\hspace{-0.03cm}h^{i+2}_{(2)}\hspace{-0.1cm}\cdot\hspace{-0.1cm}(\hspace{-0.03cm} ...\hspace{-0.1cm} \cdot\hspace{-0.1cm}(\hspace{-0.03cm} h^{n+1}_{(2)}\hspace{-0.1cm}\cdot\hspace{-0.1cm} 1_A\hspace{-0.03cm})...\hspace{-0.03cm})\hspace{-0.03cm})\hspace{-0.03cm}]\hspace{-0.03cm})(\hspace{-0.03cm}h^i_{(2)}\hspace{-0.1cm}\cdot\hspace{-0.1cm} 1_A\hspace{-0.03cm})\hspace{-0.03cm}]\hspace{-0.03cm})...\hspace{-0.03cm}) \\   
  &\hspace{-0.3cm} \stackrel{(PA2)}{=} &\hspace{-0.4cm} h^1\cdot( ... \cdot(h^{i-1}\cdot[(h^{i}_{(1)}  h^{i+1}\cdot[h^{i+2}\cdot(...\cdot( h^{n+1}\cdot 1_A)...)])(h^i_{(2)}\cdot 1_A)])...) \\   
  &\hspace{-0.3cm} \stackrel{(PA3)}{=} &\hspace{-0.4cm} h^1\cdot( ... \cdot(h^{i-1}\cdot(h^{i}\cdot(  h^{i+1}\cdot(h^{i+2}\cdot(...\cdot( h^{n+1}\cdot 1_A)...)))))...) \\
  &\hspace{-0.3cm} = &\hspace{-0.4cm} e_{n+1}(h^1,\dots, h^{n+1})
   \end{eqnarray*}

\vspace{0.3cm}
\noindent\prove{Lemma \ref{p2}} Let $f\in C^{n-1}_{par} (H,A)$ and $h^1\tens\dots\tens h^{n+1} \in H^{\tens n+1}$.

\rmr{i} $E^n(i_{n-1,n}(f)) = i_{n,n+1}(E^{n-1}(f))$.
\begin{eqnarray*} & & E^n(i_{n-1,n}(f))(h^1, \dots, \ h^{n+1}) = h^1\cdot(i_{n-1,n}(f)(h^2, \dots, h^{n+1}))\\
& = & h^1\cdot(f(h^2, \dots, \ h^n)\epsilon(h^{n+1}))
 =  (h^1\cdot f(h^2, \dots, \ h^n))\epsilon(h^{n+1})\\
& = & (E^{n-1}(f))(h^1,\dots, h^{n})\epsilon(h^{n+1})
 =  i_{n,n+1}(E^{n-1}(f))(h^1,\dots, h^{n+1})
\end{eqnarray*}

\rmr{ii} $(e_{n-1}\circ\mu_i\circ\mu_{i+1})*e_{n+1} = e_{n+1}$, $\forall\ i\in\{1,\cdots, n-1\}$.
\vspace{0.3cm}

For every $i\in \{1,\ldots, n-1\}$,
\begin{eqnarray*}
& & \hspace{-0.3cm}  (e_{n-1}\circ \mu_i\circ \mu_{i+1})*e_{n+1} (h^1, \ldots ,  h^{n+1}) =\\
& = &\hspace{-0.3cm} e_{n-1}(h^1_{(1)}, \ldots, h^i_{(1)} h^{i+1}_{(1)} h^{i+2}_{(1)}, \ldots,  h^{n+1}_{(1)})e_{n+1}(h^1_{(2)}, \ldots, h^{(n+1)}_{(2)})\\
& = &\hspace{-0.3cm}  h^1_{(1)}\cdot(... \cdot( h^{i}_{(1)} h^{i+1}_{(1)}h^{i+2}_{(1)}\cdot (... \cdot(h^{n+1}_{(1)}\cdot 1_A)...))...)\\
& & \hspace{-0.3cm}  (h^1_{(2)}\cdot ( ... \cdot ( h^{i}_{(2)} \cdot ( h^{i+1}_{(2)} \cdot (h^{i+2}_{(2)}\cdot (... \cdot( h^{n+1}_{(2)}\cdot 1_A)...))))...)\\
& \stackrel{(PA2)}{=} &\hspace{-0.3cm}  h^1\cdot (h^2\cdot (... \cdot[( h^{i}_{(1)} h^{i+1}_{(1)}h^{i+2}_{(1)}\cdot (... \cdot(h^{n+1}_{(1)}\cdot 1_A)...)) \\
& & \hspace{-0.3cm} ( h^{i}_{(2)} \cdot ( h^{i+1}_{(2)} \cdot (h^{i+2}_{(2)}\cdot (... \cdot ( h^{n+1}_{(2)}\cdot 1_A)...))))]...))\\ 
& \stackrel{(PA3)}{=} &\hspace{-0.3cm}  h^1\cdot (h^2\cdot (... \cdot [( h^{i}_{(1)} h^{i+1}_{(1)}h^{i+2}_{(1)}\cdot (... \cdot (h^{n+1}_{(1)}\cdot 1_A)...)) \\
& & \hspace{-0.3cm} ( h^{i}_{(2)} h^{i+1}_{(2)} \cdot (h^{i+2}_{(2)}\cdot (... \cdot ( h^{n+1}_{(2)}\cdot 1_A)...)))(h^{i}_{(3)}\cdot 1_A)]...))\\
& \stackrel{(PA3)}{=} &\hspace{-0.3cm}  h^1\cdot (h^2\cdot (... \cdot [( h^{i}_{(1)} h^{i+1}_{(1)}h^{i+2}_{(1)}\cdot (... \cdot (h^{n+1}_{(1)}\cdot 1_A)...)) \\
& & \hspace{-0.3cm} ( h^{i}_{(2)} h^{i+1}_{(2)} h^{i+2}_{(2)}\cdot (... \cdot( h^{n+1}_{(2)}\cdot 1_A)...)) (h^{i}_{(3)}h^{i+1}_{(3)}\cdot1_A )(h^{i}_{(4)}\cdot 1_A )]...))\\
& \stackrel{(PA3)}{=} &\hspace{-0.3cm}  h^1\cdot(h^2\cdot(... \cdot[( h^{i}_{(1)} h^{i+1}_{(1)}h^{i+2}\cdot (... \cdot(h^{n+1}\cdot 1_A)...))
 (h^{i}_{(2)}h^{i+1}_{(2)}\cdot 1_A)(h^{i}_{(3)}\cdot 1_A)]...))\\
& \stackrel{(PA3)}{=} &\hspace{-0.3cm}  h^1\cdot(h^2\cdot(... \cdot[( h^{i}_{(1)} h^{i+1}\cdot(h^{i+2}\cdot (... \cdot(h^{n+1}\cdot 1_A)...)))(h^{i}_{(2)}\cdot 1_A)]...))\\
& = &\hspace{-0.3cm}  h^1\cdot(h^2\cdot(... \cdot( h^{i}\cdot( h^{i+1}\cdot(h^{i+2}\cdot (... \cdot(h^{n+1}\cdot 1_A)...))))...))\\
& = &\hspace{-0.3cm} e_{n+1}(h^1,\ldots, h^{n+1}) . 
\end{eqnarray*}

\rmr{iii} $(e_{n-1}\circ\mu_i\circ\mu_{i+j})*e_{n+1} = e_{n+1}$, $\forall \ i\in\{1,\cdots, n-1\},\ j\in \{2,\cdots, n-i\}$.
\vspace{0.3cm}

For $i\in\{1,\cdots, n-1\},\ j\in \{2,\cdots, n-i\}$, we have
\begin{eqnarray*}
  &\hspace{-0.3cm} &\hspace{-0.35cm}  (e_{n-1}\circ \mu_i \circ \mu_{i+j})*e_{n+1} (h^1, \ldots,  h^{n+1}) =\\
  &\hspace{-0.3cm} = &\hspace{-0.35cm}  e_{n-1}(h^1_{(1)}, \ldots, h^i_{(1)} h^{i+1}_{(1)}, \ldots, h^{i+j}_{(1)} h^{i+j+1}_{(1)}, \ldots,  h^{n+1}_{(1)})\\
  & & e_{n+1}(h^1_{(2)}, \ldots, h^{i}_{(2)},  h^{i+1}_{(2)},\ldots, h^{i+j}_{(2)},  h^{i+j+1}_{(2)}, \ldots, h^{n+1}_{(2)})\\
  &\hspace{-0.3cm}= &\hspace{-0.35cm}  h^1_{(1)}\cdot( ... \cdot( h^{i}_{(1)}  h^{i+1}_{(1)}\cdot(... \cdot(h^{i+j}_{(1)}  h^{i+j+1}_{(1)}\cdot(... \cdot( h^{n+1}_{(1)}\cdot 1_A)...))...))...)\\
  &\hspace{-0.35cm} &\hspace{-0.35cm} h^1_{(2)}\cdot( ... \cdot ( h^{i}_{(2)} \cdot ( h^{i+1}_{(2)} \cdot (... \cdot( h^{i+j}_{(2)} \cdot ( h^{i+j+1}_{(2)}\cdot(...\cdot( h^{n+1}_{(2)}\cdot 1_A)...)))...)))...)\\
 &\hspace{-0.3cm}\stackrel{(PA2)}{=} &\hspace{-0.35cm}  h^1\cdot[(h^2_{(1)}\cdot( ...  \cdot( h^{i}_{(1)}  h^{i+1}_{(1)}\cdot(... \cdot(h^{i+j}_{(1)}  h^{i+j+1}_{(1)}\cdot(...\cdot( h^{n+1}_{(1)}\cdot 1_A)...))...))...))\\
  &\hspace{-0.3cm} &\hspace{-0.35cm} (h^2_{(2)}\cdot( ... \cdot ( h^{i}_{(2)} \cdot ( h^{i+1}_{(2)}\cdot(... \cdot(h^{i+j}_{(2)}\cdot(  h^{i+j+1}_{(2)} \cdot (... \cdot( h^{n+1}_{(2)}\cdot 1_A)...)))...)))...))]\\
 \end{eqnarray*}

\begin{eqnarray*} 
  &\hspace{-0.3cm}= &\hspace{-0.35cm} h^1\cdot( ... \cdot(h^{i-1}\cdot[(h^{i}_{(1)}  h^{i+1}_{(1)}\cdot(... \cdot(h^{i+j}_{(1)}  h^{i+j+1}_{(1)} \cdot(...\cdot( h^{n+1}_{(1)}\cdot 1_A)...))...))\\
  &\hspace{-0.3cm} &\hspace{-0.35cm} ( h^{i}_{(2)} \cdot ( h^{i+1}_{(2)}\cdot(... \cdot(h^{i+j}_{(2)}\cdot(  h^{i+j+1}_{(2)} \cdot (... \cdot( h^{n+1}_{(2)}\cdot 1_A)...)))...)))])...) \\ 
 &\hspace{-0.3cm}\stackrel{(PA3)}{=} & \hspace{-0.35cm} h^1\cdot( ... \cdot(h^{i-1}\cdot[(h^{i}_{(1)}  h^{i+1}_{(1)}\cdot(h^{i+2}_{(1)}\cdot(... \cdot(h^{i+j}_{(1)}  h^{i+j+1}_{(1)}\cdot(...
  \cdot( h^{n+1}_{(1)}\cdot 1_A)...))...)))\\
 &\hspace{-0.3cm} &\hspace{-0.35cm}( h^{i}_{(2)} h^{i+1}_{(2)}\cdot(h^{i+2}_{(2)}\cdot(... \cdot(h^{i+j}_{(2)}\cdot(  h^{i+j+1}_{(2)} \cdot (...
  \cdot( h^{n+1}_{(2)}\cdot 1_A)...)))...)))(h^i_{(3)}\cdot1_A)])...) \\ 
 &\hspace{-0.3cm}\stackrel{(PA2)}{=} &\hspace{-0.35cm} h^1\cdot( ... \cdot(h^{i-1}\cdot[(h^{i}_{(1)}  h^{i+1}\cdot[(h^{i+2}_{(1)}\cdot(... \cdot(h^{i+j}_{(1)}  h^{i+j+1}_{(1)}\cdot(...
  \hspace{-0.05cm}\cdot\hspace{-0.05cm}( h^{n+1}_{(1)}\hspace{-0.05cm}\cdot\hspace{-0.05cm} 1_A)...))...))\\
 &\hspace{-0.3cm} &\hspace{-0.35cm} ( h^{i+2}_{(2)}\cdot(... \cdot(h^{i+j}_{(2)}\cdot(  h^{i+j+1}_{(2)} \cdot (... \cdot( h^{n+1}_{(2)}\cdot 1_A)...)))
   ...))])(h^i_{(2)}\cdot1_A)])...) \\ 
 &\hspace{-0.3cm}\stackrel{(PA2)}{=} &\hspace{-0.35cm}  h^1\hspace{-0.05cm}\cdot\hspace{-0.05cm}( ... \hspace{-0.05cm}\cdot\hspace{-0.05cm}(h^{i-1}\hspace{-0.05cm}\cdot\hspace{-0.05cm}[(h^{i}_{(1)}  h^{i+1}\hspace{-0.05cm}\cdot\hspace{-0.05cm}[(h^{i+2}\hspace{-0.05cm}\cdot\hspace{-0.05cm}(... \hspace{-0.05cm}\cdot\hspace{-0.05cm} (h^{i+j-1}\hspace{-0.05cm}\cdot\hspace{-0.05cm}[(h^{i+j}_{(1)}  h^{i+j+1}_{(1)}\hspace{-0.05cm}\cdot\hspace{-0.05cm}(...
   \cdot( h^{n+1}_{(1)}\cdot 1_A)...))\\
 &\hspace{-0.3cm} &\hspace{-0.35cm} (h^{i+j}_{(2)}\cdot(  h^{i+j+1}_{(2)} \cdot (... \cdot( h^{n+1}_{(2)}\cdot 1_A)...)))])...))])
  (h^i_{(2)}\cdot1_A)])...) \\ 
 &\hspace{-0.3cm}\stackrel{(PA3)}{=} &\hspace{-0.35cm}  h^1\hspace{-0.05cm}\cdot\hspace{-0.05cm}( ... \hspace{-0.05cm}\cdot\hspace{-0.05cm}(h^{i-1}\hspace{-0.05cm}\cdot\hspace{-0.05cm}[(h^{i}_{(1)}  h^{i+1}\hspace{-0.05cm}\cdot\hspace{-0.05cm}[... \hspace{-0.05cm}\cdot\hspace{-0.05cm} (h^{i+j-1}\hspace{-0.05cm}\cdot\hspace{-0.05cm}[(h^{i+j}_{(1)}  h^{i+j+1}_{(1)}\hspace{-0.05cm}\cdot\hspace{-0.05cm}(h^{i+j+2}_{(1)}\hspace{-0.05cm}\cdot\hspace{-0.05cm} 
(... \hspace{-0.05cm}\cdot\hspace{-0.05cm}( h^{n+1}_{(1)}\hspace{-0.05cm}\cdot\hspace{-0.05cm} 1_A)...)))\\
 &\hspace{-0.3cm} &\hspace{-0.35cm} (h^{i+j}_{(2)}h^{i+j+1}_{(2)} \cdot(h^{i+j+2}_{(2)}\cdot (... \cdot( h^{n+1}_{(2)}\cdot 1_A)...)))
  (h^{i+j}_{(3)}\cdot1_A)])...])(h^i_{(2)}\cdot1_A)])...) \\ 
 &\hspace{-0.3cm} \stackrel{(PA2)}{=} &\hspace{-0.35cm}  h^1\hspace{-0.05cm}\cdot\hspace{-0.05cm}( ... \hspace{-0.05cm}\cdot\hspace{-0.05cm}(h^{i-1}\hspace{-0.05cm}\cdot\hspace{-0.05cm}[(h^{i}_{(1)}  h^{i+1}\hspace{-0.05cm}\cdot\hspace{-0.05cm}[... \hspace{-0.05cm}\cdot\hspace{-0.05cm} (h^{i+j-1}\hspace{-0.05cm}\cdot\hspace{-0.05cm}[(h^{i+j}_{(1)}  h^{i+j+1}\hspace{-0.05cm}\cdot\hspace{-0.05cm}[(h^{i+j+2}_{(1)}\hspace{-0.05cm}\cdot\hspace{-0.05cm}(...
  \hspace{-0.05cm}\cdot\hspace{-0.05cm}( h^{n+1}_{(1)}\hspace{-0.05cm}\cdot\hspace{-0.05cm} 1_A)...))\\
 &\hspace{-0.3cm} &\hspace{-0.35cm} (h^{i+j+2}_{(2)}\cdot (... \cdot( h^{n+1}_{(2)}\cdot 1_A)...))])(h^{i+j}_{(2)}\cdot1_A)])...])
  (h^i_{(2)}\cdot1_A)])...) \\ 
 &\hspace{-0.3cm}\stackrel{(PA2)}{=} &\hspace{-0.35cm}  h^1\hspace{-0.05cm}\cdot\hspace{-0.05cm}( ... \hspace{-0.05cm}\cdot\hspace{-0.05cm}(h^{i-1}\hspace{-0.05cm}\cdot\hspace{-0.05cm}[(h^{i}_{(1)}  h^{i+1}\hspace{-0.05cm}\cdot\hspace{-0.05cm}[... \hspace{-0.05cm}\cdot\hspace{-0.05cm} (h^{i+j-1}\hspace{-0.05cm}\cdot\hspace{-0.05cm}[(h^{i+j}_{(1)}  h^{i+j+1}\hspace{-0.05cm}\cdot\hspace{-0.05cm}[(h^{i+j+2}\hspace{-0.05cm}\cdot\hspace{-0.05cm}(...
  \hspace{-0.05cm}\cdot\hspace{-0.05cm}( h^{n+1}\hspace{-0.05cm}\cdot\hspace{-0.05cm} 1_A)...))])\\
 &\hspace{-0.3cm} &\hspace{-0.35cm} (h^{i+j}_{(2)}\cdot1_A)]...])(h^i_{(2)}\cdot1_A)])...) \\ 
 &\hspace{-0.3cm}\stackrel{(PA3)}{=} &\hspace{-0.35cm}  h^1\hspace{-0.1cm}\cdot\hspace{-0.1cm}(\hspace{-0.03cm} ... \hspace{-0.1cm}\cdot\hspace{-0.1cm}(\hspace{-0.03cm}h^{i-1}\hspace{-0.1cm}\cdot\hspace{-0.1cm}[\hspace{-0.03cm}(\hspace{-0.03cm}h^{i}_{(1)}  h^{i+1}\hspace{-0.1cm}\cdot\hspace{-0.1cm}[... \hspace{-0.1cm}\cdot\hspace{-0.1cm}[\hspace{-0.03cm}(\hspace{-0.03cm}h^{i+j}\hspace{-0.1cm}\cdot\hspace{-0.1cm}(\hspace{-0.03cm} h^{i+j+1}\hspace{-0.1cm}\cdot\hspace{-0.1cm}(\hspace{-0.03cm}h^{i+j+2}\hspace{-0.1cm}\cdot\hspace{-0.1cm} (... 
  \hspace{-0.1cm}\cdot\hspace{-0.1cm}(\hspace{-0.03cm} h^{n+1}\hspace{-0.1cm}\cdot\hspace{-0.1cm} 1_A\hspace{-0.03cm})...\hspace{-0.03cm})\hspace{-0.03cm})\hspace{-0.03cm})\hspace{-0.03cm})\hspace{-0.03cm}]...])
   (h^i_{(2)}\hspace{-0.1cm}\cdot\hspace{-0.1cm} 1_A\hspace{-0.03cm})\hspace{-0.03cm}]\hspace{-0.03cm})...\hspace{-0.03cm}) \\ 
 &\hspace{-0.3cm}\stackrel{(PA3)}{=} &\hspace{-0.35cm}  h^1\hspace{-0.07cm}\cdot\hspace{-0.05cm}( ... \hspace{-0.07cm}\cdot\hspace{-0.05cm}(h^{i-1}\hspace{-0.07cm}\cdot\hspace{-0.05cm}[(h^{i}\hspace{-0.07cm}\cdot\hspace{-0.05cm}(  h^{i+1}\hspace{-0.07cm}\cdot\hspace{-0.05cm}(... \hspace{-0.07cm}\cdot\hspace{-0.05cm}(h^{i+j}\hspace{-0.07cm}\cdot\hspace{-0.05cm}(  h^{i+j+1}\hspace{-0.07cm}\cdot\hspace{-0.05cm}(h^{i+j+2}\hspace{-0.07cm}\cdot\hspace{-0.05cm}(...
  \hspace{-0.07cm}\cdot\hspace{-0.05cm}( h^{n+1}\hspace{-0.07cm}\cdot\hspace{-0.05cm} 1_A)...)\hspace{-0.03cm})\hspace{-0.03cm})\hspace{-0.03cm})\hspace{-0.03cm})...)\hspace{-0.03cm})\hspace{-0.03cm})\hspace{-0.03cm}]\hspace{-0.03cm})...) \\ 
 &\hspace{-0.3cm} = &\hspace{-0.35cm}   e_{n+1}(h^1,\ldots, h^{n+1}).
\end{eqnarray*}

\rmr{iv} $E^n(f\circ \mu_i) = E^{n-1}(f)\circ \mu_{i+1}$, $\forall\ i\in\{1,\ldots, n-1\}$.
\vspace{0.3cm}

In fact, for every $i\in \{1,\dots, n-1\}$, we have  
\begin{eqnarray*} & & E^n(f\circ \mu_i)(h^1,\dots,h^{n+1})  =  h^1\hspace{-0.08cm}\cdot\hspace{-0.06cm}((f\circ \mu_i)(h^2,\dots, h^i,\hspace{-0.15cm} \underset{i-esima \ coord.}{\underbrace{h^{i+1}}}\hspace{-0.15cm}, h^{i+2}, \dots, h^{n+1})\hspace{-0.03cm})\\
          & = & h^1 \cdot(f(h^2,\dots, h^i, h^{i+1}h^{i+2}, \dots, h^{n+1}))
         \end{eqnarray*}

On the other hand,
 \begin{eqnarray*}& & E^{n-1}(f)\circ \mu_{i+1}(h^1,\dots,h^{n+1})  =  E^{n-1}(f)(h^1 ,\dots, h^i, h^{i+1} h^{i+2}, \dots, h^{n+1})\\
          & = & h^1 \cdot(f(h^2,\dots, h^i, h^{i+1}h^{i+2}, \dots, h^{n+1}))
        \end{eqnarray*}

         So, we prove the equality.

\vspace{0.3cm}
\rmr{v} $E^n\circ E^{n-1}(f) = i_{1,n+1}(\widetilde{e}_1)*(E^{n-1}(f)\circ\mu_1)$.
\begin{eqnarray*} & & i_{1,n+1}(\widetilde{e}_1)*(E^{n-1}(f)\circ\mu_1)(h^1,\dots,h^{n+1})=\\
   & = & i_{1,n+1}(\widetilde{e}_1)(h^1_{(1)},\dots, h^{n+1}_{(1)})E^{n-1}(f)\circ\mu_1(h^1_{(2)},\dots,h^{n+1}_{(2)})\\
\end{eqnarray*}

\begin{eqnarray*}
   & = & \widetilde{e}_1(h^1_{(1)})\varepsilon(h^2_{(1)})\dots\varepsilon( h^{n+1}_{(1)})E^{n-1}(f)(h^1_{(2)}h^2_{(2)},h^3_{(2)},\dots,h^{n+1}_{(2)})\\
   & = & (h^1_{(1)}\cdot 1_A)\varepsilon(h^2_{(1)})\dots\varepsilon( h^{n+1}_{(1)})(h^1_{(2)}h^2_{(2)}\cdot(f(h^3_{(2)},\dots,h^{n+1}_{(2)})))\\
   & = & (h^1_{(1)}\cdot 1_A)\varepsilon(h^2_{(1)})\dots\varepsilon( h^{n+1}_{(1)})(h^1_{(2)}h^2_{(2)}\cdot(f(h^3_{(2)},\dots,h^{n+1}_{(2)})))\\
   & = & (h^1_{(1)}\cdot 1_A)(h^1_{(2)}h^2\cdot(f(h^3,\dots,h^{n+1})))\\
   & = & h^1\cdot(h^2\cdot(f(h^3,\dots,h^{n+1})))\\
   & = & h^1\cdot(E^{n-1}(f)(h^2, h^3, \ldots, h^{n+1}))\\
   & = & E^n(E^{n-1}(f))(h^1,\dots,h^{n+1})\\
   & = & E^n\circ E^{n-1}(f)(h^1,\dots,h^{n+1}) 
    \end{eqnarray*}

\rmr{vi} $i_{n,n+1}(f\circ \mu_i)*i_{n-1,n}(f^{-1})\circ\mu_i = i_{n,n+1}(e_{n-1}\circ\mu_i)$, $\forall\ i\in\{1,\ldots, n-1\}$.
\vspace{0.3cm}

For $i\in \{1,\ldots, n-1\}$, we have
\begin{eqnarray*} 
& & i_{n,n+1}(f\circ \mu_i)*i_{n-1,n}(f^{-1})\circ\mu_i(h^1,\ldots,h^{n+1})= \\
& = &  i_{n,n+1}(f\circ \mu_i)(h^1_{(1)},\ldots, h^{n+1}_{(1)})i_{n-1,n}(f^{-1})\circ\mu_i(h^1_{(2)},\ldots, h^{n+1}_{(2)})\\
& = & f\circ \mu_i (h^1_{(1)},\ldots ,h^n_{(1)})\epsilon( h^{n+1}_{(1)}) 
 i_{n-1,n} (f^{-1})(h^1_{(2)},\ldots, h^i_{(2)}h^{i+1}_{(2)},\ldots, \underset{n-th\ coord}{\underbrace{h^{n+1}_{(2)}}})\\
& = & f(h^1_{(1)},\ldots, h^i_{(1)}h^{i+1}_{(1)},\ldots,h^n_{(1)})\epsilon( h^{n+1}_{(1)}) 
 f^{-1}(h^1_{(2)},\ldots, h^i_{(2)}h^{i+1}_{(2)},\ldots,h^n_{(2)})\epsilon(h^{n+1}_{(2)})\\
& = & f\circ\mu_i(h^1_{(1)},\ldots, h^i_{(1)},h^{i+1}_{(1)},\ldots,h^n_{(1)}) 
 f^{-1}\circ\mu_i(h^1_{(2)},\ldots, h^i_{(2)},h^{i+1}_{(2)},\ldots,h^n_{(2)})\epsilon (h^{n+1})\\
& = & (f\circ \mu_i)*(f^{-1}\circ\mu_i)(h^1,\ldots,h^n)\epsilon(h^{n+1})\\
& = & (f*f^{-1})\circ\mu_i(h^1,\dots,h^n)\epsilon(h^{n+1})\\
& = & e_{n-1}\circ\mu_i(h^1,\dots, h^n)\epsilon(h^{n+1}) \\
& = & i_{n,n+1}(e_n\circ\mu_i)(h^1,\ldots,h^{n+1}).
\end{eqnarray*}

\rmr{vii} $(i_{n-1,n}(f)\circ \mu_n)*i_{n-1,n+1}(f^{-1}) = i_{n-1,n+1}(e_{n-1})$.
\begin{eqnarray*} & & (i_{n-1,n}(f)\circ \mu_n)*i_{n-1,n+1}(f^{-1})(h^1,\dots,h^{n+1})= \\
   & = & (i_{n-1,n}(f)\circ \mu_n)(h^1_{(1)},\dots, h^{n+1}_{(1)})i_{n-1,n+1}(f^{-1})(h^1_{(2)},\dots, h^{n+1}_{(2)})\\
   & = & i_{n-1,n}(f)(h^1_{(1)},\dots,h^n_{(1)} h^{n+1}_{(1)})f^{-1}(h^1_{(2)},\dots, h^{n-1}_{(2)})\varepsilon(h^{n}_{(2)})\varepsilon(h^{n+1}_{(2)})\\
   & = & f(h^1_{(1)},\dots, h^{n-1}_{(1)})\varepsilon(h^n_{(1)} h^{n+1}_{(1)})f^{-1}(h^1_{(2)},\dots, h^{n-1}_{(2)})\varepsilon(h^{n}_{(2)})\varepsilon(h^{n+1}_{(2)})\\
   & = & f(h^1_{(1)},\dots, h^{n-1}_{(1)})f^{-1}(h^1_{(2)},\dots, h^{n-1}_{(2)})\varepsilon(h^n_{(1)} h^{n+1}_{(1)})\varepsilon(h^{n}_{(2)})\varepsilon(h^{n+1}_{(2)})\\
   & = & (f*f^{-1})(h^1,\dots, h^{n-1})\varepsilon(h^{n})\varepsilon(h^{n+1})\\
   & = & (e_{n-1})(h^1,\dots, h^{n-1})\varepsilon(h^{n})\varepsilon(h^{n+1})\\
   & = & i_{n-1,n+1}(e_{n-1})(h^1,\dots,h^{n+1})
\end{eqnarray*}

\rmr{viii} $(f\circ\mu_i\circ\mu_i)*(f^{-1}\circ\mu_i\circ\mu_{i+1}) = e_{n-1}\circ\mu_i\circ\mu_{i}$, $\forall\ i\in\{1,\ldots, n-1\}$.
\vspace{0.3cm}

For $i\in \{1,\ldots, n-1\}$, we have

\begin{eqnarray*} & &\hspace{-0.35cm} (f\circ\mu_i\circ\mu_i)*(f^{-1}\circ\mu_i\circ\mu_{i+1})(h^1,\dots,h^{n+1})=\\
   & = &\hspace{-0.35cm} f\circ\mu_i\circ\mu_i(h^1_{(1)},\dots,h^{n+1}_{(1)})f^{-1}\circ\mu_i\circ\mu_{i+1}(h^1_{(2)},\dots,h^{n+1}_{(2)})\\
   & = &\hspace{-0.35cm} f\hspace{-0.03cm}\circ\hspace{-0.03cm}\mu_i\hspace{-0.03cm}(h^1_{(1)},\dots, h^i_{(1)}h^{i+1}_{(1)}, h^{i+2}_{(1)} ,\dots,h^{n+1}_{(1)})f^{-1}\hspace{-0.03cm}\circ\hspace{-0.03cm}\mu_i(h^1_{(2)},\dots,h^i_{(2)}, h^{i+1}_{(2)}h^{i+2}_{(2)}\hspace{-0.03cm},\dots,h^{n+1}_{(2)})\\
   & = &\hspace{-0.35cm} f(h^1_{(1)},\dots, h^i_{(1)}h^{i+1}_{(1)} h^{i+2}_{(1)} ,\dots,h^{n+1}_{(1)})f^{-1}(h^1_{(2)},\dots,h^i_{(2)} h^{i+1}_{(2)}h^{i+2}_{(2)},\dots,h^{n+1}_{(2)})\\
   & = &\hspace{-0.35cm} (f*f^{-1})(h^1,\dots, h^ih^{i+1} h^{i+2} ,\dots,h^{n+1})\\
   & = &\hspace{-0.35cm} e_{n-1}(h^1,\dots, h^ih^{i+1} h^{i+2} ,\dots,h^{n+1})\\
   & = &\hspace{-0.35cm} e_{n-1}\circ\mu_i\circ\mu_i(h^1,\dots,h^{n+1})
  \end{eqnarray*}

\rmr{ix} $(f\circ\mu_i\circ\mu_{i+j})*(f^{-1}\circ\mu_{i+j-1}\circ\mu_i) = e_{n-1}\circ\mu_i\circ\mu_{i+j}$, $\forall\ i\in\{1,\ldots, n-2\}$, $\forall\ j \in \{2,\ldots, n-i\}$.
\vspace{0.3cm}

For $ i\in\{1,\ldots, n-2\}$ e $ j \in \{2,\ldots, n-i\}$, we have
\begin{eqnarray*}
& & (f\circ\mu_i\circ\mu_{i+j})*(f^{-1}\circ\mu_{i+j-1}\circ\mu_{i})(h^1,\ldots,h^{n+1})=\\
& = & f\circ\mu_i\circ\mu_{i+j}(h^1_{(1)},\ldots,h^{n+1}_{(1)})f^{-1}\circ\mu_{i+j-1}\circ\mu_{i}(h^1_{(2)},\ldots ,h^{n+1}_{(2)})\\
& = & f\circ\mu_i(h^1_{(1)},\ldots, h^{i+j}_{(1)}h^{i+j+1}_{(1)},\ldots,h^{n+1}_{(1)})f^{-1}\circ\mu_{i+j-1}(h^1_{(2)},
 \ldots,h^i_{(2)} h^{i+1}_{(2)},\ldots,h^{n+1}_{(2)})\\
& = & f(h^1_{(1)},\ldots, h^i_{(1)}h^{i+1}_{(1)},\ldots, h^{i+j}_{(1)}h^{i+j+1}_{(1)},\ldots,h^{n+1}_{(1)})\\
& & \; f^{-1}\circ\mu_{i+j-1}(h^1_{(2)},\ldots,\underset{i-th\ coord}{\underbrace{h^i_{(2)} h^{i+1}_{(2)}}},h^{i+2}_{(2)},\ldots, \underset{(i+j-1)-th\ coord}{\underbrace{h^{i+j}_{(2)}}},\ldots,h^{n+1}_{(2)})\\
& = & f(h^1_{(1)},\ldots, h^i_{(1)}h^{i+1}_{(1)},\ldots, h^{i+j}_{(1)}h^{i+j+1}_{(1)} ,\ldots,h^{n+1}_{(1)}) \\
& & \; f^{-1}(h^1_{(2)},\ldots,h^i_{(2)} h^{i+1}_{(2)},\ldots,h^{i+j}_{(2)}h^{i+j+1}_{(2)},\ldots,h^{n+1}_{(2)})\\
& = &(f*f^{-1})(h^1,\ldots, h^ih^{i+1},\ldots, h^{i+j}h^{i+j+1} ,\ldots,h^{n+1})\\
& = & e_{n-1}(h^1,\ldots, h^ih^{i+1},\ldots, h^{i+j}h^{i+j+1} ,\ldots,h^{n+1})\\
& = & e_{n-1}\circ\mu_i\circ\mu_{i+j}(h^1,\ldots,h^{n+1}) . 
  \end{eqnarray*} \find

\noindent\prove{Theorem \ref{deltaisnilpotent}} For any  $f\in C^{n}_{par}(H,A)$, 
we have that $\delta_{n+1}\circ\delta_n (f) = e_{n+2}$.

Indeed, take any $f\in C^{n}_{par}(H,A)$, then 
\begin{eqnarray*}
&\hspace{-0.3cm} &\hspace{-0.35cm} \delta_{n+1}(\delta_n(f)) = 
 E^{n+1}(\delta_n(f))* \prodt{i=1}{n+1} (\delta_n (f) )^{(-1)^{i}}\circ \mu_i *i_{n+1,n+2}((\delta_n (f))^{(-1)^{n+2}})\\
&\hspace{-0.3cm} =& \hspace{-0.35cm} E^{n+1}(\delta_n(f))* \prodt{i=1}{n+1} \delta_n(f^{(-1)^{i}})\circ \mu_i *i_{n+1,n+2}(\delta_n(f^{(-1)^{n+2}}))\\
&\hspace{-0.3cm} =& \hspace{-0.35cm} E^{n+1}\hspace{-0.03cm}(E^n(f)*\hspace{-0.05cm} \prodt{j=1}{n}\hspace{-0.05cm} f^{(-1)^{j}}\hspace{-0.03cm}\circ\hspace{-0.03cm} \mu_j * i_{n,n+1}(f^{(-1)^{n+1}}))
  *\prodt{i=1}{n+1}\hspace{-0.05cm} (E^n(f^{(-1)^{i}}) *\hspace{-0.03cm} \prodt{j=1}{n}\hspace{-0.05cm} f^{(-1)^{i+j}}\circ \mu_j\\
&\hspace{-0.3cm} &\hspace{-0.35cm} *i_{n,n+1}\hspace{-0.03cm}(f^{(-1)^{n+i+1}}\hspace{-0.03cm})\hspace{-0.03cm})\hspace{-0.04cm}\circ\hspace{-0.04cm} \mu_i\hspace{-0.035cm}  
 *\hspace{-0.03cm} i_{n+1,n+2}(E^n(f^{(-1)^{n+2}})*\hspace{-0.15cm} \prodt{j=1}{n}\hspace{-0.15cm} f^{(-1)^{n+j+2}}\hspace{-0.2cm}\circ\hspace{-0.05cm} \mu_j\hspace{-0.05cm} *\hspace{-0.02cm} i_{n,n+1}(\hspace{-0.03cm}f^{(-1)^{2n+3}})\hspace{-0.035cm})\\
\end{eqnarray*}
\begin{eqnarray*}
& = &\hspace{-0.3cm} \underset{Lemma\ \ref{p2}\ (v)}{\underbrace{E^{n+1}(E^n(f))}}*\hspace{-0.15cm} \prodt{j=1}{n}\hspace{-0.15cm} \underset{Lemma\ \ref{p2}\ (iv)}{\underbrace{E^{n+1}(f^{(-1)^{j}}\circ \mu_j)}} * \underset{Lemma\ \ref{p2}\ (i)}{\underbrace{E^{n+1}(i_{n,n+1}(f^{(-1)^{n+1}}))}}
 *\hspace{-0.15cm}\prodt{i=1}{n+1}\hspace{-0.15cm} E^n(f^{(-1)^{i}})\circ\mu_i\\ 
& & \hspace{-0.3cm} * \prodt{i=1}{n+1}\prodt{j=1}{n} f^{(-1)^{i+j}}\circ \mu_j\circ\mu_i 
 *\prodt{i=1}{n+1}i_{n,n+1}(f^{(-1)^{n+i+1}})\circ \mu_i  * i_{n+1,n+2}(E^n(f^{(-1)^{n+2}})) \\
& & \hspace{-0.3cm} *i_{n+1,n+2}(\prodt{j=1}{n} f^{(-1)^{n+j+2}}\circ \mu_j) * i_{n+1,n+2}(i_{n,n+1}(f^{(-1)^{2n+3}}))\\
& = & \hspace{-0.3cm} i_{1,n+2}(\wide{e}_1)*E^n(f)\circ\mu_1*\prodt{j=1}{n}E^n(f^{(-1)^j})\circ\mu_{j+1}*i_{n+1,n+2}(E^n(f^{(-1)^{n+1}}))\\
& &  \hspace{-0.3cm}  *\prodt{i=1}{n+1} E^n(f^{(-1)^{i}})\circ\mu_i * \prodt{i=1}{n+1}\prodt{j=1}{n} f^{(-1)^{i+j}}\circ \mu_j\circ\mu_i 
 *\prodt{i=1}{n+1}i_{n,n+1}(f^{(-1)^{n+i+1}})\circ \mu_i\\
& & \hspace{-0.3cm} * i_{n+1,n+2}(E^n(f^{(-1)^{n+2}})
 *\prodt{j=1}{n}i_{n+1,n+2}( f^{(-1)^{n+j+2}}\circ \mu_j) * i_{n,n+2}(f^{(-1)^{2n+3}}))\\
& = & \hspace{-0.3cm} i_{1,n+2}(\wide{e}_1)*\prodt{j=1}{n+1}E^n(f^{(-1)^{j+1}})\circ\mu_{j}*\prodt{i=1}{n+1} E^n(f^{(-1)^{i}})\circ\mu_i*
 *i_{n+1,n+2}(E^n(f^{(-1)^{n+1}}))\\
& & \hspace{-0.3cm} * i_{n+1,n+2}(E^n(f^{(-1)^{n+2}})
 * \prodt{i=1}{n+1}\prodt{j=1}{n} f^{(-1)^{i+j}}\circ \mu_j\circ\mu_i*\prodt{j=1}{n}i_{n+1,n+2}( f^{(-1)^{n+j+2}}\circ \mu_j )\\
& & \hspace{-0.3cm} \;  *\prodt{i=1}{n}i_{n,n+1}(f^{(-1)^{n+i+1}})\circ \mu_i *i_{n,n+1}(f^{(-1)^{2n+2}})\circ \mu_{n+1}
 * i_{n,n+2}(f^{(-1)^{2n+3}}))\\
& = & \hspace{-0.3cm} i_{1,n+2}(\wide{e}_1)*\prodt{j=1}{n+1}E^n(f^{(-1)^{j+1}}*f^{(-1)^{j}})\circ\mu_{j}
 *i_{n+1,n+2}(E^n(f^{(-1)^{n+1}}*f^{(-1)^{n+2}}))\\
& & \hspace{-0.3cm} * \hspace{-0.1cm} \prodt{i=1}{n+1}\prodt{j=1}{n} \hspace{-0.1cm} f^{(-1)^{i+j}} \hspace{-0.1cm}\circ \mu_j\circ\mu_i
 * \hspace{-0.1cm} \underset{Lemma\ \ref{p2}\ (vi)}{\underbrace{\prodt{j=1}{n} \hspace{-0.1cm}i_{n+1,n+2}( f^{(-1)^{n+j+2}}\circ \mu_j)* \hspace{-0.1cm}\prodt{i=1}{n} \hspace{-0.1cm}i_{n,n+1}(f^{(-1)^{n+i+1}})\circ \mu_i}} \\
& & \hspace{-0.3cm} *\underset{Lemma\ \ref{p2}\ (vii)}{\underbrace{i_{n,n+1}(f^{(-1)^{2n+2}})\circ \mu_{n+1}* i_{n,n+2}(f^{(-1)^{2n+3}}))}}\\
& = & \hspace{-0.3cm} i_{1,n+2}(\wide{e}_1)*\prodt{j=1}{n+1}\underset{Lema\ \ref{p1}\ (i)}{\underbrace{E^n(e_n)}}\circ\mu_{j}*i_{n+1,n+2}(\underset{Lemma\ \ref{p1}\ (i)}{\underbrace{E^n(e_n)}})
 * \prodt{i=1}{n+1}\prodt{j=1}{n} f^{(-1)^{i+j}}\circ \mu_j\circ\mu_i\\
 & & \hspace{-0.3cm} * \prodt{j=1}{n}i_{n+1,n+2}( e_n \circ \mu_j) * i_{n,n+2}(e_n)\\
& = &\hspace{-0.3cm} i_{1,n+2}(\wide{e}_1)*\prodt{j=1}{n+1}e_{n+1}\circ\mu_{j}*i_{n+1,n+2}(e_{n+1}) * \prodt{i=1}{n+1}\prodt{j=1}{n} f^{(-1)^{i+j}}\circ \mu_j\circ\mu_i \\
& &\hspace{-0.3cm}  * \prodt{j=1}{n}i_{n+1,n+2}( e_n \circ \mu_j) * i_{n+1,n+2}(i_{n,n+1}(e_n))\\
\end{eqnarray*}
\begin{eqnarray*}
& = &\hspace{-0.3cm} i_{1,n+2}(\wide{e}_1)*\prodt{j=1}{n+1}e_{n+1}\circ\mu_{j}*i_{n+1,n+2}(e_{n+1}) * \prodt{i=1}{n+1}\prodt{j=1}{n} f^{(-1)^{i+j}}\circ \mu_j\circ\mu_i \\
& &\hspace{-0.3cm} * \prodt{j=1}{n-1}i_{n+1,n+2}( e_n \circ \mu_j)*\underset{Lemma\ \ref{p1}\ (vi)}{\underbrace{i_{n+1,n+2}(e_n\circ\mu_n) * i_{n+1,n+2}(i_{n,n+1}(e_n))}}\\
& = &\hspace{-0.3cm} i_{1,n+2}(\wide{e}_1)*\prodt{j=1}{n+1}e_{n+1}\circ\mu_{j}*i_{n+1,n+2}(e_{n+1}) * \prodt{i=1}{n+1}\prodt{j=1}{n} f^{(-1)^{i+j}}\circ \mu_j\circ\mu_i \\ 
& &\hspace{-0.3cm}  *\underset{Lemma\ \ref{p1}\ (vii)}{\underbrace{\prodt{j=1}{n-1}i_{n+1,n+2}( e_n \circ \mu_j)*i_{n+1,n+2}(e_{n+1})}}\\
& = &\hspace{-0.3cm} i_{1,n+2}(\wide{e}_1)*\hspace{-0.1cm}\prodt{j=1}{n+1}\hspace{-0.1cm} e_{n+1}\circ\mu_{j}*i_{n+1,n+2}(e_{n+1}) *\hspace{-0.1cm} \prodt{i=1}{n+1}\prodt{j=1}{n}\hspace{-0.1cm} f^{(-1)^{i+j}}\hspace{-0.1cm}\circ \mu_j\circ\mu_i
 *i_{n+1,n+2}(e_{n+1})\\
& = &\hspace{-0.3cm} i_{1,n+2}(\wide{e}_1)*\prodt{j=1}{n}e_{n+1}\circ\mu_{j}*\underset{Lemma\ \ref{p1}\ (vi)}{\underbrace{e_{n+1}\circ\mu_{n+1}*i_{n+1,n+2}(e_{n+1})}}
 * \prodt{i=1}{n+1}\prodt{j=1}{n} f^{(-1)^{i+j}}\circ \mu_j\circ\mu_i\\
& = &\hspace{-0.3cm} i_{1,n+2}(\wide{e}_1)*\underset{Lemma\ \ref{p1}\ (vii)}{\underbrace{\prodt{j=1}{n}e_{n+1}\circ\mu_{j}*e_{n+2}}} * \prodt{i=1}{n+1}\prodt{j=1}{n} f^{(-1)^{i+j}}\circ \mu_j\circ\mu_i\\
& = &\hspace{-0.3cm} i_{1,n+2}(\wide{e}_1)*e_{n+2} * \prodt{i=1}{n+1}\prodt{j=1}{n} f^{(-1)^{i+j}}\circ \mu_j\circ\mu_i\\
& = &\hspace{-0.3cm} e_{n+2} * \prodt{i=1}{n+1}\prodt{j=1}{n} f^{(-1)^{i+j}}\circ \mu_j\circ\mu_i \\
& = &\hspace{-0.3cm} e_{n+2} * \underset{Lemma\ \ref{p2}\ (viii)}{\underbrace{\prodt{i=1}{n} f^{(-1)^{2i}}\circ \mu_i\circ\mu_i*\prodt{i=1}{n} f^{(-1)^{2i+1}}\circ \mu_i\circ\mu_{i+1}}}*\\
& &\hspace{-0.3cm} \underset{Lemma\ \ref{p2}\ (ix)}{\underbrace{*\prodt{i=1}{n-1}\prodt{j=2}{n-i+1} f^{(-1)^{2i+j}}\circ \mu_i\circ\mu_{i+j}  *\prodt{i=1}{n-1}\prodt{j=2}{n-i+1} f^{(-1)^{2i+j-1}}\circ \mu_{i+j-1}\circ\mu_i}}\\
& = &\hspace{-0.3cm} e_{n+2}* \prodt{i=1}{n} e_n\circ \mu_i\circ\mu_i*\prodt{i=1}{n-1}\prodt{j=2}{n-i+1} e_n\circ \mu_i\circ\mu_{i+j}\\
& = &\hspace{-0.3cm} e_{n+2}
\end{eqnarray*}

\vspace{-0.3cm}
Once $e_{n+2}$ absorbs $e_n\circ \mu_i\circ\mu_i$ for all $i\in \{1,\ldots,n\}$ and $e_n\circ \mu_i\circ\mu_{i+j}$, for all
$i\in \{1,\ldots,n-1\}$, $j\in\{2,\ldots,n-i-1\}$, by Lemma \ref{p2} (ii) and (iii) respectively, we conclude that $\delta_{n+1}\circ\delta_n(f) = e_{n+2}$ as we wanted to proof.
\find


\begin{thebibliography}{50}

\bibitem{Ade} A. Adem, R. Milgram: ``Cohomology of finite groups'', Springer Verlag (1994).

\bibitem{AB} M.M.S. Alves, E. Batista: ``Enveloping actions for partial Hopf actions'', Comm. Algebra 38 (2010), 2872-2902.

\bibitem{AB2} M. M. S. Alves,   E. Batista: ``Partial Hopf actions, partial invariants and a Morita context'', J. Algebra Discrete Math. 3 (2009), 1-19.

\bibitem{AB3}  M. M. S. Alves,   E. Batista:  ``Globalization theorems for partial Hopf (co)actions and some of their applications'', Contemp. Math.  537  (2011), 13-30.

\bibitem{ABDP1} M.M.S. Alves, E. Batista, M. Dokuchaev, A. Paques: ``Twisted partial actions of Hopf Algebras'', 
Israel Journal of Mathematics, 197 (2013) 263-308

\bibitem{ABDP2} M.M.S. Alves, E. Batista, M. Dokuchaev, A. Paques: ``Globalization of Twisted partial
  Hopf actions'', Journal of the Australian Mathematical Society 101 (2016) 1-28.
  
\bibitem{ABV} M.M.S. Alves, E. Batista, J. Vercruysse: ``Partial representations of Hopf algebras'', Journal of algebra 426 (2015) 137-187.  

\bibitem{Bat} E. Batista: ``Partial actions: what they are and why we care'', Bulletin of the Belgian Mathematical Society Simon Stevin 24 (2017) 35-71.
  
\bibitem{BV} E. Batista, J. Vercruysse: ``Dual constructions for partial actions of Hopf algebras'', J. Pure Appl. Algebra 220.2 (2016) 518-559.  

\bibitem{B} G. Böhm: ``Hopf algebroids'', Handbook of Algebra Vol 6, ed. by M. Hazewinkel, Elsevier (2009) 173-236. arXiv:0805.3806

\bibitem{BB} G. Böhm, T. Brzezinski: ``Cleft extensions of Hopf algebroids'', Appl. Categorical Structures 14.5-6 (2006) 431-469.
  
\bibitem{BS1} G. Böhm, D. Stefan: ``(Co) cyclic (co) homology of bialgebroids: an approach via (co) monads'', Communications in Mathematical Physics 282.1 (2008) 239-286.


\bibitem{CJ} S. Caenepeel, K. Janssen: ``Partial (co)actions of Hopf algebras and partial Hopf-Galois theory'', 
Communications in Algebra 36, (2008) 2923-2946.


\bibitem{DE} M. Dokuchaev, R. Exel:``Associativity of Crossed Products by Partial Actions, Enveloping Actions and Partial Representations'', {\sl Trans. Amer. Math. Soc.} {\bf 357} (5) (2005) 1931-1952.

\bibitem{DES1} M. Dokuchaev, R. Exel, J. J. Simon: ``Crossed products by twisted partial actions and graded algebras'',  J. Algebra 320 No. 8 (2008) 3278-3310.

\bibitem{DES2} M. Dokuchaev, R. Exel, J. J. Simon: ``Globalization of  twisted partial actions'', Trans. Amer. Math. Soc. 362 No. 8 (2010) 4137-4160.

\bibitem{DK1} M. Dokuchaev, M. Khrypchenko: ``Partial cohomology of groups'', Journal of Algebra 427 (2015) 142-182.

\bibitem{DK2} M. Dokuchaev, M. Khrypchenko: ``Twisted partial actions and extensions of 
semilattices of groups by groups'' preprint arXiv:1602.02424 (2016).


\bibitem{Ex1} R. Exel: ``Circle actions on $C^{*}$-algebras, partial automorphisms and generalized Pimsner-Voiculescu 
exact sequences'', J. Funct. Anal. 122.3 (1994) 361-401.

\bibitem{Ex2} R. Exel: ``Partial actions of groups and actions of inverse semigroups'',
Proc. Amer. Math. Soc. 126.12 (1998) 3481-3494.

\bibitem{Hopf} H. Hopf: ``\"{U}ber die topologie der Gruppen-Mannigfaltigkeiten und ihrer verallgemeinerungen'', in Collected Papers - Gesammelte Abhandlungen (Springer Collected Works in Mathematics), Springer-Verlag (2013).

\bibitem{KoP} N. Kowalzig, H. Posthuma: ``The cyclic theory of Hopf algebroids'', Journal of Noncommutative Geometry 5.3 (2011) 423-476.


\bibitem{Lau} H. Lausch: ``Cohomology of inverse semigroups'', J. Algebra 35 (1975), 273-303.

\bibitem{Law} M. Lawson: ''Inverse Semigroups'', World Scientific (1998).

\bibitem{Sch} P. Schauenburg: ``Cohomological obstructions to cleft extensions over cocommutative Hopf algebras'',
K-Theory. 24 (2001) 227-242.

\bibitem{Swe} M.E. Sweedler: ``Cohomology of algebras over Hopf algebras'',  Transactions of the 
American Mathematical Society 133.1 (1968) 205-239.

\bibitem{W} E. Weiss: ``Cohomology of groups'', Academic Press (1969). 

\end{thebibliography}
\end{document}